\newif\ifreport
\reporttrue

\newif\ifhighlighted
\highlightedtrue

\ifreport
    \documentclass{article}\def\zibreport{1}

    \usepackage{authblk}
    \makeatletter
    \renewcommand\AB@authnote[1]{\rlap{\textsuperscript{\normalfont#1}}}

    \makeatother

    \usepackage{natbib}
    \bibpunct[, ]{(}{)}{,}{a}{}{,}%
\else
    \documentclass[ijoc,nonblindrev]{informs3}\def\zibreport{0}

    \OneAndAHalfSpacedXII %
    \usepackage{natbib}
    \bibpunct[, ]{(}{)}{,}{a}{}{,}%
    \TheoremsNumberedThrough     %
    \EquationsNumberedThrough    %
\fi

\usepackage[utf8x]{inputenc}
\usepackage[english]{babel}

\usepackage{pdflscape}

\usepackage{tikz}
\usetikzlibrary{decorations.pathmorphing}

\usepackage{pgfplots, pgfplotstable}
\pgfkeys{/pgf/number format/.cd,1000 sep={}}
\usepgfplotslibrary{groupplots}
\pgfplotsset{every axis/.append style={label style={font=\footnotesize}, tick label style={font=\footnotesize}}}
\pgfplotsset{compat=1.8}

\ifreport
    \usepackage[font=small]{caption}
\fi

\usepackage{subcaption}
\usepackage[textsize=small]{todonotes}
\usepackage{enumitem}
\usepackage[linesnumbered,ruled,vlined]{algorithm2e}

\usepackage{siunitx}

\usepackage{natbib}

\newcommand{\thetitle}{Conflict-Driven Heuristics for Mixed Integer Programming}

\newcommand{\theabstract}{
  Two essential ingredients of modern mixed-integer programming (MIP) solvers are diving heuristics that simulate
  a partial depth-first search in a branch-and-bound search tree and conflict
  analysis of infeasible subproblems to learn valid constraints.
  So far, these techniques have mostly been studied independently:
  primal heuristics under the aspect of finding high-quality feasible solutions early during the
  solving process and conflict analysis for fathoming nodes of the search tree
  and improving the dual bound.
  Here, we combine both concepts in two different ways.
  First, we develop a diving heuristic that targets the generation of valid conflict constraints from the Farkas dual.
  We show that in the primal this is equivalent to the optimistic strategy of diving towards the best bound
  \wrt the objective function.
  Secondly, we use information derived from conflict analysis to enhance the search of a diving heuristic
  akin to classical coefficient diving.
  The computational performance of both methods is evaluated using an implementation in the source-open MIP solver \scip.
  Experiments are carried out on publicly available test sets including \miplib~2010 and \coral.
}

\clearpage{}%

\ifreport
 \usepackage{amsmath}
 \usepackage{amssymb}
 \usepackage{ifthen}
 \usepackage{xifthen}
 \usepackage{amsthm}
\fi

\usepackage{xspace}
\usepackage{xcolor}

\usepackage{multirow}

\usepackage{booktabs}
\usepackage{tabularx}

\usepackage{longtable}

\usepackage{microtype}

\ifreport
    \newtheorem{Def}{Definition}
    \newtheorem{Prop}[Def]{Proposition}
    \newtheorem{Lem}[Def]{Lemma}
    \newtheorem{Thm}[Def]{Theorem}
    \newtheorem{Obs}[Def]{Observation}
    \newtheorem{Ex}[Def]{Example}
\else
    \theoremstyle{THkey}
    \newtheorem{Def}{Definition}

    \newtheorem{Lem}{Lemma}

\fi

\newcommand{\ie}{i.e.,\xspace}
\newcommand{\eg}{e.g.,\xspace}

\newcommand{\wrt}{with respect to\xspace}

\newcolumntype{L}{>{\raggedright\arraybackslash}X}%
\newcolumntype{C}{>{\centering\arraybackslash}X}%
\newcolumntype{R}{>{\raggedleft\arraybackslash}X}%

\newcommand{\cleaninst}{all\xspace}
\newcommand{\affected}{affected\xspace}
\newcommand{\affectedsols}{afterroot\xspace}
\newcommand{\bracket}[2]{[#1,#2]\xspace}
\newcommand{\better}[1]{%
    \ifhighlighted%
        \textbf{\textcolor{blue}{#1}}%
    \else%
        {#1}%
    \fi%
}%

\newcommand{\quot}{\textbf{\textsubscript{\hspace*{-.3em}Q}}}
\newcommand{\virt}{\textbf{\textsubscript{\hspace*{-.3em}\text{virt}}}}

\newcommand{\solved}{\ifreport\textbf{S}\else\textbf{solved}\fi\xspace}

\renewcommand{\time}{\ifreport\textbf{T}\else\textbf{time}\fi\xspace}
\newcommand{\timeQ}{\time\quot\xspace}

\newcommand{\nodes}{\ifreport\textbf{N}\else\textbf{nodes}\fi\xspace}
\newcommand{\nodesQ}{\nodes\quot\xspace}

\newcommand{\sols}{\ifreport\textbf{F}\else\textbf{feas}\fi\xspace}

\newcommand{\virtualbestdsols}{\sols\virt\xspace}

\newcommand{\bestsols}{\ifreport\textbf{I}\else\textbf{impr}\fi\xspace}

\newcommand{\virtualbestdbsols}{\bestsols\virt\xspace}

\newcommand{\confs}{\ifreport\textbf{C}\else\textbf{confs}\fi\xspace}

\newcommand{\percent}[1]{\SI{#1}{\%}\xspace}

\newcommand{\testset}[1]{\textsc{#1}\xspace}
\newcommand{\miplib}{\testset{Miplib}}
\newcommand{\coral}{\testset{Cor@l}}
\newcommand{\MMMC}{\testset{MMMC}}

\newcommand{\farkasdiving}{Farkas diving\xspace}
\newcommand{\conflictdiving}{conflict diving\xspace}
\newcommand{\coefficientdiving}{coefficient diving\xspace}

\newcommand{\setting}[1]{\texttt{#1}}

\newcommand{\default}{\setting{default}\xspace}

\newcommand{\confdiving}{\texttt{conf\-diving}\xspace}
\newcommand{\confdivingnoconfs}{\texttt{conf\-diving\discretionary{-}{}{-}noconfs}\xspace}
\newcommand{\confdivingnosols}{\texttt{conf\-diving\discretionary{-}{}{-}nosols}\xspace}

\newcommand{\confdivinglikecoef}{\texttt{conf\discretionary{-}{}{-}like\discretionary{-}{}{-}coefdiving}\xspace}

\newcommand{\farkdiving}{\texttt{fark\-diving}\xspace}
\newcommand{\farkdivingnoconfs}{\texttt{fark\-diving\discretionary{-}{}{-}noconfs}\xspace}
\newcommand{\farkdivingnosols}{\texttt{fark\-diving\discretionary{-}{}{-}nosols}\xspace}
\newcommand{\farkdivinglp}{\texttt{fark\-diving\discretionary{-}{}{-}lp}\xspace}

\newcommand{\coefdiving}{\texttt{coef\-diving}\xspace}

\newcommand{\nolockdiving}{\texttt{no\-lock\-diving}\xspace}

\newcommand{\scip}{\texttt{SCIP}\xspace}
\newcommand{\scipv}{\scip~6.0\xspace}
\newcommand{\soplex}{\texttt{SoPlex}\xspace}
\newcommand{\soplexv}{\soplex~4.0\xspace}

\newcommand{\bandb}{branch-and-bound\xspace}
\newcommand{\Bandb}{Branch-and-bound\xspace}

\newcommand{\Abs}[1]{\vert #1\vert}

\newcommand{\intvars}{\mathcal{I}}
\newcommand{\allvars}{\mathcal{N}}

\newcommand{\opt}[1]{#1^{\star}}

\newcommand{\floor}[1]{\lfloor #1\rfloor}
\newcommand{\ceil}[1]{\lceil #1\rceil}

\newcommand{\UP}{\texttt{up}}
\newcommand{\DOWN}{\texttt{down}}

\newcommand{\divecands}{\mathcal{D}}
\newcommand{\roundfunc}[1]{\phi_{#1}}
\newcommand{\scorefunc}[1]{\psi_{#1}}

\newcommand{\uplock}[1]{\zeta^{-}_{#1}}
\newcommand{\downlock}[1]{\zeta^{+}_{#1}}
\newcommand{\confuplock}[1]{\xi^{-}_{#1}}
\newcommand{\confdownlock}[1]{\xi^{+}_{#1}}
\newcommand{\lpfrac}[1]{f_{#1}}
\newcommand{\rellpfrac}[1]{f^{rel}_{#1}}
\newcommand{\dualimpact}[1]{d_{#1}}

\newcommand{\confweight}{\kappa}
\newcommand{\upweight}[1]{\rho^{-}_{#1}}
\newcommand{\downweight}[1]{\rho^{+}_{#1}}

\newcommand{\lb}{\ell}
\newcommand{\ub}{u}
\newcommand{\loclb}{\lb^\prime}
\newcommand{\locub}{\ub^\prime}

\newcommand{\varindex}{\ensuremath{j}}
\newcommand{\rowindex}{\ensuremath{i}}

\newcommand{\rhsvector}{\ensuremath{b}\xspace}
\newcommand{\objvector}{\ensuremath{c}\xspace}

\newcommand{\coefmatrix}{\ensuremath{A}\xspace}

\newcommand{\var}{\ensuremath{x}\xspace}

\newcommand{\T}{\mathsf{T}}

\newcommand{\dualvector}{\ensuremath{y}\xspace}
\newcommand{\redcosts}{\ensuremath{r}\xspace}
\newcommand{\dualsol}{(\dualvector,\redcosts)}

\newcommand{\Z}{\mathbb{Z}\xspace}
\newcommand{\R}{\mathbb{R}\xspace}

\clearpage{}%
\ifthenelse{\zibreport = 1}{
  \usepackage{zibtitlepage}

  \ZTPTitle{\thetitle}
  \ZTPAuthor{
    \ZTPHasOrcid{Jakob Witzig}{0000-0003-2698-0767}
    \hspace{2cm}
    \ZTPHasOrcid{Ambros Gleixner}{0000-0003-0391-5903}
  }
  \ZTPPreprint
  \ZTPNumber{19-08}
  \ZTPMonth{February}
  \ZTPYear{2019}
}{}

\begin{document}

\ifthenelse{\zibreport = 0}{
    \RUNAUTHOR{Witzig and Gleixner}

    \RUNTITLE{\thetitle}

    \TITLE{\thetitle}

    \ARTICLEAUTHORS{%
    \AUTHOR{Jakob Witzig}
    \AFF{Zuse Institute Berlin, Takustr.~7, 14195~Berlin, Germany, \EMAIL{witzig@zib.de}, ORCID: \href{https://orcid.org/0000-0003-2698-0767}{0000-0003-2698-0767}}
    \AUTHOR{Ambros Gleixner}
    \AFF{Zuse Institute Berlin, Takustr.~7, 14195~Berlin, Germany, \EMAIL{gleixner@zib.de}, ORCID: \href{https://orcid.org/0000-0003-0391-5903}{0000-0003-0391-5903}}
    } %

    \ABSTRACT{%
        \theabstract%
    }%


    \maketitle
}%
{
  \author{Jakob~Witzig}
  \author{Ambros~Gleixner}

  \affil{Zuse Institute Berlin, Takustr.~7, 14195~Berlin, Germany \protect\\ \texttt{\{witzig,gleixner\}@zib.de}}

  \zibtitlepage

  \title{\thetitle}

  \maketitle

  \begin{abstract}
    \theabstract
  \end{abstract}
}

\section{Introduction}
\label{sec:introduction}

The most commonly used method to solve \emph{mixed-integer programs} (MIPs) is the
\emph{linear programming}-based (LP-based) \bandb algorithm~\citep{Dakin1965,LandDoig1960}.
In modern MIP solvers, this procedure is accelerated by various extensions~\citep[see \eg][]{bixby1999mip,laundy2009solving}.
Two examples of those extensions are primal heuristics~\citep[see \eg][]{fischetti2010heuristics,lodi2013heuristic,Berthold2014}
and conflict analysis~\citep[see \eg][]{DaveyBolandStuckey2002,SandholmShields2006,achterberg2007conflict,witzig2017experiments}.
A primal heuristic is an incomplete method without any guarantee of success, which is used to find feasible and improving solutions.
Computational studies indicate that within a MIP solver, disabling all primal heuristics would lead to a deterioration
of solving time by approximately $11\%$ -- $32\%$~\citep{Berthold2014} and $5\%$ -- $15\%$~\citep{AchterbergWunderling2013}.
Conflict analysis denotes a collection of techniques to learn from infeasible subproblems encountered during the MIP solve.
The outcome of conflict analysis is a set of so-called conflict constraints that are used in the remainder of the solving process, \eg for propagation.
As a consequence, the proof of global optimality can be accelerated, mostly by reducing the number of subproblems that
need to be explored, according to the study of~\cite{AchterbergWunderling2013} by as much as $28\%$ on affected instances.

In this paper, we propose two complementary ways of combining both concepts.
Firstly, we develop a primal diving heuristic that explicitly aims to generate conflict constraints.
As we show by elementary calculations, this amounts to an optimistic fixing of variables to their best bound \wrt the objective function.
Previous approaches that are targeted to gain additional conflict information starting from a feasible subproblem are based on,
for example, an involved random sampling approach~\citep{dickerson2013throwing} or use a black-box solver
to perform a hybrid constraint programming and MIP search~\citep{BertholdFeydyStuckey2010,BertholdStuckeyWitzig2018}.
In contrast to that, our approach constitutes a more direct method.

Secondly, we use the information obtained by conflict analysis in order to guide the LP relaxation towards feasibility.
To this end, we apply the concept of variable locks to conflict constraints
and show that this type of locks is richer on information and yield a more dynamic criterion
compared to variable locks as known from the literature~\citep{achterberg2007constraint}.
We use this to develop a new diving heuristic that harnesses variable locks implied by conflict constraints.
Our experiments indicate that this heuristic outperforms the well-known \coefficientdiving heuristic~\citep{Berthold2008}.

To show how a MIP solver can benefit from the techniques presented in this paper as supplementary features,
we carry out a detailed computational study for which both heuristics were implemented within the academic MIP solver \scip~6.0~\citep{GleixnerBastubbeEifleretal.2018}.
The heuristics presented in this paper are -- to the best of our knowledge -- the first LP-based heuristics for MIP which explicitly produce and exploit conflict constraints.\\

This paper is organized as follows.
In Section~\ref{sec:background} we give a brief overview of all the background we need in the remainder of this paper:
LP-based \bandb, diving heuristics, and conflict analysis.
In Section~\ref{sec:farkasdiving} we discuss how a diving heuristic can be used to generate additional conflict information explicitly.
Afterward, we present a modification and extension of the well-known diving heuristic \coefficientdiving by using conflict information in Section~\ref{sec:conflictdiving}.
Finally, an intense computational study of the individual impact of both presented approaches is presented in Section~\ref{sec:experiments}.
In Section~\ref{sec:conclusionandoutlook} we conclude.

\section{Background}
\label{sec:background}

We consider MIPs of the form
\begin{align}
  \min \{ \objvector^\T\var\;|\;\coefmatrix\var \geq \rhsvector,\ \lb \leq \var \leq \ub,\ \var_\varindex \in \Z\, \forall \varindex \in \intvars \}, \label{eq:mip}
\end{align}
with objective coefficient vector $\objvector \in \R^{n}$, constraint coefficient matrix $\coefmatrix \in \R^{m\times n}$,
constraint right-hand side $\rhsvector \in \R^{m}$, and variable bounds $\lb, \ub \in \overline{\R}^{n}$,
where $\overline{\R} := \R \cup \{\pm\infty\}$.
Moreover, let $\intvars \subseteq \allvars := \{1,\ldots,n\}$ be the index set of integer variables.

\paragraph{LP-based \Bandb.}

\Bandb~\citep{Dakin1965,LandDoig1960} is a divide-and-conquer method which splits the search space sequentially into smaller
subproblems that are ideally easier to solve.
For each subproblem a lower bound is computed.
To this end, the integrality requirements are omitted and the LP relaxation
\begin{align}
  \min \{ \objvector^\T\var\;|\;\coefmatrix\var \geq \rhsvector,\ \lb \leq \var \leq \ub,\ \var \in \R^{n} \} \label{eq:lprelax}
\end{align}
is solved.
On the other hand, an upper bound on the global problem is given by the objective value of the incumbent solution,
\ie the best solution found so far, if available.

During the \bandb procedure subproblems are regularly fathomed, either due to bounding or infeasibility.
In the first case, subproblems whose lower bound exceeds the global upper bound are disregarded because they cannot contain an improving solution.
Therefore, it is evident that \bandb algorithms benefit directly from finding good solutions as early as possible.
These solutions either originate directly from the LP relaxation when all variables fulfill the integrality conditions in the LP solution,
or are constructed by so-called primal heuristics~\citep{fischetti2010heuristics,lodi2013heuristic,Berthold2014}.
In the second case, the infeasibility of a subproblem is either proven by contradicting variable bound changes or by an infeasible LP relaxation.
If a node is fathomed due to infeasibility modern MIP solvers use conflict analysis~\citep{DaveyBolandStuckey2002,achterberg2007constraint,witzig2017experiments}
to ``learn'' from those subproblems.
Note that every subproblem fathomed due to bounding can be interpreted as an infeasible subproblem after adding a cutoff constraint on the objective function that restricts the feasible region to improving solutions.

\paragraph{Diving Heuristics.}

{
\begin{algorithm}[t]
    \DontPrintSemicolon
    \SetKwInOut{Input}{Input}
    \SetKwInOut{Output}{Output}

    \Input{LP solution $\var^{LP}$, rounding function $\roundfunc{}$, score function $\scorefunc{}$}
    \Output{Solution candidate $\hat{\var}$ or \texttt{NULL}}
    $\hat{\var} \gets$ \texttt{NULL}, $\tilde{\var} \gets \var^{LP}$\\
    $\divecands \gets \{ \varindex \in \intvars\,|\, \tilde{\var}_\varindex \notin \Z\}$ \tcp*{diving candidates}
    \While{ $\hat{\var} == $ \texttt{NULL} and $\divecands \neq \emptyset$ }
    {
      \ForAll{ $i \in \divecands$ }
      {
        \begin{enumerate}[leftmargin=.25cm]
          \item Determine rounding direction: $d_\varindex \gets \roundfunc{}(\varindex)$ \ifreport\vspace{-.5em}\fi
          \item Calculate variable score: $s_\varindex \gets \scorefunc{}(\varindex)$ \ifreport\vspace{-.5em}\fi
        \end{enumerate}
      } \label{line:direction-and-score}
      Select candidate $\var_\varindex$ with maximal score $s_\varindex$ and current local bounds $\lb_\varindex$ and $\ub_\varindex$ \label{line:select} \\
      Update $\divecands \gets \divecands \setminus \varindex$ \\
      \leIf{ $d_\varindex == $ \UP }
      {
        $\ell_\varindex \gets \ceil{\tilde{\var}_{\varindex}}$ %
      }
      {
        $u_\varindex \gets \floor{\tilde{\var}_{\varindex}}$ %
      }
      \texttt{(optional)} Propagate this bound change\label{line:propagation}\\
      Update $\divecands$ if propagation fixed some $\varindex \in \divecands$ \\
      \If{ Infeasibility detected }
      {
        Analyze infeasibility, add conflict constraints\label{line:conflict1}, perform 1-level backtrack \\
        If $\divecands = \emptyset$ \texttt{goto}~\ref{line:return} or~\ref{line:select} otherwise
      }
      \texttt{(optional)} Solve local LP relaxation \label{line:solveLP}\\
      \If{ Infeasibility detected }
      {
        Analyze infeasibility, add conflict constraints\label{line:conflict2}, perform 1-level backtrack \\
        If $\divecands = \emptyset$ \texttt{goto}~\ref{line:return} or~\ref{line:select} otherwise
      }
      Update $\tilde{\var}$ and $\divecands$ if LP was solved \\
      \If{ $\tilde{\var}_\varindex \in \Z$ for all $\varindex \in \intvars$ or $\divecands == \emptyset$ }{ $\hat{\var} \gets \tilde{\var}$ }
    }

    \Return $\hat{\var}$ \label{line:return}
    \caption{\textsc{GenericDivingProcedure}}\label{alg:divingalgo}
\end{algorithm}
}
 
A special type of primal heuristics are so-called \emph{diving heuristics}
such as fractional-diving and pseudo-cost diving~\citep{Berthold2008}.
The principle idea of diving heuristics comes from the \bandb procedure itself.
Starting from a feasible but fractional LP solution, diving heuristics alternate between fixing some integer variables to a rounded value based on a fractional LP solution
and reoptimizing the LP relaxation.
This procedure can be viewed as a partial tree search along one path from the current subproblem to a leaf.
Diving heuristics use a special branching rule that usually tends towards feasibility.
By contrast, branching rules for complete tree search such as reliability branching~\citep{achterberg2005branching} aim at a good subdivision of the problem.
In modern MIP solvers, the basic and simple idea of diving heuristics (see Algorithm~\ref{alg:divingalgo}) is extended by
constraint propagation (see Algorithm~\ref{alg:divingalgo} Line~\ref{line:propagation}) and
conflict analysis (see Algorithm~\ref{alg:divingalgo} Line~\ref{line:conflict1} and~\ref{line:conflict2}).

\paragraph{Conflict Analysis for Infeasible LP Relaxations.}

If the LP relaxation of a subproblem of~\eqref{eq:mip} with local bounds $\loclb, \locub$ is proven to be infeasible
there exists a dual ray $(\dualvector, s) \in \R^m_+ \times \R^n$ by the Lemma of Farkas~\citep{Farkas1902,Dinh2014} such that
\begin{align}
  \dualvector^\T\coefmatrix + s &= 0 \label{eq:redcost}\\
  \dualvector^\T\rhsvector + s\{\loclb,\locub\} &> 0, \label{eq:dualproof}
\end{align}
where we use the notation%
\ifreport
    $s\{\loclb,\locub\} := \sum_{j \in \allvars\colon s_\varindex > 0} s_\varindex\loclb_\varindex + \sum_{j \in \allvars\colon s_\varindex < 0} s_\varindex \locub_\varindex$%
\else
$s\{\loclb,\locub\} := \sum\limits_{j \in \allvars\colon s_\varindex > 0} s_\varindex\loclb_\varindex + \sum\limits_{j \in \allvars\colon s_\varindex < 0} s_\varindex \locub_\varindex$%
\fi.
Therefore, the inequality
\begin{align}\label{eq:farkasproof}
  (\dualvector^\T\coefmatrix)\var \geq \dualvector^\T\rhsvector
\end{align}
is globally valid.
We refer to~\eqref{eq:farkasproof} also as \emph{Farkas proof}.
\cite{polik15ismp} and \cite{witzig2017experiments} describe how these constraints can be collected, managed,
and used for constraint propagation to deduce tighter variable bounds in modern MIP solvers.

\section{Farkas Diving}
\label{sec:farkasdiving}

Our first aim is the design of a diving procedure such that the dual solution of the LP relaxation moves towards a valid Farkas proof,
\ie constraints~\eqref{eq:redcost} and~\eqref{eq:dualproof} are satisfied.
Suppose $\opt{\var}$ is an optimal but fractional solution of a local LP relaxation \wrt bounds $\loclb$ and $\locub$.
Let $(\opt{\dualvector}, \opt{\redcosts})$ be an optimal solution of the dual LP
\begin{align}
  \max \{ \dualvector^\T\rhsvector + \redcosts\{\loclb, \locub\}\;|\;\dualvector^\T\coefmatrix + \redcosts = \objvector, \dualvector \in \R^m_+, \redcosts \in \R^n \},
\end{align}
where $\redcosts_\varindex$ is the \emph{reduced cost} of $\var_\varindex$, for all $\varindex \in \allvars$.
The dual solution $(\opt{\dualvector}, \opt{\redcosts})$ is not feasible for~\eqref{eq:redcost}
and~\eqref{eq:dualproof} with $(\dualvector, s) = (\opt{\dualvector}, \opt{\redcosts})$,
but note that $(\dualvector, s) = (\opt{\dualvector}, \opt{\redcosts} - \objvector)$ fulfills~\eqref{eq:redcost}.

In order to reduce the violation of~\eqref{eq:dualproof}, we need to increase the lower bound $\loclb_\varindex$
of $\var_\varindex$ if $\opt{\redcosts}_\varindex - \objvector_\varindex > 0$
and decrease the upper bound $\locub_\varindex$ if $\opt{\redcosts}_\varindex - \objvector_\varindex < 0$.
By complementary slackness, all integer variables with fractional LP solution value have reduced costs zero, \ie $\opt{\redcosts}_\varindex = 0$.
Here, without loss of generality we assume $\loclb_\varindex, \locub_\varindex \in \Z$ for all $\varindex \in \intvars$.
Hence, tightening the lower or upper bound of variable $\varindex$ reduces the violation of~\eqref{eq:dualproof} by $\Abs{\objvector_\varindex}\cdot\dualimpact{\varindex}$, where
\begin{align}
  \dualimpact{\varindex} :=
  \begin{cases}
        \ceil{\opt{\var}_\varindex} - \loclb_\varindex      &\text{ if } \objvector_\varindex < 0,\\
        \locub_\varindex - \floor{\opt{\var}_\varindex}     &\text{ if } \objvector_\varindex > 0.\\
  \end{cases} \label{eq:dualimpact}
\end{align}
Therefore, the rounding direction of a variable $\varindex$ is implied by the sign of the objective coefficient $\objvector_\varindex$,
\ie upwards if $\objvector_\varindex < 0$, downwards if $\objvector_\varindex > 0$.

The previous part of this section took a strictly dual point of view.
When switching the perspective to the primal side, the above rounding procedure can be interpreted as follows.
On variables with negative objective coefficients we always tighten the lower bounds,
\ie the solution values of those variables are pushed towards the upper bound, and vice versa.
Since~\eqref{eq:mip} is a minimization problem, all variables are rounded into the best direction \wrt the objective function $c$.
If $\objvector_\varindex = 0$ neither pushing to the lower nor upper bound has a direct impact on the objective function.
In that case, a natural choice for breaking this tie is to consider the \emph{fractionality} $\lpfrac{\varindex} := \opt{\var}_\varindex - \floor{\opt{\var}_\varindex}$ of the corresponding LP solution value.
This leads to the following diving heuristic, defined by a rounding function $\roundfunc{F}$ and scoring function $\scorefunc{F}$.
For every variable $\varindex \in \intvars$, let
\begin{align}
  \roundfunc{F}(\varindex) :=
  \begin{cases}
        \UP    &\hspace{-.75em}\text{ if } \objvector_\varindex < 0 \text{ or } \objvector_\varindex = 0 \wedge \lpfrac{\varindex} \geq \frac{1}{2},\\
        \DOWN  &\hspace{-.75em}\text{ if } \objvector_\varindex > 0 \text{ or } \objvector_\varindex = 0 \wedge \lpfrac{\varindex} < \frac{1}{2}.\\
  \end{cases}
\end{align}
While we use the fractionality only for tie-breaking, other diving heuristics, \eg fractionality diving~\citep{achterberg2007constraint},
use this criterion solely to determining the rounding directions.

In addition to the rounding direction $\roundfunc{F}$, an order needs to be defined in which the diving candidates are explored.
As discussed before, the violation of the dual infeasibility constraint~\eqref{eq:dualproof} can be reduced by $\Abs{\objvector_\varindex} \cdot \dualimpact{\varindex}$
when tightening the lower or upper bound of variable~$\varindex$.
From the primal point of view, the potential change in the objective function depends on how much a variable can be pushed until it reaches one of its bounds.
This change can be approximated by $\Abs{\objvector_\varindex} \cdot \rellpfrac{\varindex}$, where
\begin{align}
        \rellpfrac{\varindex} :=
                \begin{cases}
                        1 - \lpfrac{\varindex} & \text{if } \roundfunc{F(\varindex)} = \UP, \\
                        \lpfrac{\varindex} & \text{if } \roundfunc{F}(\varindex) = \DOWN,
                \end{cases}
\end{align}
is the \emph{relative fractionality} of $\opt{\var}_\varindex$, for all $\varindex \in \intvars$.
This measure is used by, \eg pseudo cost diving~\citep{achterberg2007constraint}, to define an ordering in which the variables are explored during diving.

Combining both criteria, let the score of $\varindex \in \intvars$ be given by
\begin{align}
  \scorefunc{F}(\varindex) := \Abs{\objvector_\varindex} \cdot \dualimpact{\varindex} \cdot \rellpfrac{\varindex}.
\end{align}
A higher score is preferred.

In the following, we refer to this diving heuristic as \emph{\farkasdiving}.
Its rounding procedure pushes all variable towards the so-called \emph{pseudo solution}~\citep{achterberg2007constraint}.
In the pseudo solution, each variable takes the best bound \wrt its objective coefficient as solution value.
This relaxation solution is overly optimistic and most of the time not feasible for the constraints~$\coefmatrix\var \geq \rhsvector$.
However, if this strategy leads to a feasible solution, it may be expected to be very good.

\section{Conflict Diving}
\label{sec:conflictdiving}

We continue by exploring the complementary question of how conflict constraints can be used to guide the search of a diving heuristic.
First, consider the following well-known concept.

\ifreport
\begin{Def}[Variable Locks~\citep{achterberg2007constraint}]\label{def:variablelock}
\else
\smallskip
\begin{definition}[Variable Locks~\citep{achterberg2007constraint}]\label{def:variablelock}
\fi
  For a mixed integer program of form~\eqref{eq:mip}, the number of \emph{down-locks} and \emph{up-locks}
  of variable~$\varindex \in \allvars$ is defined as the number of
  positive coefficients per column $\downlock{\varindex} := \Abs{\{\rowindex\;|\; \coefmatrix_{\rowindex\varindex} > 0\}}$
  and the number of negative coefficients per column $\uplock{\varindex} := \Abs{\{\rowindex\;|\; \coefmatrix_{\rowindex\varindex} < 0\}}$, respectively.
\ifreport
\end{Def}
\else
\end{definition}
\smallskip
\fi

Starting at an LP solution $\opt{\var}$ satisfying $\coefmatrix\var \geq \rhsvector$,
a variable $\varindex \in \intvars$ with zero down-locks~$\downlock{\varindex}$ (with zero up-locks $\uplock{\varindex}$)
can always be set to its lower bound (upper bound) without increasing the violation of any constraint.
On the other hand, if a variable has down-locks (up-locks), rounding the variable downwards (upwards) might increase the violation of at least one constraint.
If a constraint down-locks variable~$\varindex$ and is tight \wrt the LP solution candidate $\opt{\var}$, rounding $\opt{\var}_\varindex$ downwards will violate the constraint.
Hence, the variable locks measure the ``risk'' that rounding a variable leads to additional constraint violations.

Usually, only variable locks that are implied by model constraints (short: variable locks) are used during a MIP solve, \eg during presolving, propagation, or primal heuristics.
A well-known diving heuristic that solely relies on variable locks is \emph{\coefficientdiving}~\citep{Berthold2008}.
\coefficientdiving as implemented in \scip follows the diving scheme of Algorithm~\ref{alg:divingalgo}.
For every variable the rounding direction is determined based on the variable locks only.
For a variable~$\varindex$ that is locked in both directions, \ie $\uplock{\varindex} > 0$ and $\downlock{\varindex} > 0$,
\coefficientdiving prefers the direction with less variable locks, \ie the ``safe'' direction.
The order to process the variables during diving is given by the number of variable locks, too.
Here, variables with many variable locks on the chosen direction are preferred.

Processing variables first that tend to lead to infeasibilities is often called a \emph{fail fast strategy}.
This strategy was introduced by~\cite{haralick1980increasing} in the context of artificial intelligence and constraint programming.
Later, \cite{Berthold2014} verified that taking the most critical decisions first is a good strategy for MIP.

Variable locks are a very static criterion and include also model constraints that either do not propagate frequently or are not tight at the current LP relaxation.
For this reason, we propose to consider also variable locks implied only by conflict constraints (short: \emph{conflict locks}),
which can be defined analogous to Definition~\ref{def:variablelock}.
The following lemma shows that conflict locks can have the effect of measuring the ``risk'' more accurately.

\begin{Lem}\label{lem:conflictlocks}
  Let $(\dualvector^\T\coefmatrix)\var \geq \dualvector^\T\rhsvector$ be a Farkas proof derived from an infeasible LP
  and $\dualsol$ a dual ray proving the infeasibility of this local subproblem.
  If the conflict contributes to the conflict up-locks (conflict down-locks) of $\varindex$, then,
  there exists at least one model constraint that contributes to the variable up-locks (variable down-locks) of $\varindex$.
\end{Lem}

\ifreport
\begin{proof}
\else
\proof{Proof.}
\fi
  The coefficient of $\var_\varindex$ in the conflict constraint $(\dualvector^\T\coefmatrix)\var \geq \dualvector^\T\rhsvector$
  is $ \bar{a}_\varindex := \sum_{k = 1}^{m} a_{k\varindex} \cdot \dualvector_\rowindex$.
  If $\bar{a}_\varindex < 0$, \ie the conflict up-locks $\var_\varindex$ then there must exist at least one constraint $k$ with $a_{k\varindex} < 0$
  because $\dualvector \geq 0$.
  Hence, constraint $k$ contributes to the variable up-locks of $\var_\varindex$.
  The analogous argument holds for $\bar{a}_\varindex > 0$.
\ifreport
\end{proof}
\else
{\hfill\Halmos}
\endproof
\fi
\smallskip

This motivates the use of conflict locks: they measure the ``risk'' of violating both model and conflict constraints but focus on constraints that have been actively involved in pruning \bandb nodes.

A similar concept which is used in SAT is VSIDS (variable state independent decaying sum)~\citep{moskewicz2001chaff}.
The VSIDS score takes the contribution of every variable (and its negated complement) into account.
For every variable the number of clauses (in MIP speaking: conflict constraints) the variable is part of is counted.
During continuing the search the VSIDS are periodically scaled by a predefined constant.
With this periodically scaling the weight of older clauses is reduced over the time and more recent observation are weighted higher.
In contrast to VSIDS, conflict locks are not periodically scaled
and only respect conflict constraints that are part at the current subproblem.
This is especially interesting within a MIP solver that uses a pool-based approach to maintain the conflict constraints~\citep{witzig2017experiments}.
Therefore, conflict locks give rise to the current set of variables that are involved in conflict constraints, whereas VSIDS also incorporate past conflict information.

While classical \coefficientdiving solely relies on variable locks, we exploit Lemma~\ref{lem:conflictlocks} and use a combination of both variable and conflict locks.
Given a weight $\confweight \in [0,1]$, the \emph{up-weight} $\upweight{\varindex}$ and \emph{down-weight} $\downweight{\varindex}$ of a variable $\varindex$ is given as
a convex combination of both lock types, \ie $\upweight{\varindex} := \confweight \cdot \confuplock{\varindex} + (1 - \confweight) \cdot \uplock{\varindex}$
and $\downweight{\varindex} := \confweight \cdot \confdownlock{\varindex} + (1 - \confweight) \cdot \downlock{\varindex}$,
where $\confuplock{\varindex}$ and $\confdownlock{\varindex}$ denote the number of conflict up-locks and conflict down-locks, respectively.

Furthermore, we revert the rounding strategy of \coefficientdiving: we use a rounding function preferring the direction
that is more likely to lead to infeasibilities whenever variable $\varindex$ has locks in at most one direction,
\begin{align}
  \roundfunc{C}(\varindex) :=
  \begin{cases}
        \UP    &\hspace{-.75em}\text{ if } \upweight{\varindex} > \downweight{\varindex} \text{ or } \upweight{\varindex} = \downweight{\varindex} \wedge \lpfrac{\varindex} \geq \frac{1}{2},\\
        \DOWN  &\hspace{-.75em}\text{ if } \downweight{\varindex} > \upweight{\varindex} \text{ or } \upweight{\varindex} = \downweight{\varindex} \wedge \lpfrac{\varindex} < \frac{1}{2},\\
  \end{cases}
\end{align}
Here we may assume that all variables that have no locks at all, \ie free variables, are already set to their best bound \wrt their objective coefficient.
With this strategy we aim to guide the heuristic into parts of the search tree that are usually not explored by other heuristics.
The order in which the diving candidates are explored is identical to the one used in \coefficientdiving,
\begin{align}
  \scorefunc{C}(\varindex) :=
  \begin{cases}
    \upweight{\varindex} & \text{ if } \roundfunc{C}(\varindex) = \UP, \\
    \downweight{\varindex} & \text{ if } \roundfunc{C}(\varindex) = \DOWN. \\
  \end{cases}
\end{align}
\ie variables that have a large number of locks on the chosen rounding direction are preferred.
With this scoring function we pursue a fail fast strategy in order to reduce the time spent by this heuristic.
We use a fail fast strategy that is even more aggressive then \coefficientdiving since we already choose the critical direction.
As a result, the heuristic tries to process variables that tend to infeasibility at the beginning of the diving path before
the degree of freedom is further reduced during the dive.
In the following, we refer to this diving heuristic as \emph{\conflictdiving}.

\section{Computational Results}
\label{sec:experiments}

In order to investigate whether and how MIP solvers can benefit from the methods presented in this paper,
we carried out an extensive computational study regarding the individual impact of each heuristic.
For each heuristic, we present computational results on its general performance impact, followed by additional experiments that provide more detailed insights
regarding their particular properties.

In the first part of this section, we compare \scip in its default configuration (\default) to \scip extended by \farkasdiving.
We will refer to the latter setting by \farkdiving.
In the second part of this section, we compare the individual impact of \coefficientdiving and \conflictdiving, to which we will refer by \coefdiving and \confdiving, respectively.
As a baseline we use \scip without both lock-exploiting diving heuristics (\nolockdiving).

All experiments were performed with the academic MIP solver \scip~\citep{GleixnerBastubbeEifleretal.2018} (git hash bf6a486, based on
\scipv), using \soplexv as LP solver.
To evaluate the generated data the \emph{interactive performance evaluation tool} (IPET)~\citep{ipet} was used.
The experiments were run on a cluster of identical machines equipped with Intel Xeon E5-2690 CPUs with 2.6\,GHz and 128\,GB of RAM; a
time limit of $7200$ seconds was set.
To account for the effect of performance variability~\citep{danna2008performance,lodi2013performance} all experiments were performed with four different global random seeds.
As test set we used a union of \miplib~3~\citep{BixbyCeriaMcZealSavelsbergh1998}, \miplib~2003~\citep{AchterbergKochMartin2006}, \miplib~2010~\citep{KochEtAl2011},
and the \coral~\citep{linderoth2005noncommercial} benchmark set.
After removing all duplicates and problems that are known to be numerically unstable, the test set consists of 488 publicly available MIP problems,
which we will refer to as \MMMC.
Every pair of MIP problem and seed is treated as an individual observation, effectively resulting in a test set of $1952$ instances.
We will use the term ``instance'' when referring to a problem-seed combination.

Aggregated results over all random seeds are shown in Table~\ref{tab:farkasdiving_aggregated} and Table~\ref{tab:conflictdiving_aggregated}.
Here, $5$ and $11$ instances, respectively, are excluded because at least one setting finished with numerical violations.
Besides the results on \MMMC, the tables state the impact on affected instances,
\ie instances for which the solving path differs among settings.
Further, the subset of affected instances is grouped into a hierarchy of increasingly harder classes \bracket{$k$}{tilim}.
Class \bracket{$k$}{tilim} contains all instances for which all settings need at least $k$ seconds
and can be solved by at least one setting within the time limit.
As explained by~\cite{AchterbergWunderling2013}, this excludes instances that are ``easy'' for all settings in an unbiased manner.
Detailed tables with instance-wise computational results can be found in the electronic supplement.

\subsection{Farkas Diving}
\label{subsec:compresults-farkas}

We first present computational experiments regarding the general impact of \farkasdiving as proposed in Section~\ref{sec:farkasdiving}.
In this setup, \farkasdiving solved the LP relaxation at every diving node (cf.~Line~\ref{line:solveLP}).
Since this strategy is costly compared to other diving heuristics in \scip, \farkasdiving was executed at local nodes of the search tree,
if the heuristic could already find a feasible solution at the root node.
Moreover, since \farkasdiving relies on the objective function, we executed the heuristic on instances with a nonzero objective function only.
In the second part, we discuss additional computational experiments used to analyze the effect of the individual components of \farkasdiving.

\medskip

\begin{table*}[t]
\begin{centering}
\ifreport
    \scriptsize
\else
    \footnotesize
\fi
\setlength{\tabcolsep}{.75pt}
\caption{Aggregated computational results for Farkas diving on \MMMC over four different random seeds.
          \ifhighlighted Relative changes by at least $5\%$ are highlighted in bold. %
          \fi}\label{tab:farkasdiving_aggregated}
\begin{tabularx}{\textwidth}{LR *{10}{R}}
\toprule
 &  & \multicolumn{5}{c}{\default} & \multicolumn{5}{c}{\farkdiving}\\ 
  \cmidrule(lr){3-7}  \cmidrule(lr){8-12}
 & \# & \solved & \time & \nodes & \virtualbestdsols & \virtualbestdbsols & \solved & \timeQ & \nodesQ & \sols & \bestsols \\ 
\midrule
\cleaninst                & 1947 &   1401 &       183 &      3366 &     2.653 &     0.806 &   1408  &              0.982   & \better{     0.947}  &     3.725  &     0.335  \\
\miplib~2010              &  348 &    305 &       349 &      6225 &     3.454 &     0.828 &    307  &              0.966   & \better{     0.940}  &     4.862  &     0.379  \\
\midrule 
\affected                 &  706 &    686 &        66 &      1204 &     4.238 &     0.564 &    693  & \better{     0.946}  & \better{     0.852}  &     5.924  &     0.585  \\
\bracket{10}{tilim}       &  529 &    509 &       186 &      2420 &     5.168 &     0.733 &    516  & \better{     0.918}  & \better{     0.808}  &     6.112  &     0.573  \\
\bracket{100}{tilim}      &  297 &    277 &       697 &      4444 &     6.997 &     0.980 &    284  & \better{     0.852}  & \better{     0.737}  &     6.862  &     0.690  \\
\bracket{1000}{tilim}     &  127 &    107 &      2599 &      8072 &    11.764 &     1.291 &    114  & \better{     0.837}  & \better{     0.715}  &    10.803  &     1.150  \\
\affectedsols             &  273 &    270 &        34 &      1095 &     9.747 &     0.689 &    272  & \better{     0.879}  & \better{     0.783}  &    15.319  &     1.513  \\
\bottomrule\\
\end{tabularx}
\end{centering}
\end{table*}
 
\paragraph{Overall Impact of Farkas Diving.}

Aggregated computational results comparing \scip in its default configuration (\default)
and \scip extended by \farkasdiving (\farkdiving) on \MMMC are shown in Table~\ref{tab:farkasdiving_aggregated}.
Due to the restricted execution strategy described above, \farkdiving affected only $706$ out of $1947$ instances.
However, on these $706$~instances \farkdiving led to a speed-up of $\percent{5.4}$ (\timeQ) and reduced the tree size by $\percent{14.8}$ (\nodesQ).

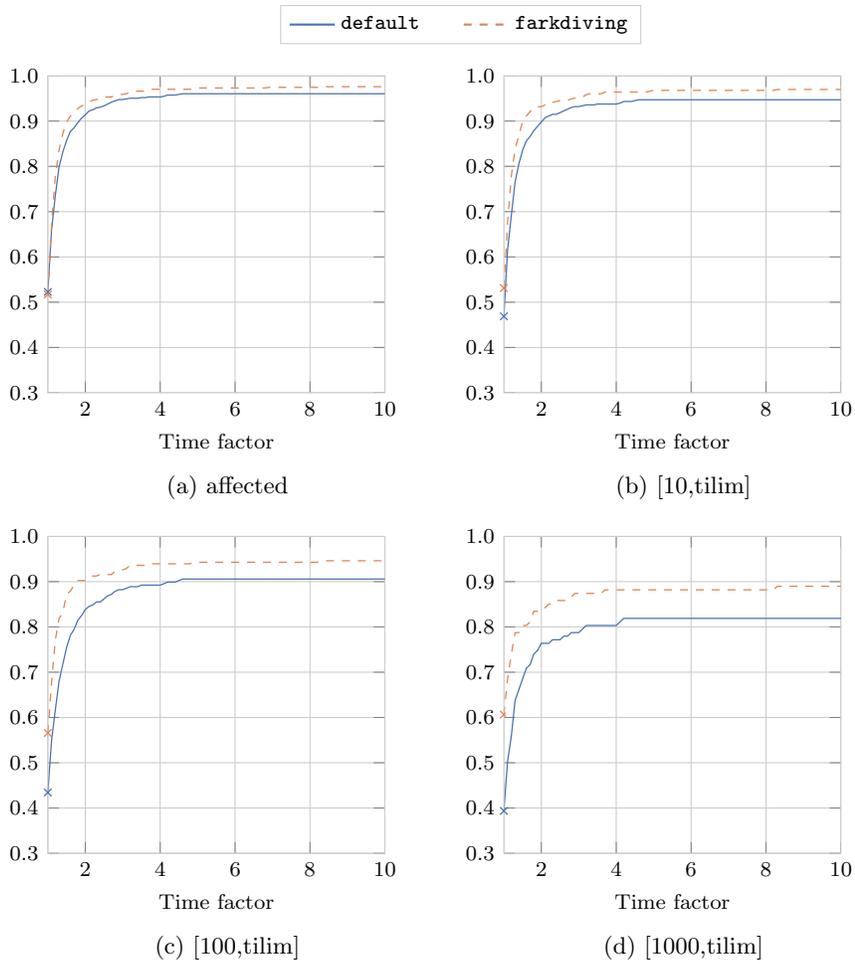
\begin{figure}%
  \begin{center}%
    \begin{subfigure}[t][][t]{\ifreport.5\else.25\fi\textwidth}
        {
        \newlength\figureheight
        \newlength\figurewidth
        \setlength\figureheight{.3\textheight}
        \setlength\figurewidth{\textwidth}
\begin{tikzpicture}

\definecolor{color0}{rgb}{0.298039215686275,0.447058823529412,0.690196078431373}
\definecolor{color1}{rgb}{0.866666666666667,0.517647058823529,0.32156862745098}

\begin{axis}[
axis line style={white!80.0!black},
height=\figureheight,
legend cell align={left},
legend entries={{\default},{\farkdiving}},
legend style={at={\ifreport(1.225,1.1)\else(3,1.10)\fi}, anchor=south, draw=white!80.0!black, font=\footnotesize, /tikz/every even column/.append style={column sep=0.5cm}},
legend columns=2,
tick pos=both,
width=\figurewidth,
x grid style={white!80.0!black},
xlabel={Time factor},
xmajorgrids,
xmin=1, xmax=10,
y grid style={white!80.0!black},
ymajorgrids,
ymin=0.3, ymax=1,
ytick={0.2,0.3,0.4,0.5,0.6,0.7,0.8,0.9,1},
yticklabels={,0.3,0.4,0.5,0.6,0.7,0.8,0.9,1.0}
]
\addlegendimage{no markers, thick, color0}
\addlegendimage{no markers, thick, color1, dashed}
\addplot [line width=0.4800000000000001pt, color0, forget plot, mark=x]
table [row sep=\\]{%
1	0.522663 \\
};
\addplot [line width=0.4800000000000001pt, color0, forget plot]
table [row sep=\\]{%
1	0.522663 \\
1.1	0.660057 \\
1.2	0.735127 \\
1.3	0.798867 \\
1.4	0.831445 \\
1.5	0.856941 \\
1.6	0.876771 \\
1.7	0.885269 \\
1.8	0.896601 \\
1.9	0.906516 \\
2	0.913598 \\
2.1	0.922096 \\
2.2	0.924929 \\
2.3	0.929178 \\
2.4	0.930595 \\
2.5	0.933428 \\
2.6	0.937677 \\
2.7	0.941926 \\
2.8	0.944759 \\
2.9	0.947592 \\
3	0.947592 \\
3.1	0.949008 \\
3.2	0.950425 \\
3.3	0.950425 \\
3.4	0.950425 \\
3.5	0.951841 \\
3.6	0.951841 \\
3.7	0.953258 \\
3.8	0.953258 \\
3.9	0.953258 \\
4	0.953258 \\
4.1	0.954674 \\
4.2	0.957507 \\
4.3	0.957507 \\
4.4	0.957507 \\
4.5	0.958924 \\
4.6	0.96034 \\
4.7	0.96034 \\
4.8	0.96034 \\
4.9	0.96034 \\
5	0.96034 \\
5.1	0.96034 \\
5.2	0.96034 \\
5.3	0.96034 \\
5.4	0.96034 \\
5.5	0.96034 \\
5.6	0.96034 \\
5.7	0.96034 \\
5.8	0.96034 \\
5.9	0.96034 \\
6	0.96034 \\
6.1	0.96034 \\
6.19999999999999	0.96034 \\
6.29999999999999	0.96034 \\
6.39999999999999	0.96034 \\
6.49999999999999	0.96034 \\
6.59999999999999	0.96034 \\
6.69999999999999	0.96034 \\
6.79999999999999	0.96034 \\
6.89999999999999	0.96034 \\
6.99999999999999	0.96034 \\
7.09999999999999	0.96034 \\
7.19999999999999	0.96034 \\
7.29999999999999	0.96034 \\
7.39999999999999	0.96034 \\
7.49999999999999	0.96034 \\
7.59999999999999	0.96034 \\
7.69999999999999	0.96034 \\
7.79999999999999	0.96034 \\
7.89999999999999	0.96034 \\
7.99999999999999	0.96034 \\
8.09999999999999	0.96034 \\
8.19999999999999	0.96034 \\
8.29999999999999	0.96034 \\
8.39999999999999	0.96034 \\
8.49999999999999	0.96034 \\
8.59999999999999	0.96034 \\
8.69999999999999	0.96034 \\
8.79999999999999	0.96034 \\
8.89999999999999	0.96034 \\
8.99999999999999	0.96034 \\
9.09999999999999	0.96034 \\
9.19999999999999	0.96034 \\
9.29999999999998	0.96034 \\
9.39999999999998	0.96034 \\
9.49999999999998	0.96034 \\
9.59999999999998	0.96034 \\
9.69999999999998	0.96034 \\
9.79999999999998	0.96034 \\
9.89999999999998	0.96034 \\
9.99999999999998	0.96034 \\
};
\addplot [line width=0.4800000000000001pt, color1, forget plot, mark=x]
table [row sep=\\]{%
1	0.516997 \\
};
\addplot [line width=0.4800000000000001pt, color1, dashed, forget plot]
table [row sep=\\]{%
1	0.516997 \\
1.1	0.671388 \\
1.2	0.779037 \\
1.3	0.835694 \\
1.4	0.872521 \\
1.5	0.898017 \\
1.6	0.910765 \\
1.7	0.919263 \\
1.8	0.927762 \\
1.9	0.933428 \\
2	0.936261 \\
2.1	0.941926 \\
2.2	0.946176 \\
2.3	0.947592 \\
2.4	0.950425 \\
2.5	0.953258 \\
2.6	0.953258 \\
2.7	0.953258 \\
2.8	0.956091 \\
2.9	0.958924 \\
3	0.958924 \\
3.1	0.96034 \\
3.2	0.964589 \\
3.3	0.966006 \\
3.4	0.966006 \\
3.5	0.966006 \\
3.6	0.966006 \\
3.7	0.968839 \\
3.8	0.970255 \\
3.9	0.970255 \\
4	0.970255 \\
4.1	0.970255 \\
4.2	0.970255 \\
4.3	0.970255 \\
4.4	0.970255 \\
4.5	0.970255 \\
4.6	0.970255 \\
4.7	0.970255 \\
4.8	0.970255 \\
4.9	0.971671 \\
5	0.973088 \\
5.1	0.973088 \\
5.2	0.973088 \\
5.3	0.973088 \\
5.4	0.973088 \\
5.5	0.973088 \\
5.6	0.973088 \\
5.7	0.973088 \\
5.8	0.973088 \\
5.9	0.973088 \\
6	0.973088 \\
6.1	0.973088 \\
6.19999999999999	0.973088 \\
6.29999999999999	0.973088 \\
6.39999999999999	0.973088 \\
6.49999999999999	0.973088 \\
6.59999999999999	0.973088 \\
6.69999999999999	0.973088 \\
6.79999999999999	0.973088 \\
6.89999999999999	0.974504 \\
6.99999999999999	0.974504 \\
7.09999999999999	0.974504 \\
7.19999999999999	0.974504 \\
7.29999999999999	0.974504 \\
7.39999999999999	0.974504 \\
7.49999999999999	0.974504 \\
7.59999999999999	0.974504 \\
7.69999999999999	0.974504 \\
7.79999999999999	0.974504 \\
7.89999999999999	0.974504 \\
7.99999999999999	0.974504 \\
8.09999999999999	0.974504 \\
8.19999999999999	0.974504 \\
8.29999999999999	0.975921 \\
8.39999999999999	0.975921 \\
8.49999999999999	0.975921 \\
8.59999999999999	0.975921 \\
8.69999999999999	0.975921 \\
8.79999999999999	0.975921 \\
8.89999999999999	0.975921 \\
8.99999999999999	0.975921 \\
9.09999999999999	0.975921 \\
9.19999999999999	0.975921 \\
9.29999999999998	0.975921 \\
9.39999999999998	0.975921 \\
9.49999999999998	0.975921 \\
9.59999999999998	0.975921 \\
9.69999999999998	0.975921 \\
9.79999999999998	0.975921 \\
9.89999999999998	0.975921 \\
9.99999999999998	0.975921 \\
};
\path [draw=white!80.0!black, fill opacity=0] (axis cs:0,0.3)
--(axis cs:0,1);

\path [draw=white!80.0!black, fill opacity=0] (axis cs:1,0.3)
--(axis cs:1,1);

\path [draw=white!80.0!black, fill opacity=0] (axis cs:1,0)
--(axis cs:10,0);

\path [draw=white!80.0!black, fill opacity=0] (axis cs:1,1)
--(axis cs:10,1);

\end{axis}

\end{tikzpicture}
        \caption{\affected}
        }
    \end{subfigure}%
    \hfill%
    \begin{subfigure}[t][][t]{\ifreport.5\else.25\fi\textwidth}
        {
        \setlength\figureheight{.3\textheight}
        \setlength\figurewidth{\textwidth}
\begin{tikzpicture}

\definecolor{color0}{rgb}{0.298039215686275,0.447058823529412,0.690196078431373}
\definecolor{color1}{rgb}{0.866666666666667,0.517647058823529,0.32156862745098}

\begin{axis}[
axis line style={white!80.0!black},
height=\figureheight,
tick pos=both,
width=\figurewidth,
x grid style={white!80.0!black},
xlabel={Time factor},
xmajorgrids,
xmin=1, xmax=10,
y grid style={white!80.0!black},
ymajorgrids,
ymin=0.3, ymax=1,
ytick={0.2,0.3,0.4,0.5,0.6,0.7,0.8,0.9,1},
yticklabels={,0.3,0.4,0.5,0.6,0.7,0.8,0.9,1.0}
]
\addlegendimage{no markers, color0}
\addlegendimage{no markers, color1}
\addplot [line width=0.4800000000000001pt, color0, forget plot, mark=x]
table [row sep=\\]{%
1	0.468809 \\
};
\addplot [line width=0.4800000000000001pt, color0, forget plot]
table [row sep=\\]{%
1	0.468809 \\
1.1	0.616257 \\
1.2	0.691871 \\
1.3	0.765595 \\
1.4	0.805293 \\
1.5	0.835539 \\
1.6	0.856333 \\
1.7	0.865784 \\
1.8	0.879017 \\
1.9	0.888469 \\
2	0.897921 \\
2.1	0.907372 \\
2.2	0.911153 \\
2.3	0.914934 \\
2.4	0.914934 \\
2.5	0.918715 \\
2.6	0.922495 \\
2.7	0.926276 \\
2.8	0.930057 \\
2.9	0.931947 \\
3	0.931947 \\
3.1	0.933837 \\
3.2	0.935728 \\
3.3	0.935728 \\
3.4	0.935728 \\
3.5	0.937618 \\
3.6	0.937618 \\
3.7	0.937618 \\
3.8	0.937618 \\
3.9	0.937618 \\
4	0.937618 \\
4.1	0.939509 \\
4.2	0.943289 \\
4.3	0.943289 \\
4.4	0.943289 \\
4.5	0.94518 \\
4.6	0.94707 \\
4.7	0.94707 \\
4.8	0.94707 \\
4.9	0.94707 \\
5	0.94707 \\
5.1	0.94707 \\
5.2	0.94707 \\
5.3	0.94707 \\
5.4	0.94707 \\
5.5	0.94707 \\
5.6	0.94707 \\
5.7	0.94707 \\
5.8	0.94707 \\
5.9	0.94707 \\
6	0.94707 \\
6.1	0.94707 \\
6.19999999999999	0.94707 \\
6.29999999999999	0.94707 \\
6.39999999999999	0.94707 \\
6.49999999999999	0.94707 \\
6.59999999999999	0.94707 \\
6.69999999999999	0.94707 \\
6.79999999999999	0.94707 \\
6.89999999999999	0.94707 \\
6.99999999999999	0.94707 \\
7.09999999999999	0.94707 \\
7.19999999999999	0.94707 \\
7.29999999999999	0.94707 \\
7.39999999999999	0.94707 \\
7.49999999999999	0.94707 \\
7.59999999999999	0.94707 \\
7.69999999999999	0.94707 \\
7.79999999999999	0.94707 \\
7.89999999999999	0.94707 \\
7.99999999999999	0.94707 \\
8.09999999999999	0.94707 \\
8.19999999999999	0.94707 \\
8.29999999999999	0.94707 \\
8.39999999999999	0.94707 \\
8.49999999999999	0.94707 \\
8.59999999999999	0.94707 \\
8.69999999999999	0.94707 \\
8.79999999999999	0.94707 \\
8.89999999999999	0.94707 \\
8.99999999999999	0.94707 \\
9.09999999999999	0.94707 \\
9.19999999999999	0.94707 \\
9.29999999999998	0.94707 \\
9.39999999999998	0.94707 \\
9.49999999999998	0.94707 \\
9.59999999999998	0.94707 \\
9.69999999999998	0.94707 \\
9.79999999999998	0.94707 \\
9.89999999999998	0.94707 \\
9.99999999999998	0.94707 \\
};
\addplot [line width=0.4800000000000001pt, color1, forget plot, mark=x]
table [row sep=\\]{%
1	0.531191 \\
};
\addplot [line width=0.4800000000000001pt, color1, dashed, forget plot]
table [row sep=\\]{%
1	0.531191 \\
1.1	0.676749 \\
1.2	0.778828 \\
1.3	0.839319 \\
1.4	0.867675 \\
1.5	0.899811 \\
1.6	0.911153 \\
1.7	0.920605 \\
1.8	0.930057 \\
1.9	0.931947 \\
2	0.931947 \\
2.1	0.935728 \\
2.2	0.941399 \\
2.3	0.941399 \\
2.4	0.943289 \\
2.5	0.94518 \\
2.6	0.94518 \\
2.7	0.94518 \\
2.8	0.94896 \\
2.9	0.950851 \\
3	0.950851 \\
3.1	0.952741 \\
3.2	0.958412 \\
3.3	0.960302 \\
3.4	0.960302 \\
3.5	0.960302 \\
3.6	0.960302 \\
3.7	0.964083 \\
3.8	0.964083 \\
3.9	0.964083 \\
4	0.964083 \\
4.1	0.964083 \\
4.2	0.964083 \\
4.3	0.964083 \\
4.4	0.964083 \\
4.5	0.964083 \\
4.6	0.964083 \\
4.7	0.964083 \\
4.8	0.964083 \\
4.9	0.965974 \\
5	0.967864 \\
5.1	0.967864 \\
5.2	0.967864 \\
5.3	0.967864 \\
5.4	0.967864 \\
5.5	0.967864 \\
5.6	0.967864 \\
5.7	0.967864 \\
5.8	0.967864 \\
5.9	0.967864 \\
6	0.967864 \\
6.1	0.967864 \\
6.19999999999999	0.967864 \\
6.29999999999999	0.967864 \\
6.39999999999999	0.967864 \\
6.49999999999999	0.967864 \\
6.59999999999999	0.967864 \\
6.69999999999999	0.967864 \\
6.79999999999999	0.967864 \\
6.89999999999999	0.967864 \\
6.99999999999999	0.967864 \\
7.09999999999999	0.967864 \\
7.19999999999999	0.967864 \\
7.29999999999999	0.967864 \\
7.39999999999999	0.967864 \\
7.49999999999999	0.967864 \\
7.59999999999999	0.967864 \\
7.69999999999999	0.967864 \\
7.79999999999999	0.967864 \\
7.89999999999999	0.967864 \\
7.99999999999999	0.967864 \\
8.09999999999999	0.967864 \\
8.19999999999999	0.967864 \\
8.29999999999999	0.969754 \\
8.39999999999999	0.969754 \\
8.49999999999999	0.969754 \\
8.59999999999999	0.969754 \\
8.69999999999999	0.969754 \\
8.79999999999999	0.969754 \\
8.89999999999999	0.969754 \\
8.99999999999999	0.969754 \\
9.09999999999999	0.969754 \\
9.19999999999999	0.969754 \\
9.29999999999998	0.969754 \\
9.39999999999998	0.969754 \\
9.49999999999998	0.969754 \\
9.59999999999998	0.969754 \\
9.69999999999998	0.969754 \\
9.79999999999998	0.969754 \\
9.89999999999998	0.969754 \\
9.99999999999998	0.969754 \\
};
\path [draw=white!80.0!black, fill opacity=0] (axis cs:0,0.3)
--(axis cs:0,1);

\path [draw=white!80.0!black, fill opacity=0] (axis cs:1,0.3)
--(axis cs:1,1);

\path [draw=white!80.0!black, fill opacity=0] (axis cs:1,0)
--(axis cs:10,0);

\path [draw=white!80.0!black, fill opacity=0] (axis cs:1,1)
--(axis cs:10,1);

\end{axis}

\end{tikzpicture}
        \caption{\bracket{10}{tilim}}
        }
    \end{subfigure}%
    \ifreport
        \medskip

    \else
        \hfill%
    \fi
    \begin{subfigure}[t][][t]{\ifreport.5\else.25\fi\textwidth}
        {
        \setlength\figureheight{.3\textheight}
        \setlength\figurewidth{\textwidth}
\begin{tikzpicture}

\definecolor{color0}{rgb}{0.298039215686275,0.447058823529412,0.690196078431373}
\definecolor{color1}{rgb}{0.866666666666667,0.517647058823529,0.32156862745098}

\begin{axis}[
axis line style={white!80.0!black},
height=\figureheight,
tick pos=both,
width=\figurewidth,
x grid style={white!80.0!black},
xlabel={Time factor},
xmajorgrids,
xmin=1, xmax=10,
y grid style={white!80.0!black},
ymajorgrids,
ymin=0.3, ymax=1,
ytick={0.2,0.3,0.4,0.5,0.6,0.7,0.8,0.9,1},
yticklabels={,0.3,0.4,0.5,0.6,0.7,0.8,0.9,1.0}
]
\addlegendimage{no markers, color0}
\addlegendimage{no markers, color1}
\addplot [line width=0.4800000000000001pt, color0, forget plot, mark=x]
table [row sep=\\]{%
1	0.434343 \\
};
\addplot [line width=0.4800000000000001pt, color0, forget plot]
table [row sep=\\]{%
1	0.434343 \\
1.1	0.548822 \\
1.2	0.616162 \\
1.3	0.680135 \\
1.4	0.717172 \\
1.5	0.754209 \\
1.6	0.781145 \\
1.7	0.794613 \\
1.8	0.814815 \\
1.9	0.824916 \\
2	0.838384 \\
2.1	0.845118 \\
2.2	0.848485 \\
2.3	0.855219 \\
2.4	0.855219 \\
2.5	0.861953 \\
2.6	0.868687 \\
2.7	0.872054 \\
2.8	0.878788 \\
2.9	0.882155 \\
3	0.882155 \\
3.1	0.885522 \\
3.2	0.888889 \\
3.3	0.888889 \\
3.4	0.888889 \\
3.5	0.892256 \\
3.6	0.892256 \\
3.7	0.892256 \\
3.8	0.892256 \\
3.9	0.892256 \\
4	0.892256 \\
4.1	0.895623 \\
4.2	0.89899 \\
4.3	0.89899 \\
4.4	0.89899 \\
4.5	0.902357 \\
4.6	0.905724 \\
4.7	0.905724 \\
4.8	0.905724 \\
4.9	0.905724 \\
5	0.905724 \\
5.1	0.905724 \\
5.2	0.905724 \\
5.3	0.905724 \\
5.4	0.905724 \\
5.5	0.905724 \\
5.6	0.905724 \\
5.7	0.905724 \\
5.8	0.905724 \\
5.9	0.905724 \\
6	0.905724 \\
6.1	0.905724 \\
6.19999999999999	0.905724 \\
6.29999999999999	0.905724 \\
6.39999999999999	0.905724 \\
6.49999999999999	0.905724 \\
6.59999999999999	0.905724 \\
6.69999999999999	0.905724 \\
6.79999999999999	0.905724 \\
6.89999999999999	0.905724 \\
6.99999999999999	0.905724 \\
7.09999999999999	0.905724 \\
7.19999999999999	0.905724 \\
7.29999999999999	0.905724 \\
7.39999999999999	0.905724 \\
7.49999999999999	0.905724 \\
7.59999999999999	0.905724 \\
7.69999999999999	0.905724 \\
7.79999999999999	0.905724 \\
7.89999999999999	0.905724 \\
7.99999999999999	0.905724 \\
8.09999999999999	0.905724 \\
8.19999999999999	0.905724 \\
8.29999999999999	0.905724 \\
8.39999999999999	0.905724 \\
8.49999999999999	0.905724 \\
8.59999999999999	0.905724 \\
8.69999999999999	0.905724 \\
8.79999999999999	0.905724 \\
8.89999999999999	0.905724 \\
8.99999999999999	0.905724 \\
9.09999999999999	0.905724 \\
9.19999999999999	0.905724 \\
9.29999999999998	0.905724 \\
9.39999999999998	0.905724 \\
9.49999999999998	0.905724 \\
9.59999999999998	0.905724 \\
9.69999999999998	0.905724 \\
9.79999999999998	0.905724 \\
9.89999999999998	0.905724 \\
9.99999999999998	0.905724 \\
};
\addplot [line width=0.4800000000000001pt, color1, forget plot, mark=x]
table [row sep=\\]{%
1	0.565657 \\
};
\addplot [line width=0.4800000000000001pt, color1, dashed, forget plot]
table [row sep=\\]{%
1	0.565657 \\
1.1	0.676768 \\
1.2	0.771044 \\
1.3	0.818182 \\
1.4	0.83165 \\
1.5	0.872054 \\
1.6	0.878788 \\
1.7	0.892256 \\
1.8	0.902357 \\
1.9	0.902357 \\
2	0.902357 \\
2.1	0.905724 \\
2.2	0.912458 \\
2.3	0.912458 \\
2.4	0.915825 \\
2.5	0.915825 \\
2.6	0.915825 \\
2.7	0.915825 \\
2.8	0.922559 \\
2.9	0.925926 \\
3	0.925926 \\
3.1	0.929293 \\
3.2	0.936027 \\
3.3	0.936027 \\
3.4	0.936027 \\
3.5	0.936027 \\
3.6	0.936027 \\
3.7	0.939394 \\
3.8	0.939394 \\
3.9	0.939394 \\
4	0.939394 \\
4.1	0.939394 \\
4.2	0.939394 \\
4.3	0.939394 \\
4.4	0.939394 \\
4.5	0.939394 \\
4.6	0.939394 \\
4.7	0.939394 \\
4.8	0.939394 \\
4.9	0.942761 \\
5	0.942761 \\
5.1	0.942761 \\
5.2	0.942761 \\
5.3	0.942761 \\
5.4	0.942761 \\
5.5	0.942761 \\
5.6	0.942761 \\
5.7	0.942761 \\
5.8	0.942761 \\
5.9	0.942761 \\
6	0.942761 \\
6.1	0.942761 \\
6.19999999999999	0.942761 \\
6.29999999999999	0.942761 \\
6.39999999999999	0.942761 \\
6.49999999999999	0.942761 \\
6.59999999999999	0.942761 \\
6.69999999999999	0.942761 \\
6.79999999999999	0.942761 \\
6.89999999999999	0.942761 \\
6.99999999999999	0.942761 \\
7.09999999999999	0.942761 \\
7.19999999999999	0.942761 \\
7.29999999999999	0.942761 \\
7.39999999999999	0.942761 \\
7.49999999999999	0.942761 \\
7.59999999999999	0.942761 \\
7.69999999999999	0.942761 \\
7.79999999999999	0.942761 \\
7.89999999999999	0.942761 \\
7.99999999999999	0.942761 \\
8.09999999999999	0.942761 \\
8.19999999999999	0.942761 \\
8.29999999999999	0.946128 \\
8.39999999999999	0.946128 \\
8.49999999999999	0.946128 \\
8.59999999999999	0.946128 \\
8.69999999999999	0.946128 \\
8.79999999999999	0.946128 \\
8.89999999999999	0.946128 \\
8.99999999999999	0.946128 \\
9.09999999999999	0.946128 \\
9.19999999999999	0.946128 \\
9.29999999999998	0.946128 \\
9.39999999999998	0.946128 \\
9.49999999999998	0.946128 \\
9.59999999999998	0.946128 \\
9.69999999999998	0.946128 \\
9.79999999999998	0.946128 \\
9.89999999999998	0.946128 \\
9.99999999999998	0.946128 \\
};
\path [draw=white!80.0!black, fill opacity=0] (axis cs:0,0.3)
--(axis cs:0,1);

\path [draw=white!80.0!black, fill opacity=0] (axis cs:1,0.3)
--(axis cs:1,1);

\path [draw=white!80.0!black, fill opacity=0] (axis cs:1,0)
--(axis cs:10,0);

\path [draw=white!80.0!black, fill opacity=0] (axis cs:1,1)
--(axis cs:10,1);

\end{axis}

\end{tikzpicture}
        \caption{\bracket{100}{tilim}}
        }
    \end{subfigure}%
    \hfill%
    \begin{subfigure}[t][][t]{\ifreport.5\else.25\fi\textwidth}
        {
        \setlength\figureheight{.3\textheight}
        \setlength\figurewidth{\textwidth}
\begin{tikzpicture}

\definecolor{color0}{rgb}{0.298039215686275,0.447058823529412,0.690196078431373}
\definecolor{color1}{rgb}{0.866666666666667,0.517647058823529,0.32156862745098}

\begin{axis}[
axis line style={white!80.0!black},
height=\figureheight,
tick pos=both,
width=\figurewidth,
x grid style={white!80.0!black},
xlabel={Time factor},
xmajorgrids,
xmin=1, xmax=10,
y grid style={white!80.0!black},
ymajorgrids,
ymin=0.3, ymax=1,
ytick={0.2,0.3,0.4,0.5,0.6,0.7,0.8,0.9,1},
yticklabels={,0.3,0.4,0.5,0.6,0.7,0.8,0.9,1.0}
]
\addlegendimage{no markers, color0}
\addlegendimage{no markers, color1}
\addplot [line width=0.4800000000000001pt, color0, forget plot, mark=x]
table [row sep=\\]{%
1	0.393701 \\
};
\addplot [line width=0.4800000000000001pt, color0, forget plot]
table [row sep=\\]{%
1	0.393701 \\
1.1	0.503937 \\
1.2	0.559055 \\
1.3	0.637795 \\
1.4	0.661417 \\
1.5	0.685039 \\
1.6	0.708661 \\
1.7	0.716535 \\
1.8	0.740157 \\
1.9	0.748031 \\
2	0.76378 \\
2.1	0.76378 \\
2.2	0.76378 \\
2.3	0.771654 \\
2.4	0.771654 \\
2.5	0.771654 \\
2.6	0.779528 \\
2.7	0.779528 \\
2.8	0.787402 \\
2.9	0.787402 \\
3	0.787402 \\
3.1	0.795276 \\
3.2	0.80315 \\
3.3	0.80315 \\
3.4	0.80315 \\
3.5	0.80315 \\
3.6	0.80315 \\
3.7	0.80315 \\
3.8	0.80315 \\
3.9	0.80315 \\
4	0.80315 \\
4.1	0.811024 \\
4.2	0.818898 \\
4.3	0.818898 \\
4.4	0.818898 \\
4.5	0.818898 \\
4.6	0.818898 \\
4.7	0.818898 \\
4.8	0.818898 \\
4.9	0.818898 \\
5	0.818898 \\
5.1	0.818898 \\
5.2	0.818898 \\
5.3	0.818898 \\
5.4	0.818898 \\
5.5	0.818898 \\
5.6	0.818898 \\
5.7	0.818898 \\
5.8	0.818898 \\
5.9	0.818898 \\
6	0.818898 \\
6.1	0.818898 \\
6.19999999999999	0.818898 \\
6.29999999999999	0.818898 \\
6.39999999999999	0.818898 \\
6.49999999999999	0.818898 \\
6.59999999999999	0.818898 \\
6.69999999999999	0.818898 \\
6.79999999999999	0.818898 \\
6.89999999999999	0.818898 \\
6.99999999999999	0.818898 \\
7.09999999999999	0.818898 \\
7.19999999999999	0.818898 \\
7.29999999999999	0.818898 \\
7.39999999999999	0.818898 \\
7.49999999999999	0.818898 \\
7.59999999999999	0.818898 \\
7.69999999999999	0.818898 \\
7.79999999999999	0.818898 \\
7.89999999999999	0.818898 \\
7.99999999999999	0.818898 \\
8.09999999999999	0.818898 \\
8.19999999999999	0.818898 \\
8.29999999999999	0.818898 \\
8.39999999999999	0.818898 \\
8.49999999999999	0.818898 \\
8.59999999999999	0.818898 \\
8.69999999999999	0.818898 \\
8.79999999999999	0.818898 \\
8.89999999999999	0.818898 \\
8.99999999999999	0.818898 \\
9.09999999999999	0.818898 \\
9.19999999999999	0.818898 \\
9.29999999999998	0.818898 \\
9.39999999999998	0.818898 \\
9.49999999999998	0.818898 \\
9.59999999999998	0.818898 \\
9.69999999999998	0.818898 \\
9.79999999999998	0.818898 \\
9.89999999999998	0.818898 \\
9.99999999999998	0.818898 \\
};
\addplot [line width=0.4800000000000001pt, color1, dashed, forget plot, mark=x]
table [row sep=\\]{%
1	0.606299 \\
};
\addplot [line width=0.4800000000000001pt, color1, dashed, forget plot]
table [row sep=\\]{%
1	0.606299 \\
1.1	0.685039 \\
1.2	0.740157 \\
1.3	0.787402 \\
1.4	0.787402 \\
1.5	0.80315 \\
1.6	0.80315 \\
1.7	0.811024 \\
1.8	0.834646 \\
1.9	0.834646 \\
2	0.834646 \\
2.1	0.84252 \\
2.2	0.850394 \\
2.3	0.850394 \\
2.4	0.858268 \\
2.5	0.858268 \\
2.6	0.858268 \\
2.7	0.858268 \\
2.8	0.866142 \\
2.9	0.874016 \\
3	0.874016 \\
3.1	0.874016 \\
3.2	0.874016 \\
3.3	0.874016 \\
3.4	0.874016 \\
3.5	0.874016 \\
3.6	0.874016 \\
3.7	0.88189 \\
3.8	0.88189 \\
3.9	0.88189 \\
4	0.88189 \\
4.1	0.88189 \\
4.2	0.88189 \\
4.3	0.88189 \\
4.4	0.88189 \\
4.5	0.88189 \\
4.6	0.88189 \\
4.7	0.88189 \\
4.8	0.88189 \\
4.9	0.88189 \\
5	0.88189 \\
5.1	0.88189 \\
5.2	0.88189 \\
5.3	0.88189 \\
5.4	0.88189 \\
5.5	0.88189 \\
5.6	0.88189 \\
5.7	0.88189 \\
5.8	0.88189 \\
5.9	0.88189 \\
6	0.88189 \\
6.1	0.88189 \\
6.19999999999999	0.88189 \\
6.29999999999999	0.88189 \\
6.39999999999999	0.88189 \\
6.49999999999999	0.88189 \\
6.59999999999999	0.88189 \\
6.69999999999999	0.88189 \\
6.79999999999999	0.88189 \\
6.89999999999999	0.88189 \\
6.99999999999999	0.88189 \\
7.09999999999999	0.88189 \\
7.19999999999999	0.88189 \\
7.29999999999999	0.88189 \\
7.39999999999999	0.88189 \\
7.49999999999999	0.88189 \\
7.59999999999999	0.88189 \\
7.69999999999999	0.88189 \\
7.79999999999999	0.88189 \\
7.89999999999999	0.88189 \\
7.99999999999999	0.88189 \\
8.09999999999999	0.88189 \\
8.19999999999999	0.88189 \\
8.29999999999999	0.889764 \\
8.39999999999999	0.889764 \\
8.49999999999999	0.889764 \\
8.59999999999999	0.889764 \\
8.69999999999999	0.889764 \\
8.79999999999999	0.889764 \\
8.89999999999999	0.889764 \\
8.99999999999999	0.889764 \\
9.09999999999999	0.889764 \\
9.19999999999999	0.889764 \\
9.29999999999998	0.889764 \\
9.39999999999998	0.889764 \\
9.49999999999998	0.889764 \\
9.59999999999998	0.889764 \\
9.69999999999998	0.889764 \\
9.79999999999998	0.889764 \\
9.89999999999998	0.889764 \\
9.99999999999998	0.889764 \\
};
\path [draw=white!80.0!black, fill opacity=0] (axis cs:0,0.3)
--(axis cs:0,1);

\path [draw=white!80.0!black, fill opacity=0] (axis cs:1,0.3)
--(axis cs:1,1);

\path [draw=white!80.0!black, fill opacity=0] (axis cs:1,0)
--(axis cs:10,0);

\path [draw=white!80.0!black, fill opacity=0] (axis cs:1,1)
--(axis cs:10,1);

\end{axis}

\end{tikzpicture}
        \caption{\bracket{1000}{tilim}}
        }
    \end{subfigure}%
  \end{center}
  \smallskip

  \caption{Performance profiles of \default and \farkdiving for four hierarchical groups of increasingly hard, affected instances.}
  \label{fig:farkdiving_performance}
\end{figure}

On the subset of harder instances \bracket{100}{tilim} the overall performance could even be improved by $\percent{14.8}$
and the size of the search tree could be reduced by more than $\percent{25}$.
The observed performance improvement is spread over the complete group of affected instances.
This becomes apparent also in the performance profiles~\citep{dolan2002benchmarking,gould2016note} displayed in Figure~\ref{fig:farkdiving_performance}.
For a time factor of $1.0$ the respective amount of instances that could be solved best with the respective setting is marked with a colored cross.
In our experiments, we observed that on the set of affected instances \default performs best on $369$ instances,
whereas \farkdiving was marginal worse and performs best on $365$ affected instances.
On this group of instances, both profiles cross exactly once at time factor $1.021$.
Afterward, \farkdiving is always superior to \default.
This fact indicates that \farkdiving is especially superior to \default on harder instances and
is confirmed by the performance profiles over instances where both settings need at least $10$, $100$ or $1000$ settings.
On these groups of affected instances, \farkdiving is clearly superior to \default.
For a time factor of $1.0$ \farkdiving performs best the harder the instances are.
On the group of affected instances for which both settings need at least $10$ seconds $248$ can be solved best by \default,
whereas \farkdiving performs best on $281$ instances.
These results indicate that \farkasdiving is expensive on easy instances compared to \default, but superior on harder instances.

\paragraph{Success in Generating Conflicts.}

\begin{figure}%
    \centering
    \setlength\figureheight{.15\textheight}
    \setlength\figurewidth{\ifreport.75\else.5\fi\textwidth}
\begin{tikzpicture}

\definecolor{color0}{rgb}{0.347058823529412,0.458823529411765,0.641176470588235}
\definecolor{color1}{rgb}{0.798529411764706,0.536764705882353,0.389705882352941}
\definecolor{color2}{rgb}{0.298039215686275,0.447058823529412,0.690196078431373}

\begin{groupplot}[group style={group size=1 by 2, vertical sep=.25cm}]
\nextgroupplot[
axis line style={white!80.0!black},
height=\figureheight,
legend cell align={left},
legend entries={{\default},{\farkdiving}},
legend style={draw=white!80.0!black, at={(0.5,1.1)},anchor=south, /tikz/every even column/.append style={column sep=0.5cm}, font=\footnotesize},
legend columns=2,
tick pos=both,
width=\figurewidth,
x grid style={white!80.0!black},
xmin=-0.5, xmax=3.5,
xtick={0,1,2,3},
xticklabels={},
y grid style={white!80.0!black},
ylabel={Generated Conflicts},
ylabel style={at={(axis description cs:-0.15,-.01)},anchor=south},
ymajorgrids,
ymin=60, ymax=50000,
ymode=log,
ytick={0.1,1,10,100,1000,10000},
yticklabels={,,$10$,$100$,$1000$,$10000$}
]

\addlegendimage{no markers,color0,thick}
\addlegendimage{no markers,color1,thick}

\path [draw=white!29.80392156862745!black, fill=color0, opacity=.65] (axis cs:-0.22275,0)
--(axis cs:-0.00225,0)
--(axis cs:-0.00225,46)
--(axis cs:-0.22275,46)
--(axis cs:-0.22275,0)
--cycle;

\path [draw=white!29.80392156862745!black, fill=color1, opacity=.65] (axis cs:0.00225,0)
--(axis cs:0.22275,0)
--(axis cs:0.22275,89)
--(axis cs:0.00225,89)
--(axis cs:0.00225,0)
--cycle;

\path [draw=white!29.80392156862745!black, fill=color0, opacity=.65] (axis cs:0.77725,0)
--(axis cs:0.99775,0)
--(axis cs:0.99775,81)
--(axis cs:0.77725,81)
--(axis cs:0.77725,0)
--cycle;

\path [draw=white!29.80392156862745!black, fill=color1, opacity=.65] (axis cs:1.00225,0)
--(axis cs:1.22275,0)
--(axis cs:1.22275,136)
--(axis cs:1.00225,136)
--(axis cs:1.00225,0)
--cycle;

\path [draw=white!29.80392156862745!black, fill=color0, opacity=.65] (axis cs:1.77725,3.25)
--(axis cs:1.99775,3.25)
--(axis cs:1.99775,122.5)
--(axis cs:1.77725,122.5)
--(axis cs:1.77725,3.25)
--cycle;

\path [draw=white!29.80392156862745!black, fill=color1, opacity=.65] (axis cs:2.00225,5)
--(axis cs:2.22275,5)
--(axis cs:2.22275,289.25)
--(axis cs:2.00225,289.25)
--(axis cs:2.00225,5)
--cycle;

\path [draw=white!29.80392156862745!black, fill=color0, opacity=.65] (axis cs:2.77725,30.5)
--(axis cs:2.99775,30.5)
--(axis cs:2.99775,162.5)
--(axis cs:2.77725,162.5)
--(axis cs:2.77725,30.5)
--cycle;

\path [draw=white!29.80392156862745!black, fill=color1, opacity=.65] (axis cs:3.00225,31.5)
--(axis cs:3.22275,31.5)
--(axis cs:3.22275,659)
--(axis cs:3.00225,659)
--(axis cs:3.00225,31.5)
--cycle;

\addplot [line width=0.4800000000000001pt, white!29.80392156862745!black, forget plot, dashed, thick]
table [row sep=\\]{%
-0.22275	1 \\
-0.00225	1 \\
};
\addplot [color0, mark=diamond*, mark size=1, mark options={solid,draw=white!29.80392156862745!black,fill=color0}, opacity=.3, only marks, forget plot]
table [row sep=\\]{%
-0.1125	133 \\
-0.1125	144 \\
-0.1125	65 \\
-0.1125	66 \\
-0.1125	51 \\
-0.1125	50 \\
-0.1125	98 \\
-0.1125	71 \\
-0.1125	174 \\
-0.1125	110 \\
-0.1125	89 \\
-0.1125	104 \\
-0.1125	81 \\
-0.1125	81 \\
-0.1125	1066 \\
-0.1125	1253 \\
-0.1125	2447 \\
-0.1125	2896 \\
-0.1125	50 \\
-0.1125	83 \\
-0.1125	62 \\
-0.1125	59 \\
-0.1125	61 \\
-0.1125	550 \\
-0.1125	765 \\
-0.1125	202 \\
-0.1125	223 \\
-0.1125	742 \\
-0.1125	6168 \\
-0.1125	4953 \\
-0.1125	954 \\
-0.1125	1216 \\
-0.1125	4404 \\
-0.1125	93 \\
-0.1125	64 \\
-0.1125	1095 \\
-0.1125	154 \\
-0.1125	365 \\
-0.1125	307 \\
-0.1125	242 \\
-0.1125	284 \\
-0.1125	366 \\
-0.1125	275 \\
-0.1125	372 \\
-0.1125	1339 \\
-0.1125	598 \\
-0.1125	92 \\
-0.1125	64 \\
-0.1125	94 \\
-0.1125	103 \\
-0.1125	91 \\
-0.1125	179 \\
-0.1125	206 \\
-0.1125	61 \\
-0.1125	80 \\
-0.1125	79 \\
-0.1125	125 \\
-0.1125	47 \\
-0.1125	52 \\
-0.1125	100 \\
-0.1125	57 \\
-0.1125	251 \\
-0.1125	754 \\
-0.1125	584 \\
-0.1125	115 \\
-0.1125	62 \\
-0.1125	200 \\
};
\addplot [line width=0.4800000000000001pt, white!29.80392156862745!black, forget plot, dashed, thick]
table [row sep=\\]{%
0.00225	13 \\
0.22275	13 \\
};
\addplot [color1, mark=diamond*, mark size=1, mark options={solid,draw=white!29.80392156862745!black,fill=color1}, opacity=.3, only marks, forget plot]
table [row sep=\\]{%
0.1125	137 \\
0.1125	152 \\
0.1125	224 \\
0.1125	100 \\
0.1125	97 \\
0.1125	98 \\
0.1125	164 \\
0.1125	128 \\
0.1125	128 \\
0.1125	163 \\
0.1125	141 \\
0.1125	132 \\
0.1125	178 \\
0.1125	213 \\
0.1125	332 \\
0.1125	208 \\
0.1125	968 \\
0.1125	1009 \\
0.1125	1034 \\
0.1125	963 \\
0.1125	116 \\
0.1125	328 \\
0.1125	391 \\
0.1125	211 \\
0.1125	156 \\
0.1125	1002 \\
0.1125	2700 \\
0.1125	1231 \\
0.1125	3309 \\
0.1125	150 \\
0.1125	2064 \\
0.1125	315 \\
0.1125	385 \\
0.1125	474 \\
0.1125	287 \\
0.1125	313 \\
0.1125	9266 \\
0.1125	116 \\
0.1125	93 \\
0.1125	140 \\
0.1125	106 \\
0.1125	93 \\
0.1125	309 \\
0.1125	319 \\
0.1125	103 \\
0.1125	142 \\
0.1125	521 \\
0.1125	646 \\
0.1125	672 \\
0.1125	526 \\
0.1125	332 \\
0.1125	312 \\
0.1125	343 \\
0.1125	228 \\
0.1125	100 \\
0.1125	135 \\
0.1125	126 \\
0.1125	121 \\
0.1125	593 \\
0.1125	733 \\
0.1125	536 \\
0.1125	464 \\
0.1125	118 \\
0.1125	199 \\
0.1125	1724 \\
0.1125	209 \\
0.1125	230 \\
};
\addplot [line width=0.4800000000000001pt, white!29.80392156862745!black, forget plot, dashed, thick]
table [row sep=\\]{%
0.77725	15 \\
0.99775	15 \\
};
\addplot [color0, mark=diamond*, mark size=1, mark options={solid,draw=white!29.80392156862745!black,fill=color0}, opacity=.3, only marks, forget plot]
table [row sep=\\]{%
0.8875	98 \\
0.8875	174 \\
0.8875	110 \\
0.8875	89 \\
0.8875	104 \\
0.8875	1066 \\
0.8875	1253 \\
0.8875	2447 \\
0.8875	2896 \\
0.8875	83 \\
0.8875	550 \\
0.8875	765 \\
0.8875	202 \\
0.8875	223 \\
0.8875	742 \\
0.8875	6168 \\
0.8875	4953 \\
0.8875	954 \\
0.8875	1216 \\
0.8875	4404 \\
0.8875	93 \\
0.8875	1095 \\
0.8875	154 \\
0.8875	365 \\
0.8875	307 \\
0.8875	242 \\
0.8875	284 \\
0.8875	366 \\
0.8875	275 \\
0.8875	372 \\
0.8875	1339 \\
0.8875	598 \\
0.8875	92 \\
0.8875	94 \\
0.8875	103 \\
0.8875	91 \\
0.8875	179 \\
0.8875	206 \\
0.8875	125 \\
0.8875	100 \\
0.8875	251 \\
0.8875	754 \\
0.8875	584 \\
0.8875	115 \\
0.8875	200 \\
};
\addplot [line width=0.4800000000000001pt, white!29.80392156862745!black, forget plot, dashed, thick]
table [row sep=\\]{%
1.00225	25 \\
1.22275	25 \\
};
\addplot [color1, mark=diamond*, mark size=1, mark options={solid,draw=white!29.80392156862745!black,fill=color1}, opacity=.3, only marks, forget plot]
table [row sep=\\]{%
1.1125	137 \\
1.1125	152 \\
1.1125	224 \\
1.1125	164 \\
1.1125	163 \\
1.1125	141 \\
1.1125	178 \\
1.1125	213 \\
1.1125	332 \\
1.1125	208 \\
1.1125	328 \\
1.1125	391 \\
1.1125	211 \\
1.1125	156 \\
1.1125	1002 \\
1.1125	2700 \\
1.1125	1231 \\
1.1125	3309 \\
1.1125	150 \\
1.1125	2064 \\
1.1125	315 \\
1.1125	385 \\
1.1125	474 \\
1.1125	287 \\
1.1125	313 \\
1.1125	9266 \\
1.1125	140 \\
1.1125	309 \\
1.1125	319 \\
1.1125	142 \\
1.1125	521 \\
1.1125	646 \\
1.1125	672 \\
1.1125	526 \\
1.1125	332 \\
1.1125	312 \\
1.1125	343 \\
1.1125	228 \\
1.1125	593 \\
1.1125	733 \\
1.1125	536 \\
1.1125	464 \\
1.1125	199 \\
1.1125	1724 \\
1.1125	209 \\
1.1125	230 \\
};
\addplot [line width=0.4800000000000001pt, white!29.80392156862745!black, forget plot, dashed, thick]
table [row sep=\\]{%
1.77725	38 \\
1.99775	38 \\
};
\addplot [color0, mark=diamond*, mark size=1, mark options={solid,draw=white!29.80392156862745!black,fill=color0}, opacity=.3, only marks, forget plot]
table [row sep=\\]{%
1.8875	0 \\
1.8875	0 \\
1.8875	1 \\
1.8875	1 \\
1.8875	3 \\
1.8875	1 \\
1.8875	0 \\
1.8875	0 \\
1.8875	0 \\
1.8875	0 \\
1.8875	0 \\
1.8875	0 \\
1.8875	0 \\
1.8875	1 \\
1.8875	0 \\
1.8875	2 \\
1.8875	1 \\
1.8875	0 \\
1.8875	0 \\
1.8875	0 \\
1.8875	0 \\
1.8875	1066 \\
1.8875	1253 \\
1.8875	2447 \\
1.8875	2896 \\
1.8875	202 \\
1.8875	742 \\
1.8875	6168 \\
1.8875	4953 \\
1.8875	954 \\
1.8875	1216 \\
1.8875	4404 \\
1.8875	1095 \\
1.8875	154 \\
1.8875	1339 \\
1.8875	179 \\
1.8875	206 \\
1.8875	125 \\
1.8875	251 \\
1.8875	754 \\
1.8875	584 \\
1.8875	200 \\
};
\addplot [line width=0.4800000000000001pt, white!29.80392156862745!black, forget plot, dashed, thick]
table [row sep=\\]{%
2.00225	48.5 \\
2.22275	48.5 \\
};
\addplot [color1, mark=diamond*, mark size=1, mark options={solid,draw=white!29.80392156862745!black,fill=color1}, opacity=.3, only marks, forget plot]
table [row sep=\\]{%
2.1125	0 \\
2.1125	0 \\
2.1125	0 \\
2.1125	3 \\
2.1125	0 \\
2.1125	0 \\
2.1125	0 \\
2.1125	0 \\
2.1125	0 \\
2.1125	0 \\
2.1125	0 \\
2.1125	0 \\
2.1125	4 \\
2.1125	0 \\
2.1125	3 \\
2.1125	2 \\
2.1125	0 \\
2.1125	2 \\
2.1125	0 \\
2.1125	332 \\
2.1125	1002 \\
2.1125	2700 \\
2.1125	1231 \\
2.1125	3309 \\
2.1125	2064 \\
2.1125	9266 \\
2.1125	309 \\
2.1125	319 \\
2.1125	521 \\
2.1125	646 \\
2.1125	672 \\
2.1125	526 \\
2.1125	332 \\
2.1125	312 \\
2.1125	343 \\
2.1125	593 \\
2.1125	733 \\
2.1125	536 \\
2.1125	464 \\
2.1125	1724 \\
};
\addplot [line width=0.4800000000000001pt, white!29.80392156862745!black, forget plot, dashed, thick]
table [row sep=\\]{%
2.77725	52 \\
2.99775	52 \\
};
\addplot [color0, mark=diamond*, mark size=1, mark options={solid,draw=white!29.80392156862745!black,fill=color0}, opacity=.3, only marks, forget plot]
table [row sep=\\]{%
2.8875	1 \\
2.8875	0 \\
2.8875	2 \\
2.8875	1 \\
2.8875	30 \\
2.8875	0 \\
2.8875	24 \\
2.8875	6168 \\
2.8875	4953 \\
2.8875	4404 \\
2.8875	1095 \\
2.8875	1339 \\
2.8875	754 \\
2.8875	200 \\
};
\addplot [line width=0.4800000000000001pt, white!29.80392156862745!black, forget plot, dashed, thick]
table [row sep=\\]{%
3.00225	228 \\
3.22275	228 \\
};
\addplot [color1, mark=diamond*, mark size=1, mark options={solid,draw=white!29.80392156862745!black,fill=color1}, opacity=.3, only marks, forget plot]
table [row sep=\\]{%
3.1125	29 \\
3.1125	0 \\
3.1125	3 \\
3.1125	2 \\
3.1125	0 \\
3.1125	9 \\
3.1125	2 \\
3.1125	2700 \\
3.1125	1231 \\
3.1125	3309 \\
3.1125	2064 \\
3.1125	9266 \\
3.1125	672 \\
3.1125	733 \\
};

\nextgroupplot[
axis line style={white!80.0!black},
height=\figureheight,
tick pos=both,
width=\figurewidth,
x grid style={white!80.0!black},
xmin=-0.5, xmax=3.5,
xtick={0,1,2,3},
xticklabels={\affectedsols,$\bracket{10}{tilim}$,$\bracket{100}{tilim}$,$\bracket{1000}{tilim}$},
xlabel={Instance Group},
y grid style={white!80.0!black},
ymajorgrids,
ymin=-0.5, ymax=60
]

\path [draw=white!29.80392156862745!black, fill=color0, opacity=.65] (axis cs:-0.22275,0)
--(axis cs:-0.00225,0)
--(axis cs:-0.00225,46)
--(axis cs:-0.22275,46)
--(axis cs:-0.22275,0)
--cycle;

\path [draw=white!29.80392156862745!black, fill=color1, opacity=.65] (axis cs:0.00225,0)
--(axis cs:0.22275,0)
--(axis cs:0.22275,89)
--(axis cs:0.00225,89)
--(axis cs:0.00225,0)
--cycle;

\path [draw=white!29.80392156862745!black, fill=color0, opacity=.65] (axis cs:0.77725,0)
--(axis cs:0.99775,0)
--(axis cs:0.99775,81)
--(axis cs:0.77725,81)
--(axis cs:0.77725,0)
--cycle;

\path [draw=white!29.80392156862745!black, fill=color1, opacity=.65] (axis cs:1.00225,0)
--(axis cs:1.22275,0)
--(axis cs:1.22275,136)
--(axis cs:1.00225,136)
--(axis cs:1.00225,0)
--cycle;

\path [draw=white!29.80392156862745!black, fill=color0, opacity=.65] (axis cs:1.77725,3.25)
--(axis cs:1.99775,3.25)
--(axis cs:1.99775,122.5)
--(axis cs:1.77725,122.5)
--(axis cs:1.77725,3.25)
--cycle;

\path [draw=white!29.80392156862745!black, fill=color1, opacity=.65] (axis cs:2.00225,5)
--(axis cs:2.22275,5)
--(axis cs:2.22275,289.25)
--(axis cs:2.00225,289.25)
--(axis cs:2.00225,5)
--cycle;

\path [draw=white!29.80392156862745!black, fill=color0, opacity=.65] (axis cs:2.77725,30.5)
--(axis cs:2.99775,30.5)
--(axis cs:2.99775,162.5)
--(axis cs:2.77725,162.5)
--(axis cs:2.77725,30.5)
--cycle;

\path [draw=white!29.80392156862745!black, fill=color1, opacity=.65] (axis cs:3.00225,31.5)
--(axis cs:3.22275,31.5)
--(axis cs:3.22275,659)
--(axis cs:3.00225,659)
--(axis cs:3.00225,31.5)
--cycle;

\draw[draw=white!29.80392156862745!black,fill=color0] (axis cs:0,0) rectangle (axis cs:0,0);
\draw[draw=white!29.80392156862745!black,fill=color1] (axis cs:0,0) rectangle (axis cs:0,0);

\addplot [line width=0.4800000000000001pt, white!29.80392156862745!black, forget plot, dashed, thick]
table [row sep=\\]{%
-0.22275	1 \\
-0.00225	1 \\
};
\addplot [color0, mark=diamond*, mark size=1, mark options={solid,draw=white!29.80392156862745!black,fill=color0}, opacity=.3, only marks, forget plot]
table [row sep=\\]{%
-0.1125	133 \\
-0.1125	144 \\
-0.1125	65 \\
-0.1125	66 \\
-0.1125	51 \\
-0.1125	50 \\
-0.1125	98 \\
-0.1125	71 \\
-0.1125	174 \\
-0.1125	110 \\
-0.1125	89 \\
-0.1125	104 \\
-0.1125	81 \\
-0.1125	81 \\
-0.1125	1066 \\
-0.1125	1253 \\
-0.1125	2447 \\
-0.1125	2896 \\
-0.1125	50 \\
-0.1125	83 \\
-0.1125	62 \\
-0.1125	59 \\
-0.1125	61 \\
-0.1125	550 \\
-0.1125	765 \\
-0.1125	202 \\
-0.1125	223 \\
-0.1125	742 \\
-0.1125	6168 \\
-0.1125	4953 \\
-0.1125	954 \\
-0.1125	1216 \\
-0.1125	4404 \\
-0.1125	93 \\
-0.1125	64 \\
-0.1125	1095 \\
-0.1125	154 \\
-0.1125	365 \\
-0.1125	307 \\
-0.1125	242 \\
-0.1125	284 \\
-0.1125	366 \\
-0.1125	275 \\
-0.1125	372 \\
-0.1125	1339 \\
-0.1125	598 \\
-0.1125	92 \\
-0.1125	64 \\
-0.1125	94 \\
-0.1125	103 \\
-0.1125	91 \\
-0.1125	179 \\
-0.1125	206 \\
-0.1125	61 \\
-0.1125	80 \\
-0.1125	79 \\
-0.1125	125 \\
-0.1125	47 \\
-0.1125	52 \\
-0.1125	100 \\
-0.1125	57 \\
-0.1125	251 \\
-0.1125	754 \\
-0.1125	584 \\
-0.1125	115 \\
-0.1125	62 \\
-0.1125	200 \\
};
\addplot [line width=0.4800000000000001pt, white!29.80392156862745!black, forget plot, dashed, thick]
table [row sep=\\]{%
0.00225	13 \\
0.22275	13 \\
};
\addplot [color1, mark=diamond*, mark size=1, mark options={solid,draw=white!29.80392156862745!black,fill=color1}, opacity=.3, only marks, forget plot]
table [row sep=\\]{%
0.1125	137 \\
0.1125	152 \\
0.1125	224 \\
0.1125	100 \\
0.1125	97 \\
0.1125	98 \\
0.1125	164 \\
0.1125	128 \\
0.1125	128 \\
0.1125	163 \\
0.1125	141 \\
0.1125	132 \\
0.1125	178 \\
0.1125	213 \\
0.1125	332 \\
0.1125	208 \\
0.1125	968 \\
0.1125	1009 \\
0.1125	1034 \\
0.1125	963 \\
0.1125	116 \\
0.1125	328 \\
0.1125	391 \\
0.1125	211 \\
0.1125	156 \\
0.1125	1002 \\
0.1125	2700 \\
0.1125	1231 \\
0.1125	3309 \\
0.1125	150 \\
0.1125	2064 \\
0.1125	315 \\
0.1125	385 \\
0.1125	474 \\
0.1125	287 \\
0.1125	313 \\
0.1125	9266 \\
0.1125	116 \\
0.1125	93 \\
0.1125	140 \\
0.1125	106 \\
0.1125	93 \\
0.1125	309 \\
0.1125	319 \\
0.1125	103 \\
0.1125	142 \\
0.1125	521 \\
0.1125	646 \\
0.1125	672 \\
0.1125	526 \\
0.1125	332 \\
0.1125	312 \\
0.1125	343 \\
0.1125	228 \\
0.1125	100 \\
0.1125	135 \\
0.1125	126 \\
0.1125	121 \\
0.1125	593 \\
0.1125	733 \\
0.1125	536 \\
0.1125	464 \\
0.1125	118 \\
0.1125	199 \\
0.1125	1724 \\
0.1125	209 \\
0.1125	230 \\
};
\addplot [line width=0.4800000000000001pt, white!29.80392156862745!black, forget plot, dashed, thick]
table [row sep=\\]{%
0.77725	15 \\
0.99775	15 \\
};
\addplot [color0, mark=diamond*, mark size=1, mark options={solid,draw=white!29.80392156862745!black,fill=color0}, opacity=.3, only marks, forget plot]
table [row sep=\\]{%
0.8875	98 \\
0.8875	174 \\
0.8875	110 \\
0.8875	89 \\
0.8875	104 \\
0.8875	1066 \\
0.8875	1253 \\
0.8875	2447 \\
0.8875	2896 \\
0.8875	83 \\
0.8875	550 \\
0.8875	765 \\
0.8875	202 \\
0.8875	223 \\
0.8875	742 \\
0.8875	6168 \\
0.8875	4953 \\
0.8875	954 \\
0.8875	1216 \\
0.8875	4404 \\
0.8875	93 \\
0.8875	1095 \\
0.8875	154 \\
0.8875	365 \\
0.8875	307 \\
0.8875	242 \\
0.8875	284 \\
0.8875	366 \\
0.8875	275 \\
0.8875	372 \\
0.8875	1339 \\
0.8875	598 \\
0.8875	92 \\
0.8875	94 \\
0.8875	103 \\
0.8875	91 \\
0.8875	179 \\
0.8875	206 \\
0.8875	125 \\
0.8875	100 \\
0.8875	251 \\
0.8875	754 \\
0.8875	584 \\
0.8875	115 \\
0.8875	200 \\
};
\addplot [line width=0.4800000000000001pt, white!29.80392156862745!black, forget plot, dashed, thick]
table [row sep=\\]{%
1.00225	25 \\
1.22275	25 \\
};
\addplot [color1, mark=diamond*, mark size=1, mark options={solid,draw=white!29.80392156862745!black,fill=color1}, opacity=.3, only marks, forget plot]
table [row sep=\\]{%
1.1125	137 \\
1.1125	152 \\
1.1125	224 \\
1.1125	164 \\
1.1125	163 \\
1.1125	141 \\
1.1125	178 \\
1.1125	213 \\
1.1125	332 \\
1.1125	208 \\
1.1125	328 \\
1.1125	391 \\
1.1125	211 \\
1.1125	156 \\
1.1125	1002 \\
1.1125	2700 \\
1.1125	1231 \\
1.1125	3309 \\
1.1125	150 \\
1.1125	2064 \\
1.1125	315 \\
1.1125	385 \\
1.1125	474 \\
1.1125	287 \\
1.1125	313 \\
1.1125	9266 \\
1.1125	140 \\
1.1125	309 \\
1.1125	319 \\
1.1125	142 \\
1.1125	521 \\
1.1125	646 \\
1.1125	672 \\
1.1125	526 \\
1.1125	332 \\
1.1125	312 \\
1.1125	343 \\
1.1125	228 \\
1.1125	593 \\
1.1125	733 \\
1.1125	536 \\
1.1125	464 \\
1.1125	199 \\
1.1125	1724 \\
1.1125	209 \\
1.1125	230 \\
};
\addplot [line width=0.4800000000000001pt, white!29.80392156862745!black, forget plot, dashed, thick]
table [row sep=\\]{%
1.77725	38 \\
1.99775	38 \\
};
\addplot [color0, mark=diamond*, mark size=1, mark options={solid,draw=white!29.80392156862745!black,fill=color0}, opacity=.3, only marks, forget plot]
table [row sep=\\]{%
1.8875	0 \\
1.8875	0 \\
1.8875	1 \\
1.8875	1 \\
1.8875	3 \\
1.8875	1 \\
1.8875	0 \\
1.8875	0 \\
1.8875	0 \\
1.8875	0 \\
1.8875	0 \\
1.8875	0 \\
1.8875	0 \\
1.8875	1 \\
1.8875	0 \\
1.8875	2 \\
1.8875	1 \\
1.8875	0 \\
1.8875	0 \\
1.8875	0 \\
1.8875	0 \\
1.8875	1066 \\
1.8875	1253 \\
1.8875	2447 \\
1.8875	2896 \\
1.8875	202 \\
1.8875	742 \\
1.8875	6168 \\
1.8875	4953 \\
1.8875	954 \\
1.8875	1216 \\
1.8875	4404 \\
1.8875	1095 \\
1.8875	154 \\
1.8875	1339 \\
1.8875	179 \\
1.8875	206 \\
1.8875	125 \\
1.8875	251 \\
1.8875	754 \\
1.8875	584 \\
1.8875	200 \\
};
\addplot [line width=0.4800000000000001pt, white!29.80392156862745!black, forget plot, dashed, thick]
table [row sep=\\]{%
2.00225	48.5 \\
2.22275	48.5 \\
};
\addplot [color1, mark=diamond*, mark size=1, mark options={solid,draw=white!29.80392156862745!black,fill=color1}, opacity=.3, only marks, forget plot]
table [row sep=\\]{%
2.1125	0 \\
2.1125	0 \\
2.1125	0 \\
2.1125	3 \\
2.1125	0 \\
2.1125	0 \\
2.1125	0 \\
2.1125	0 \\
2.1125	0 \\
2.1125	0 \\
2.1125	0 \\
2.1125	0 \\
2.1125	4 \\
2.1125	0 \\
2.1125	3 \\
2.1125	2 \\
2.1125	0 \\
2.1125	2 \\
2.1125	0 \\
2.1125	332 \\
2.1125	1002 \\
2.1125	2700 \\
2.1125	1231 \\
2.1125	3309 \\
2.1125	2064 \\
2.1125	9266 \\
2.1125	309 \\
2.1125	319 \\
2.1125	521 \\
2.1125	646 \\
2.1125	672 \\
2.1125	526 \\
2.1125	332 \\
2.1125	312 \\
2.1125	343 \\
2.1125	593 \\
2.1125	733 \\
2.1125	536 \\
2.1125	464 \\
2.1125	1724 \\
};
\addplot [line width=0.4800000000000001pt, white!29.80392156862745!black, forget plot, dashed, thick]
table [row sep=\\]{%
2.77725	52 \\
2.99775	52 \\
};
\addplot [color0, mark=diamond*, mark size=1, mark options={solid,draw=white!29.80392156862745!black,fill=color0}, opacity=.3, only marks, forget plot]
table [row sep=\\]{%
2.8875	1 \\
2.8875	0 \\
2.8875	2 \\
2.8875	1 \\
2.8875	30 \\
2.8875	0 \\
2.8875	24 \\
2.8875	6168 \\
2.8875	4953 \\
2.8875	4404 \\
2.8875	1095 \\
2.8875	1339 \\
2.8875	754 \\
2.8875	200 \\
};
\addplot [line width=0.4800000000000001pt, white!29.80392156862745!black, forget plot, dashed, thick]
table [row sep=\\]{%
3.00225	228 \\
3.22275	228 \\
};
\addplot [color1, mark=diamond*, mark size=1, mark options={solid,draw=white!29.80392156862745!black,fill=color1}, opacity=.3, only marks, forget plot]
table [row sep=\\]{%
3.1125	29 \\
3.1125	0 \\
3.1125	3 \\
3.1125	2 \\
3.1125	0 \\
3.1125	9 \\
3.1125	2 \\
3.1125	2700 \\
3.1125	1231 \\
3.1125	3309 \\
3.1125	2064 \\
3.1125	9266 \\
3.1125	672 \\
3.1125	733 \\
};

\end{groupplot}

\end{tikzpicture}
    \caption{Box plot showing the number of generated conflicts by \farkasdiving when running with \farkdiving
             and the virtual best diving heuristic when running with \default on increasingly hard subgroups of \affectedsols.
             }
    \label{fig:farkasdiving_confs}
\end{figure}
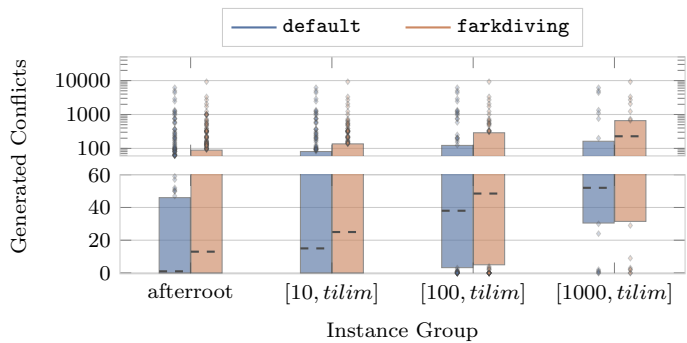

\farkasdiving is motivated by the idea of explicitly diving towards a valid Farkas proof.
To analyze how successful \farkasdiving is in generating conflict constraints we considered all affected instances
where \farkasdiving was allowed to run after the root node,
\ie instances where \farkasdiving was able to find a feasible solution at the root node.
This subset contains $273$ instances and is displayed in line ``\affectedsols'' of Table~\ref{tab:farkasdiving_aggregated}.
On this subset of instances, \farkasdiving was able to find $\percent{7.3}$ more conflict constraint
than the ``virtual best diving heuristic'' in this regard.
Here, the virtual best diving heuristic is determined by taking, for each instance run with \default settings,
the diving heuristic that generated the largest number of conflict constraints.
Hence, \farkasdiving indeed succeeded in generating an above-average number of conflict constraints.
Figure~\ref{fig:farkasdiving_confs} provides a more detailed comparison of the number of generated conflict constraints for increasingly hard subgroups of \affectedsols.
The box plots~\citep{mcgill1978variations} show for every setting and instance group the 1st and 3rd quartile (shaded box)
as well as the median (dashed line).
All observations below the 1st or above the 3rd quartile are marked with shaded diamonds.
On all four groups of instances, \farkasdiving led to more conflict constraints in the median and 3rd quartile
than the virtual best diving heuristic in this regard.
On instances where both settings needed at least $1000$ seconds, \farkasdiving produced on $\percent{50}$ of the instances more
conflict constraints than the virtual best diving heuristic of \scip with \default settings on $\percent{75}$ did.
These results indicate that the strategy of diving towards a valid Farkas proof as performed by \farkasdiving
leads to additional conflict information and succeeds over the whole set where this strategy is pursued.
Note that our analysis is conservative.
The virtual best diving heuristic strictly overestimates every single diving heuristic, hence the number of additionally generated conflict constraints would only increase when comparing to each diving heuristic individually.

\paragraph{Success in Generating Solution.}

Usually, the number of solutions found by diving heuristics is quite small~\citep{Khalil2017}.
From a primal viewpoint, \farkasdiving follows a rounding strategy which leads to overly optimistic solutions
that are expected to be infeasible most of the time.
However, if this strategy succeeds, the solutions can be expected to be very good.
A more detailed look into the success rate of \farkasdiving \wrt finding primal solutions
answers the question about the impact of these overly optimistic solutions on the overall MIP solver.
In our experiments over several seeds, we could observe that on affected instances \farkasdiving was able to find $\percent{39.7}$ more feasible solutions (\sols)
and $\percent{3.7}$ more improving solutions (\bestsols) than the virtual best diving heuristic (\virtualbestdsols and \virtualbestdbsols)
of \scip with \default settings.
Here, the virtual best diving heuristic is determined as above, for each instance selecting the heuristic with the largest number of feasible and improving solutions, respectively.
Again, both \virtualbestdsols and \virtualbestdbsols overestimate every single diving heuristic in \default setting.
Note that this evaluation is conservative also in the sense that it
includes instances for which \farkasdiving was not allowed to run at local nodes of the tree,
\ie instances where \farkasdiving did not find a solution at the root node, and, therefore, did not find a solution at all in these instances.

\begin{figure}
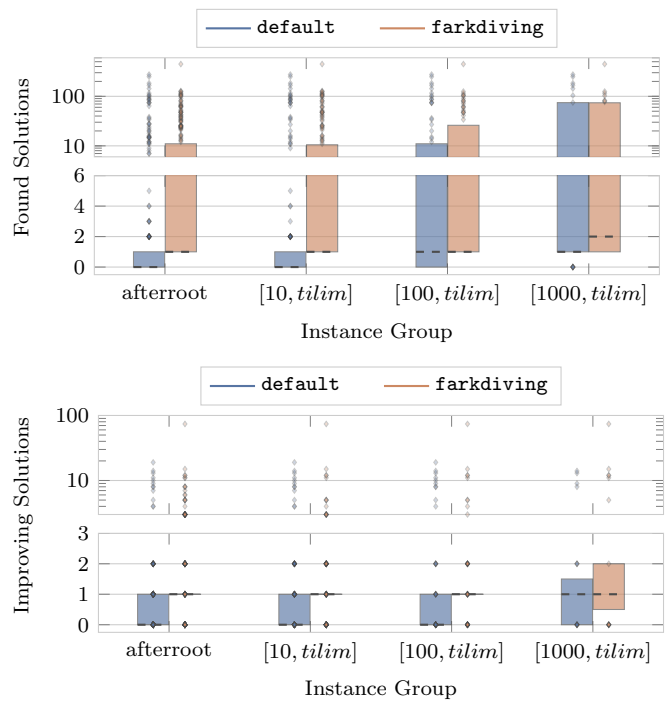

    \centering
    \begin{subfigure}[t][][t]{\ifreport.75\else.5\fi\textwidth}
        \setlength\figureheight{.15\textheight}
        \setlength\figurewidth{\textwidth}

    \end{subfigure}
    \caption{Box plots showing the number of feasible and improving solutions found by \farkasdiving when running with \farkdiving
             and the virtual best diving heuristic when running with \default.
             }
    \label{fig:farkasdiving_sols}
\end{figure}

On the group of instances where \farkasdiving was allowed to run within the tree (\affectedsols), the results are even more pronounced.
\farkasdiving was able to find $\percent{57.2}$ more feasible solutions and more than twice as many improving solutions
than the virtual best diving heuristic in \default setting.
Figure~\ref{fig:farkasdiving_sols} plots the relation of feasible and improving solutions found by \farkasdiving
and the virtual best diving heuristic of \default on increasingly hard groups of \affectedsols.
The box plots show for every setting and instance group the 1st and 3rd quartile of all numbers of feasible and improving solutions, respectively, (shaded box)
as well the median (dashed line).
All observations below the 1st or above the 3rd quartile are marked with shaded diamonds.
On all four groups of instances, \farkasdiving was able to find $10$ or more feasible solution on at least $\percent{25}$
of the instances (observations above the 3rd quartile), whereas the 3rd quartile of the virtual best diving heuristic in this regard
is always $1$ for \affected and \bracket{10}{tilim}.
On hard instances for which both setting need at least $1000$ seconds the 1st and 3rd quartiles of both \default and \farkdiving are identical.
Note, the set \affectedsols only contains instances for which \farkasdiving was able to find at least one feasible solution at the root node.
Thus, this set of instances is slightly biased; hence, it is not suitable to conclude that \farkasdiving finds more feasible solutions
than the virtual best diving heuristic in general.
However, the results on this set of instances indicate that the rule we used to decide whether \farkasdiving is allowed to run within the tree succeeds.
On the complementary set consisting of 1562 instances that could be solved by at least one setting,
\ie those for which \farkasdiving was not able to find a solution root node, \farkdiving succeeded \wrt finding at least one feasible solution within the tree
on only $\percent{0.019}$ of the instances.
Considering the improving solutions found by the virtual best diving heuristic when running \scip with \default settings and \farkasdiving
indicates that whenever \farkasdiving was allowed to run within the tree, \ie it founds a feasible solution at the root node,
\farkasdiving yields at least one improving solution on $\percent{75}$ (all observations above the 1st quartile)
of the instances of \affectedsols, \bracket{10}{tilim}, and \bracket{100}{tilim}.
By contrast, the 1st and 3rd quartile of the virtual best diving heuristic are $0$ and $1$, respectively.
Consequently, we can conclude that on instances that are cumbersome for the established diving heuristics of \scip~6.0 \wrt finding solutions,
\farkasdiving can easily find both feasible and improving solutions if it is allowed to run within the tree.

However, the pure amount of feasible and improving solutions alone is not meaningful enough to conclude
whether the solutions found by \farkasdiving have a positive impact on the overall MIP solver.
Therefore, we consider also the \emph{primal integral}~\citep{Berthold2013},
which measures the progress of the primal bound towards the optimal solution.
Compared to \default, \farkdiving improved the primal integral by $\percent{12.3}$.
This exhibits, maybe surprisingly so, that the optimistic strategy of \farkasdiving is also very successful on the primal side
and, thus, a valuable extension to the portfolio of primal heuristics.

\medskip

\paragraph{Impact of Primal Solutions and Generated Conflicts.}

The previous results indicate that \farkasdiving both produces an over-average number of solutions and conflict constraints if it is allowed to run within the tree,
but do not finally answer which component is more relevant for performance.
We tried to quantify this by testing two modified versions of \farkasdiving:
one that disables conflict analysis (\farkdivingnoconfs) and one that discards feasible solutions found by \farkasdiving (\farkdivingnosols).
The aggregated results are reported in Table~\ref{tab:farkasdiving_no_sols_no_confs}.
For this experiment, we used the subset of $273$~instances where \farkasdiving was allowed to be executed after the root node when running \farkasdiving, \ie \affectedsols.

\begin{table*}%
\begin{centering}
\footnotesize
\setlength{\tabcolsep}{.75pt}
\caption{\farkdiving with and without adding feasible solutions and generated conflict constraints, respectively.}\label{tab:farkasdiving_no_sols_no_confs}
\begin{tabularx}{\textwidth}{L *{6}{R}}
    \toprule
    & \# & \solved & \time & \nodes & \timeQ & \nodesQ \\
    \midrule
    \default                   & 273  &  270  &  33.82  &  1095  &  1.000 & 1.000 \\
    \mbox{\farkdiving}         & 273  &  272  &  29.75  &   857  &  0.879 & 0.783 \\
    \mbox{\farkdivingnoconfs}  & 273  &  266  &  31.79  &   976  &  0.940 & 0.891 \\
    \mbox{\farkdivingnosols}   & 273  &  269  &  32.73  &  1045  &  0.968 & 0.954 \\
    \bottomrule\\
\end{tabularx}
\end{centering}
\end{table*}
 
Both settings solve fewer instances than \default.
This result is not surprising since both \farkdivingnoconfs and \farkdivingnosols spend computational time for running \farkasdiving.
As we have already discussed, \farkasdiving is much more expensive than other diving heuristics in \scip since it solves the LP at every node during diving.
Consequently, by not performing conflict analysis or not adding feasible solutions during \farkasdiving one of the key ingredients is missing.
Thus, on a few instances, \farkasdiving needs the impact of both found solutions and generated conflict constraints
to compensate the computational overhead compared to \default.

When disabling conflict analysis (\farkdivingnoconfs) the performance improvement compared to \default
decreased from $\percent{12.1}$ to $\percent{6}$.
The number of nodes explored during the tree search also increased (compared to \farkdiving) but is still $\percent{10.9}$
smaller than \scip with \default settings.
Disabling the addition of primal solutions found by \farkasdiving to the main search (\farkdivingnosols)
reduced the performance improvement compared to \default from $\percent{12.1}$ to $\percent{3.2}$,
which corresponds to a slowdown by $\percent{9.2}$ compared to \farkdiving.
Interestingly, the tree increased by $\percent{17.9}$ compared to \farkdiving.
However, \farkdivingnosols still led to $\percent{4.6}$ smaller trees than \scip with \default settings.

In hindsight, it is not very surprising that disabling the addition of solutions
has a larger impact than disabling conflict analysis during \farkasdiving.
With the very optimistic rounding strategy, \farkasdiving was able to find $3.2$~times as many best solutions%
\footnote{A solution is called 'best' if it is an optimal solution or the best-known solution when reaching the time limit.}
as all remaining diving heuristics together.
Overall, on \affectedsols, $\percent{5.7}$ of all best solutions were found by \farkasdiving.
The only heuristic that was even more effective on this set of instances is the large neighborhood search heuristic \textsc{Rens}~\citep{Berthold2014rens},
which contributed $\percent{7.9}$ of all best solutions.
Thus, \farkasdiving has a remarkable success rate, which leads to a not negligible impact on the primal side.
As mentioned above, this is also reflected by the primal integral~\citep{Berthold2013}, which improved by $\percent{12.3}$
when running \farkasdiving during the tree search (\affectedsols).

These results make clear that both the primal and dual aspects of \farkasdiving contributes to the improved performance,
but that the solutions found with the strategy of \farkasdiving seem to be the main driver of the heuristic.
This may come as a surprise because our initial motivation for the design of \farkasdiving was the targeted generation of conflict constraints.

\medskip

\paragraph{Impact of LP Frequency.}

By default, all diving heuristics in \scip are configured to solve the LP dependent on the number of found bound deductions
during constraint propagation over all variables.
An LP solve is triggered whenever the number of bound changes since the last LP solve exceeds $0.15$~times the number of variables.
By contrast, \farkasdiving solves the LP at every diving node.
One motivation for this high LP frequency is to ``pull'' the variable assignments back towards the feasible region
in order to counteract the very optimistic rounding strategy of \farkasdiving,
which tends to ``push'' the variable assignments out of the feasible region.

Within \scip, every LP solution found during diving is automatically rounded to a solution satisfying all integrality requirements of the variables.
If this solution also satisfies all constraints, \ie is feasible for the entire MIP, the solution is added to the solution storage,
whereby the diving heuristic that performs the current dive is credited for finding the feasible solution.
Thus, \farkasdiving may have a slight advantage over all other diving heuristics since it solves the LP relaxation more frequently.
This fact might be one reason why \farkasdiving was able to find more than $\percent{40}$ more solutions than the virtual best (\virtualbestdsols), see Table~\ref{tab:farkasdiving_aggregated}.
Hence, in a final control experiment, we configured \farkasdiving identically to all remaining diving heuristics in order to measure the impact of the increased LP frequency.
We will refer to this configuration by \farkdivinglp.

On the set of affected instances for which our criteria for running \farkasdiving within the tree is satisfied (\affectedsols),
\farkdivinglp performed almost as ``poorly'' as \farkdivingnosols (see Table~\ref{tab:farkasdiving_no_sols_no_confs}) \wrt solving time, nodes, and number of solved instances.
Thus, when using \farkdivinglp the solving time could be improved by only $\percent{3.5}$ compared to \default.
Note, on the same set of instances \farkdiving led to a speed-up of $\percent{12.1}$.
Compared to \farkdiving, the amount of feasible solutions found by \farkdivinglp was reduced by $\percent{90}$.
These results show that solving the LP relaxation frequently gives an important boost to the degree by \emph{how much} \farkasdiving improves performance.
However, also \farkasdiving with less LP solves, \ie configured identically to all remaining diving heuristics, %
outperforms the \default settings.
Consequently, the fact \emph{that} \farkasdiving works indeed seems to be a result of its specific choice of rounding function $\roundfunc{F}$ and scoring function $\scorefunc{F}$.

\subsection{Conflict Diving}
\label{subsec:compresults-conflict}

The new \conflictdiving is closely related to the existing \coefficientdiving in the sense that they both exploit lock information
and follow the same diving framework of Algorithm~\ref{alg:divingalgo}.
Hence, we chose to include both in the computational analysis and compare their impact individually to \scip without any of these lock-exploiting diving heuristics.
In the following, we will refer to the latter setting by \nolockdiving.
We will refer to \scip with \conflictdiving and with \coefficientdiving activated by \confdiving and \coefdiving, respectively.
We first present computational experiments to quantify the general performance impact, followed by
further computational results to analyze the impact of individual components in more detail.

\paragraph{Overall Impact of Conflict Diving.}

In our computational setup, \conflictdiving used the same parameter settings as \coefficientdiving, \eg frequency of execution and LP solve frequency.
The parameter to weight variable and conflict locks within \conflictdiving was set to $\confweight = 0.75$, \ie conflict locks dominate by a factor of $3$.
Aggregated computational results of all three configurations on \MMMC are shown in Table~\ref{tab:conflictdiving_aggregated}.

While \coefdiving only showed minimally improved performance compared to \nolockdiving,
the setting \confdiving was clearly superior. %
Conflict diving increased the number of solved instances by 16 (\solved), led to an overall speed-up of $\percent{4.8}$ (\timeQ),
and reduced the tree size by $\percent{5.5}$ (\nodesQ) on affected instances.
On the subset of affected instances, \confdiving needed $\percent{4.6}$ less solving time and
up to $\percent{6.8}$ fewer nodes compared to \coefdiving.
The node reduction on the affected instances may be explained by the increased number of generated conflict constraints (\confs).

Concerning the number of solutions, %
we observed that \confdiving was only slightly more successful:
\confdiving found $\percent{11.1}$ more feasible solutions and $\percent{4.7}$ more improving solutions than \coefdiving.
The impact of these additional solutions is reflected by the primal integral, which decreased by $\percent{4.7}$ compared to \coefdiving
and by~$\percent{4.9}$ compared to \nolockdiving.
However, both heuristics had a small success rate \wrt finding solutions: $0.037$ (\conflictdiving) and $0.035$ (\coefficientdiving) solutions per dive.

\begin{table*}%
\begin{centering}
\ifreport
    \scriptsize
\else
    \footnotesize
\fi
\setlength{\tabcolsep}{.75pt}
\caption{Aggregated computational results for coefficient and conflict diving on \MMMC over four different random seeds.
        \ifhighlighted Relative changes by at least $5\%$ are highlighted in bold. %
        \fi}\label{tab:conflictdiving_aggregated}
\begin{tabularx}{\textwidth}{L@{\hskip .2cm}R *{11}{R}}
    \toprule
    &  & \multicolumn{3}{c}{\nolockdiving} & \multicolumn{4}{c}{\coefdiving} & \multicolumn{4}{c}{\confdiving}\\
    \cmidrule(lr){3-5}  \cmidrule(lr){6-9}  \cmidrule(lr){10-13}
    & \# & \solved & \time & \nodes & \solved & \timeQ & \nodesQ & \confs & \solved & \timeQ & \nodesQ & \confs \\
    \midrule
    \cleaninst                & 1941 &   1372 &       216 &      3108 &   1381  &              0.998   &              1.003   &    1386.4  &   1388  &              0.978   &              0.972   &    1364.4  \\
    \miplib~2010              &  347 &    297 &       425 &      5788 &    299  &              1.008   &              1.015   &     374.9  &    300  &              0.986   &              0.975   &     343.6  \\
    \midrule
    \affected                 &  861 &    822 &       186 &      4259 &    831  &              0.998   &              1.014   &     686.2  &    838  &              0.952   & \better{     0.945}  &     892.2  \\
    \bracket{10}{tilim}       &  775 &    736 &       287 &      5356 &    745  &              0.998   &              1.018   &     760.9  &    752  & \better{     0.948}  & \better{     0.941}  &     990.0  \\
    \bracket{100}{tilim}      &  542 &    503 &       707 &      7982 &    512  &              0.988   &              1.000   &    1045.2  &    519  & \better{     0.932}  & \better{     0.923}  &    1377.9  \\
    \bracket{1000}{tilim}     &  246 &    207 &      2365 &     19538 &    216  &              0.984   &              0.981   &    2068.5  &    223  & \better{     0.874}  & \better{     0.882}  &    2820.5  \\
    \bottomrule\\
\end{tabularx}
\end{centering}
\end{table*}
 
\begin{figure}
  \begin{center}%
    \begin{subfigure}[t][][t]{\ifreport.5\else.25\fi\textwidth}
        {
        \setlength\figureheight{.3\textheight}
        \setlength\figurewidth{\textwidth}
\begin{tikzpicture}

\definecolor{color0}{rgb}{0.298039215686275,0.447058823529412,0.690196078431373}
\definecolor{color1}{rgb}{0.866666666666667,0.517647058823529,0.32156862745098}

\begin{axis}[
axis line style={white!80.0!black},
height=\figureheight,
legend cell align={left},
legend entries={{\coefdiving},{\confdiving}},
legend style={at={\ifreport(1.225,1.1)\else(3,1.10)\fi}, anchor=south, draw=white!80.0!black, font=\footnotesize, /tikz/every even column/.append style={column sep=0.5cm}},
legend columns=2,
tick pos=both,
width=\figurewidth,
x grid style={white!80.0!black},
xlabel={Time factor},
xmajorgrids,
xmin=1, xmax=10,
y grid style={white!80.0!black},
ymajorgrids,
ymin=0.3, ymax=1,
ytick={0.2,0.3,0.4,0.5,0.6,0.7,0.8,0.9,1},
yticklabels={,0.3,0.4,0.5,0.6,0.7,0.8,0.9,1.0}
]
\addlegendimage{no markers, thick, color0}
\addlegendimage{no markers, thick, color1, dashed}

\addplot [line width=0.4800000000000001pt, color0, forget plot, mark=x]
table [row sep=\\]{%
1	0.450639 \\
};
\addplot [line width=0.4800000000000001pt, color0, forget plot]
table [row sep=\\]{%
1	0.450639 \\
1.1	0.662021 \\
1.2	0.751452 \\
1.3	0.810685 \\
1.4	0.844367 \\
1.5	0.869919 \\
1.6	0.881533 \\
1.7	0.896632 \\
1.8	0.9036 \\
1.9	0.908246 \\
2	0.914053 \\
2.1	0.915215 \\
2.2	0.919861 \\
2.3	0.927991 \\
2.4	0.930314 \\
2.5	0.933798 \\
2.6	0.937282 \\
2.7	0.938444 \\
2.8	0.938444 \\
2.9	0.940767 \\
3	0.941928 \\
3.1	0.943089 \\
3.2	0.945412 \\
3.3	0.945412 \\
3.4	0.945412 \\
3.5	0.945412 \\
3.6	0.945412 \\
3.7	0.946574 \\
3.8	0.946574 \\
3.9	0.946574 \\
4	0.946574 \\
4.1	0.946574 \\
4.2	0.948897 \\
4.3	0.948897 \\
4.4	0.950058 \\
4.5	0.953542 \\
4.6	0.953542 \\
4.7	0.953542 \\
4.8	0.953542 \\
4.9	0.953542 \\
5	0.953542 \\
5.1	0.953542 \\
5.2	0.953542 \\
5.3	0.953542 \\
5.4	0.954704 \\
5.5	0.955865 \\
5.6	0.955865 \\
5.7	0.955865 \\
5.8	0.955865 \\
5.9	0.955865 \\
6	0.955865 \\
6.1	0.955865 \\
6.19999999999999	0.955865 \\
6.29999999999999	0.955865 \\
6.39999999999999	0.955865 \\
6.49999999999999	0.955865 \\
6.59999999999999	0.957027 \\
6.69999999999999	0.958188 \\
6.79999999999999	0.95935 \\
6.89999999999999	0.95935 \\
6.99999999999999	0.95935 \\
7.09999999999999	0.95935 \\
7.19999999999999	0.95935 \\
7.29999999999999	0.960511 \\
7.39999999999999	0.960511 \\
7.49999999999999	0.960511 \\
7.59999999999999	0.960511 \\
7.69999999999999	0.960511 \\
7.79999999999999	0.960511 \\
7.89999999999999	0.960511 \\
7.99999999999999	0.960511 \\
8.09999999999999	0.960511 \\
8.19999999999999	0.960511 \\
8.29999999999999	0.960511 \\
8.39999999999999	0.960511 \\
8.49999999999999	0.960511 \\
8.59999999999999	0.960511 \\
8.69999999999999	0.960511 \\
8.79999999999999	0.960511 \\
8.89999999999999	0.960511 \\
8.99999999999999	0.960511 \\
9.09999999999999	0.960511 \\
9.19999999999999	0.960511 \\
9.29999999999998	0.960511 \\
9.39999999999998	0.961672 \\
9.49999999999998	0.961672 \\
9.59999999999998	0.961672 \\
9.69999999999998	0.961672 \\
9.79999999999998	0.961672 \\
9.89999999999998	0.961672 \\
9.99999999999998	0.961672 \\
};
\addplot [line width=0.4800000000000001pt, color1, forget plot, mark=x]
table [row sep=\\]{%
1	0.5482 \\
};
\addplot [line width=0.4800000000000001pt, color1, dashed, forget plot]
table [row sep=\\]{%
1	0.5482 \\
1.1	0.739837 \\
1.2	0.815331 \\
1.3	0.875726 \\
1.4	0.897793 \\
1.5	0.917538 \\
1.6	0.929152 \\
1.7	0.931475 \\
1.8	0.934959 \\
1.9	0.938444 \\
2	0.941928 \\
2.1	0.944251 \\
2.2	0.944251 \\
2.3	0.947735 \\
2.4	0.948897 \\
2.5	0.948897 \\
2.6	0.950058 \\
2.7	0.95122 \\
2.8	0.952381 \\
2.9	0.953542 \\
3	0.953542 \\
3.1	0.953542 \\
3.2	0.955865 \\
3.3	0.955865 \\
3.4	0.955865 \\
3.5	0.957027 \\
3.6	0.957027 \\
3.7	0.957027 \\
3.8	0.957027 \\
3.9	0.957027 \\
4	0.957027 \\
4.1	0.957027 \\
4.2	0.95935 \\
4.3	0.95935 \\
4.4	0.960511 \\
4.5	0.960511 \\
4.6	0.960511 \\
4.7	0.960511 \\
4.8	0.960511 \\
4.9	0.960511 \\
5	0.960511 \\
5.1	0.960511 \\
5.2	0.960511 \\
5.3	0.960511 \\
5.4	0.960511 \\
5.5	0.960511 \\
5.6	0.960511 \\
5.7	0.960511 \\
5.8	0.961672 \\
5.9	0.961672 \\
6	0.961672 \\
6.1	0.962834 \\
6.19999999999999	0.962834 \\
6.29999999999999	0.962834 \\
6.39999999999999	0.962834 \\
6.49999999999999	0.962834 \\
6.59999999999999	0.962834 \\
6.69999999999999	0.962834 \\
6.79999999999999	0.962834 \\
6.89999999999999	0.962834 \\
6.99999999999999	0.963995 \\
7.09999999999999	0.963995 \\
7.19999999999999	0.963995 \\
7.29999999999999	0.963995 \\
7.39999999999999	0.963995 \\
7.49999999999999	0.963995 \\
7.59999999999999	0.963995 \\
7.69999999999999	0.963995 \\
7.79999999999999	0.963995 \\
7.89999999999999	0.963995 \\
7.99999999999999	0.963995 \\
8.09999999999999	0.963995 \\
8.19999999999999	0.963995 \\
8.29999999999999	0.963995 \\
8.39999999999999	0.963995 \\
8.49999999999999	0.963995 \\
8.59999999999999	0.963995 \\
8.69999999999999	0.963995 \\
8.79999999999999	0.963995 \\
8.89999999999999	0.963995 \\
8.99999999999999	0.963995 \\
9.09999999999999	0.963995 \\
9.19999999999999	0.963995 \\
9.29999999999998	0.963995 \\
9.39999999999998	0.963995 \\
9.49999999999998	0.963995 \\
9.59999999999998	0.963995 \\
9.69999999999998	0.965157 \\
9.79999999999998	0.965157 \\
9.89999999999998	0.965157 \\
9.99999999999998	0.965157 \\
};
\path [draw=white!80.0!black, fill opacity=0] (axis cs:0,0.3)
--(axis cs:0,1);

\path [draw=white!80.0!black, fill opacity=0] (axis cs:1,0.3)
--(axis cs:1,1);

\path [draw=white!80.0!black, fill opacity=0] (axis cs:1,0)
--(axis cs:10,0);

\path [draw=white!80.0!black, fill opacity=0] (axis cs:1,1)
--(axis cs:10,1);

\end{axis}

\end{tikzpicture}
        \caption{\affected}
        }
    \end{subfigure}%
    \hfill%
    \begin{subfigure}[t][][t]{\ifreport.5\else.25\fi\textwidth}
        {
        \setlength\figureheight{.3\textheight}
        \setlength\figurewidth{\textwidth}
\begin{tikzpicture}

\definecolor{color0}{rgb}{0.298039215686275,0.447058823529412,0.690196078431373}
\definecolor{color1}{rgb}{0.866666666666667,0.517647058823529,0.32156862745098}

\begin{axis}[
axis line style={white!80.0!black},
height=\figureheight,
tick pos=both,
width=\figurewidth,
x grid style={white!80.0!black},
xlabel={Time factor},
xmajorgrids,
xmin=1, xmax=10,
y grid style={white!80.0!black},
ymajorgrids,
ymin=0.3, ymax=1,
ytick={0.2,0.3,0.4,0.5,0.6,0.7,0.8,0.9,1},
yticklabels={,0.3,0.4,0.5,0.6,0.7,0.8,0.9,1.0}
]
\addlegendimage{no markers, color0}
\addlegendimage{no markers, color1}
\addplot [line width=0.4800000000000001pt, color0, forget plot, mark=x]
table [row sep=\\]{%
1	0.447742 \\
};
\addplot [line width=0.4800000000000001pt, color0, forget plot]
table [row sep=\\]{%
1	0.447742 \\
1.1	0.645161 \\
1.2	0.736774 \\
1.3	0.796129 \\
1.4	0.832258 \\
1.5	0.858065 \\
1.6	0.870968 \\
1.7	0.885161 \\
1.8	0.892903 \\
1.9	0.898065 \\
2	0.904516 \\
2.1	0.905806 \\
2.2	0.910968 \\
2.3	0.92 \\
2.4	0.922581 \\
2.5	0.926452 \\
2.6	0.930323 \\
2.7	0.931613 \\
2.8	0.931613 \\
2.9	0.934194 \\
3	0.935484 \\
3.1	0.936774 \\
3.2	0.939355 \\
3.3	0.939355 \\
3.4	0.939355 \\
3.5	0.939355 \\
3.6	0.939355 \\
3.7	0.940645 \\
3.8	0.940645 \\
3.9	0.940645 \\
4	0.940645 \\
4.1	0.940645 \\
4.2	0.943226 \\
4.3	0.943226 \\
4.4	0.944516 \\
4.5	0.948387 \\
4.6	0.948387 \\
4.7	0.948387 \\
4.8	0.948387 \\
4.9	0.948387 \\
5	0.948387 \\
5.1	0.948387 \\
5.2	0.948387 \\
5.3	0.948387 \\
5.4	0.949677 \\
5.5	0.950968 \\
5.6	0.950968 \\
5.7	0.950968 \\
5.8	0.950968 \\
5.9	0.950968 \\
6	0.950968 \\
6.1	0.950968 \\
6.19999999999999	0.950968 \\
6.29999999999999	0.950968 \\
6.39999999999999	0.950968 \\
6.49999999999999	0.950968 \\
6.59999999999999	0.952258 \\
6.69999999999999	0.953548 \\
6.79999999999999	0.954839 \\
6.89999999999999	0.954839 \\
6.99999999999999	0.954839 \\
7.09999999999999	0.954839 \\
7.19999999999999	0.954839 \\
7.29999999999999	0.956129 \\
7.39999999999999	0.956129 \\
7.49999999999999	0.956129 \\
7.59999999999999	0.956129 \\
7.69999999999999	0.956129 \\
7.79999999999999	0.956129 \\
7.89999999999999	0.956129 \\
7.99999999999999	0.956129 \\
8.09999999999999	0.956129 \\
8.19999999999999	0.956129 \\
8.29999999999999	0.956129 \\
8.39999999999999	0.956129 \\
8.49999999999999	0.956129 \\
8.59999999999999	0.956129 \\
8.69999999999999	0.956129 \\
8.79999999999999	0.956129 \\
8.89999999999999	0.956129 \\
8.99999999999999	0.956129 \\
9.09999999999999	0.956129 \\
9.19999999999999	0.956129 \\
9.29999999999998	0.956129 \\
9.39999999999998	0.957419 \\
9.49999999999998	0.957419 \\
9.59999999999998	0.957419 \\
9.69999999999998	0.957419 \\
9.79999999999998	0.957419 \\
9.89999999999998	0.957419 \\
9.99999999999998	0.957419 \\
};
\addplot [line width=0.4800000000000001pt, color1, forget plot, mark=x]
table [row sep=\\]{%
1	0.544516 \\
};
\addplot [line width=0.4800000000000001pt, color1, dashed, forget plot]
table [row sep=\\]{%
1	0.544516 \\
1.1	0.730323 \\
1.2	0.805161 \\
1.3	0.865806 \\
1.4	0.889032 \\
1.5	0.909677 \\
1.6	0.92129 \\
1.7	0.923871 \\
1.8	0.927742 \\
1.9	0.931613 \\
2	0.935484 \\
2.1	0.938065 \\
2.2	0.938065 \\
2.3	0.941935 \\
2.4	0.943226 \\
2.5	0.943226 \\
2.6	0.944516 \\
2.7	0.945806 \\
2.8	0.947097 \\
2.9	0.948387 \\
3	0.948387 \\
3.1	0.948387 \\
3.2	0.950968 \\
3.3	0.950968 \\
3.4	0.950968 \\
3.5	0.952258 \\
3.6	0.952258 \\
3.7	0.952258 \\
3.8	0.952258 \\
3.9	0.952258 \\
4	0.952258 \\
4.1	0.952258 \\
4.2	0.954839 \\
4.3	0.954839 \\
4.4	0.956129 \\
4.5	0.956129 \\
4.6	0.956129 \\
4.7	0.956129 \\
4.8	0.956129 \\
4.9	0.956129 \\
5	0.956129 \\
5.1	0.956129 \\
5.2	0.956129 \\
5.3	0.956129 \\
5.4	0.956129 \\
5.5	0.956129 \\
5.6	0.956129 \\
5.7	0.956129 \\
5.8	0.957419 \\
5.9	0.957419 \\
6	0.957419 \\
6.1	0.95871 \\
6.19999999999999	0.95871 \\
6.29999999999999	0.95871 \\
6.39999999999999	0.95871 \\
6.49999999999999	0.95871 \\
6.59999999999999	0.95871 \\
6.69999999999999	0.95871 \\
6.79999999999999	0.95871 \\
6.89999999999999	0.95871 \\
6.99999999999999	0.96 \\
7.09999999999999	0.96 \\
7.19999999999999	0.96 \\
7.29999999999999	0.96 \\
7.39999999999999	0.96 \\
7.49999999999999	0.96 \\
7.59999999999999	0.96 \\
7.69999999999999	0.96 \\
7.79999999999999	0.96 \\
7.89999999999999	0.96 \\
7.99999999999999	0.96 \\
8.09999999999999	0.96 \\
8.19999999999999	0.96 \\
8.29999999999999	0.96 \\
8.39999999999999	0.96 \\
8.49999999999999	0.96 \\
8.59999999999999	0.96 \\
8.69999999999999	0.96 \\
8.79999999999999	0.96 \\
8.89999999999999	0.96 \\
8.99999999999999	0.96 \\
9.09999999999999	0.96 \\
9.19999999999999	0.96 \\
9.29999999999998	0.96 \\
9.39999999999998	0.96 \\
9.49999999999998	0.96 \\
9.59999999999998	0.96 \\
9.69999999999998	0.96129 \\
9.79999999999998	0.96129 \\
9.89999999999998	0.96129 \\
9.99999999999998	0.96129 \\
};
\path [draw=white!80.0!black, fill opacity=0] (axis cs:0,0.3)
--(axis cs:0,1);

\path [draw=white!80.0!black, fill opacity=0] (axis cs:1,0.3)
--(axis cs:1,1);

\path [draw=white!80.0!black, fill opacity=0] (axis cs:1,0)
--(axis cs:10,0);

\path [draw=white!80.0!black, fill opacity=0] (axis cs:1,1)
--(axis cs:10,1);

\end{axis}

\end{tikzpicture}
        \caption{\bracket{10}{tilim}}
        }
    \end{subfigure}%
    \ifreport
        \medskip
        
    \else
        \hfill%
    \fi
    \begin{subfigure}[t][][t]{\ifreport.5\else.25\fi\textwidth}
        {
        \setlength\figureheight{.3\textheight}
        \setlength\figurewidth{\textwidth}
\begin{tikzpicture}

\definecolor{color0}{rgb}{0.298039215686275,0.447058823529412,0.690196078431373}
\definecolor{color1}{rgb}{0.866666666666667,0.517647058823529,0.32156862745098}

\begin{axis}[
axis line style={white!80.0!black},
height=\figureheight,
tick pos=both,
width=\figurewidth,
x grid style={white!80.0!black},
xlabel={Time factor},
xmajorgrids,
xmin=1, xmax=10,
y grid style={white!80.0!black},
ymajorgrids,
ymin=0.3, ymax=1,
ytick={0.2,0.3,0.4,0.5,0.6,0.7,0.8,0.9,1},
yticklabels={,0.3,0.4,0.5,0.6,0.7,0.8,0.9,1.0}
]
\addlegendimage{no markers, color0}
\addlegendimage{no markers, color1}
\addplot [line width=0.4800000000000001pt, color0, forget plot, mark=x]
table [row sep=\\]{%
1	0.464945 \\
};
\addplot [line width=0.4800000000000001pt, color0, forget plot]
table [row sep=\\]{%
1	0.464945 \\
1.1	0.614391 \\
1.2	0.695572 \\
1.3	0.745387 \\
1.4	0.787823 \\
1.5	0.815498 \\
1.6	0.830258 \\
1.7	0.848708 \\
1.8	0.857934 \\
1.9	0.863469 \\
2	0.872694 \\
2.1	0.874539 \\
2.2	0.881919 \\
2.3	0.891144 \\
2.4	0.894834 \\
2.5	0.898524 \\
2.6	0.904059 \\
2.7	0.905904 \\
2.8	0.905904 \\
2.9	0.909594 \\
3	0.911439 \\
3.1	0.913284 \\
3.2	0.916974 \\
3.3	0.916974 \\
3.4	0.916974 \\
3.5	0.916974 \\
3.6	0.916974 \\
3.7	0.918819 \\
3.8	0.918819 \\
3.9	0.918819 \\
4	0.918819 \\
4.1	0.918819 \\
4.2	0.920664 \\
4.3	0.920664 \\
4.4	0.922509 \\
4.5	0.928044 \\
4.6	0.928044 \\
4.7	0.928044 \\
4.8	0.928044 \\
4.9	0.928044 \\
5	0.928044 \\
5.1	0.928044 \\
5.2	0.928044 \\
5.3	0.928044 \\
5.4	0.929889 \\
5.5	0.931734 \\
5.6	0.931734 \\
5.7	0.931734 \\
5.8	0.931734 \\
5.9	0.931734 \\
6	0.931734 \\
6.1	0.931734 \\
6.19999999999999	0.931734 \\
6.29999999999999	0.931734 \\
6.39999999999999	0.931734 \\
6.49999999999999	0.931734 \\
6.59999999999999	0.933579 \\
6.69999999999999	0.933579 \\
6.79999999999999	0.935424 \\
6.89999999999999	0.935424 \\
6.99999999999999	0.935424 \\
7.09999999999999	0.935424 \\
7.19999999999999	0.935424 \\
7.29999999999999	0.937269 \\
7.39999999999999	0.937269 \\
7.49999999999999	0.937269 \\
7.59999999999999	0.937269 \\
7.69999999999999	0.937269 \\
7.79999999999999	0.937269 \\
7.89999999999999	0.937269 \\
7.99999999999999	0.937269 \\
8.09999999999999	0.937269 \\
8.19999999999999	0.937269 \\
8.29999999999999	0.937269 \\
8.39999999999999	0.937269 \\
8.49999999999999	0.937269 \\
8.59999999999999	0.937269 \\
8.69999999999999	0.937269 \\
8.79999999999999	0.937269 \\
8.89999999999999	0.937269 \\
8.99999999999999	0.937269 \\
9.09999999999999	0.937269 \\
9.19999999999999	0.937269 \\
9.29999999999998	0.937269 \\
9.39999999999998	0.939114 \\
9.49999999999998	0.939114 \\
9.59999999999998	0.939114 \\
9.69999999999998	0.939114 \\
9.79999999999998	0.939114 \\
9.89999999999998	0.939114 \\
9.99999999999998	0.939114 \\
};
\addplot [line width=0.4800000000000001pt, color1, forget plot, mark=x]
table [row sep=\\]{%
1	0.52214 \\
};
\addplot [line width=0.4800000000000001pt, color1, dashed, forget plot]
table [row sep=\\]{%
1	0.52214 \\
1.1	0.695572 \\
1.2	0.774908 \\
1.3	0.832103 \\
1.4	0.857934 \\
1.5	0.881919 \\
1.6	0.894834 \\
1.7	0.896679 \\
1.8	0.900369 \\
1.9	0.905904 \\
2	0.909594 \\
2.1	0.913284 \\
2.2	0.913284 \\
2.3	0.918819 \\
2.4	0.920664 \\
2.5	0.920664 \\
2.6	0.920664 \\
2.7	0.922509 \\
2.8	0.924354 \\
2.9	0.926199 \\
3	0.926199 \\
3.1	0.926199 \\
3.2	0.929889 \\
3.3	0.929889 \\
3.4	0.929889 \\
3.5	0.931734 \\
3.6	0.931734 \\
3.7	0.931734 \\
3.8	0.931734 \\
3.9	0.931734 \\
4	0.931734 \\
4.1	0.931734 \\
4.2	0.935424 \\
4.3	0.935424 \\
4.4	0.937269 \\
4.5	0.937269 \\
4.6	0.937269 \\
4.7	0.937269 \\
4.8	0.937269 \\
4.9	0.937269 \\
5	0.937269 \\
5.1	0.937269 \\
5.2	0.937269 \\
5.3	0.937269 \\
5.4	0.937269 \\
5.5	0.937269 \\
5.6	0.937269 \\
5.7	0.937269 \\
5.8	0.939114 \\
5.9	0.939114 \\
6	0.939114 \\
6.1	0.940959 \\
6.19999999999999	0.940959 \\
6.29999999999999	0.940959 \\
6.39999999999999	0.940959 \\
6.49999999999999	0.940959 \\
6.59999999999999	0.940959 \\
6.69999999999999	0.940959 \\
6.79999999999999	0.940959 \\
6.89999999999999	0.940959 \\
6.99999999999999	0.942804 \\
7.09999999999999	0.942804 \\
7.19999999999999	0.942804 \\
7.29999999999999	0.942804 \\
7.39999999999999	0.942804 \\
7.49999999999999	0.942804 \\
7.59999999999999	0.942804 \\
7.69999999999999	0.942804 \\
7.79999999999999	0.942804 \\
7.89999999999999	0.942804 \\
7.99999999999999	0.942804 \\
8.09999999999999	0.942804 \\
8.19999999999999	0.942804 \\
8.29999999999999	0.942804 \\
8.39999999999999	0.942804 \\
8.49999999999999	0.942804 \\
8.59999999999999	0.942804 \\
8.69999999999999	0.942804 \\
8.79999999999999	0.942804 \\
8.89999999999999	0.942804 \\
8.99999999999999	0.942804 \\
9.09999999999999	0.942804 \\
9.19999999999999	0.942804 \\
9.29999999999998	0.942804 \\
9.39999999999998	0.942804 \\
9.49999999999998	0.942804 \\
9.59999999999998	0.942804 \\
9.69999999999998	0.944649 \\
9.79999999999998	0.944649 \\
9.89999999999998	0.944649 \\
9.99999999999998	0.944649 \\
};
\path [draw=white!80.0!black, fill opacity=0] (axis cs:0,0.3)
--(axis cs:0,1);

\path [draw=white!80.0!black, fill opacity=0] (axis cs:1,0.3)
--(axis cs:1,1);

\path [draw=white!80.0!black, fill opacity=0] (axis cs:1,0)
--(axis cs:10,0);

\path [draw=white!80.0!black, fill opacity=0] (axis cs:1,1)
--(axis cs:10,1);

\end{axis}

\end{tikzpicture}
        \caption{\bracket{100}{tilim}}
        }
    \end{subfigure}%
    \hfill%
    \begin{subfigure}[t][][t]{\ifreport.5\else.25\fi\textwidth}
        {
        \setlength\figureheight{.3\textheight}
        \setlength\figurewidth{\textwidth}
\begin{tikzpicture}

\definecolor{color0}{rgb}{0.298039215686275,0.447058823529412,0.690196078431373}
\definecolor{color1}{rgb}{0.866666666666667,0.517647058823529,0.32156862745098}

\begin{axis}[
axis line style={white!80.0!black},
height=\figureheight,
tick pos=both,
width=\figurewidth,
x grid style={white!80.0!black},
xlabel={Time factor},
xmajorgrids,
xmin=1, xmax=10,
y grid style={white!80.0!black},
ymajorgrids,
ymin=0.3, ymax=1,
ytick={0.2,0.3,0.4,0.5,0.6,0.7,0.8,0.9,1},
yticklabels={,0.3,0.4,0.5,0.6,0.7,0.8,0.9,1.0}
]
\addlegendimage{no markers, color0}
\addlegendimage{no markers, color1}
\addplot [line width=0.4800000000000001pt, color0, forget plot, mark=x]
table [row sep=\\]{%
1	0.45122 \\
};
\addplot [line width=0.4800000000000001pt, color0, forget plot]
table [row sep=\\]{%
1	0.45122 \\
1.1	0.552846 \\
1.2	0.617886 \\
1.3	0.658537 \\
1.4	0.699187 \\
1.5	0.723577 \\
1.6	0.735772 \\
1.7	0.756098 \\
1.8	0.772358 \\
1.9	0.776423 \\
2	0.776423 \\
2.1	0.776423 \\
2.2	0.784553 \\
2.3	0.792683 \\
2.4	0.800813 \\
2.5	0.804878 \\
2.6	0.808943 \\
2.7	0.813008 \\
2.8	0.813008 \\
2.9	0.821138 \\
3	0.825203 \\
3.1	0.829268 \\
3.2	0.833333 \\
3.3	0.833333 \\
3.4	0.833333 \\
3.5	0.833333 \\
3.6	0.833333 \\
3.7	0.833333 \\
3.8	0.833333 \\
3.9	0.833333 \\
4	0.833333 \\
4.1	0.833333 \\
4.2	0.837398 \\
4.3	0.837398 \\
4.4	0.837398 \\
4.5	0.845528 \\
4.6	0.845528 \\
4.7	0.845528 \\
4.8	0.845528 \\
4.9	0.845528 \\
5	0.845528 \\
5.1	0.845528 \\
5.2	0.845528 \\
5.3	0.845528 \\
5.4	0.849593 \\
5.5	0.853659 \\
5.6	0.853659 \\
5.7	0.853659 \\
5.8	0.853659 \\
5.9	0.853659 \\
6	0.853659 \\
6.1	0.853659 \\
6.19999999999999	0.853659 \\
6.29999999999999	0.853659 \\
6.39999999999999	0.853659 \\
6.49999999999999	0.853659 \\
6.59999999999999	0.857724 \\
6.69999999999999	0.857724 \\
6.79999999999999	0.857724 \\
6.89999999999999	0.857724 \\
6.99999999999999	0.857724 \\
7.09999999999999	0.857724 \\
7.19999999999999	0.857724 \\
7.29999999999999	0.861789 \\
7.39999999999999	0.861789 \\
7.49999999999999	0.861789 \\
7.59999999999999	0.861789 \\
7.69999999999999	0.861789 \\
7.79999999999999	0.861789 \\
7.89999999999999	0.861789 \\
7.99999999999999	0.861789 \\
8.09999999999999	0.861789 \\
8.19999999999999	0.861789 \\
8.29999999999999	0.861789 \\
8.39999999999999	0.861789 \\
8.49999999999999	0.861789 \\
8.59999999999999	0.861789 \\
8.69999999999999	0.861789 \\
8.79999999999999	0.861789 \\
8.89999999999999	0.861789 \\
8.99999999999999	0.861789 \\
9.09999999999999	0.861789 \\
9.19999999999999	0.861789 \\
9.29999999999998	0.861789 \\
9.39999999999998	0.865854 \\
9.49999999999998	0.865854 \\
9.59999999999998	0.865854 \\
9.69999999999998	0.865854 \\
9.79999999999998	0.865854 \\
9.89999999999998	0.865854 \\
9.99999999999998	0.865854 \\
};
\addplot [line width=0.4800000000000001pt, color1, forget plot, mark=x]
table [row sep=\\]{%
1	0.520325 \\
};
\addplot [line width=0.4800000000000001pt, color1, dashed, forget plot]
table [row sep=\\]{%
1	0.520325 \\
1.1	0.658537 \\
1.2	0.727642 \\
1.3	0.768293 \\
1.4	0.792683 \\
1.5	0.825203 \\
1.6	0.837398 \\
1.7	0.837398 \\
1.8	0.837398 \\
1.9	0.841463 \\
2	0.849593 \\
2.1	0.849593 \\
2.2	0.849593 \\
2.3	0.853659 \\
2.4	0.857724 \\
2.5	0.857724 \\
2.6	0.857724 \\
2.7	0.861789 \\
2.8	0.865854 \\
2.9	0.869919 \\
3	0.869919 \\
3.1	0.869919 \\
3.2	0.873984 \\
3.3	0.873984 \\
3.4	0.873984 \\
3.5	0.873984 \\
3.6	0.873984 \\
3.7	0.873984 \\
3.8	0.873984 \\
3.9	0.873984 \\
4	0.873984 \\
4.1	0.873984 \\
4.2	0.878049 \\
4.3	0.878049 \\
4.4	0.882114 \\
4.5	0.882114 \\
4.6	0.882114 \\
4.7	0.882114 \\
4.8	0.882114 \\
4.9	0.882114 \\
5	0.882114 \\
5.1	0.882114 \\
5.2	0.882114 \\
5.3	0.882114 \\
5.4	0.882114 \\
5.5	0.882114 \\
5.6	0.882114 \\
5.7	0.882114 \\
5.8	0.886179 \\
5.9	0.886179 \\
6	0.886179 \\
6.1	0.890244 \\
6.19999999999999	0.890244 \\
6.29999999999999	0.890244 \\
6.39999999999999	0.890244 \\
6.49999999999999	0.890244 \\
6.59999999999999	0.890244 \\
6.69999999999999	0.890244 \\
6.79999999999999	0.890244 \\
6.89999999999999	0.890244 \\
6.99999999999999	0.890244 \\
7.09999999999999	0.890244 \\
7.19999999999999	0.890244 \\
7.29999999999999	0.890244 \\
7.39999999999999	0.890244 \\
7.49999999999999	0.890244 \\
7.59999999999999	0.890244 \\
7.69999999999999	0.890244 \\
7.79999999999999	0.890244 \\
7.89999999999999	0.890244 \\
7.99999999999999	0.890244 \\
8.09999999999999	0.890244 \\
8.19999999999999	0.890244 \\
8.29999999999999	0.890244 \\
8.39999999999999	0.890244 \\
8.49999999999999	0.890244 \\
8.59999999999999	0.890244 \\
8.69999999999999	0.890244 \\
8.79999999999999	0.890244 \\
8.89999999999999	0.890244 \\
8.99999999999999	0.890244 \\
9.09999999999999	0.890244 \\
9.19999999999999	0.890244 \\
9.29999999999998	0.890244 \\
9.39999999999998	0.890244 \\
9.49999999999998	0.890244 \\
9.59999999999998	0.890244 \\
9.69999999999998	0.890244 \\
9.79999999999998	0.890244 \\
9.89999999999998	0.890244 \\
9.99999999999998	0.890244 \\
};
\path [draw=white!80.0!black, fill opacity=0] (axis cs:0,0.3)
--(axis cs:0,1);

\path [draw=white!80.0!black, fill opacity=0] (axis cs:1,0.3)
--(axis cs:1,1);

\path [draw=white!80.0!black, fill opacity=0] (axis cs:1,0)
--(axis cs:10,0);

\path [draw=white!80.0!black, fill opacity=0] (axis cs:1,1)
--(axis cs:10,1);

\end{axis}

\end{tikzpicture}
        \caption{\bracket{1000}{tilim}}
        }
    \end{subfigure}%
  \end{center}
  \smallskip
  
  \caption{Performance profiles of \coefdiving and \confdiving for four hierarchical groups of increasingly hard affected instances.}
  \label{fig:coef_vs_conf_performance}
\end{figure}
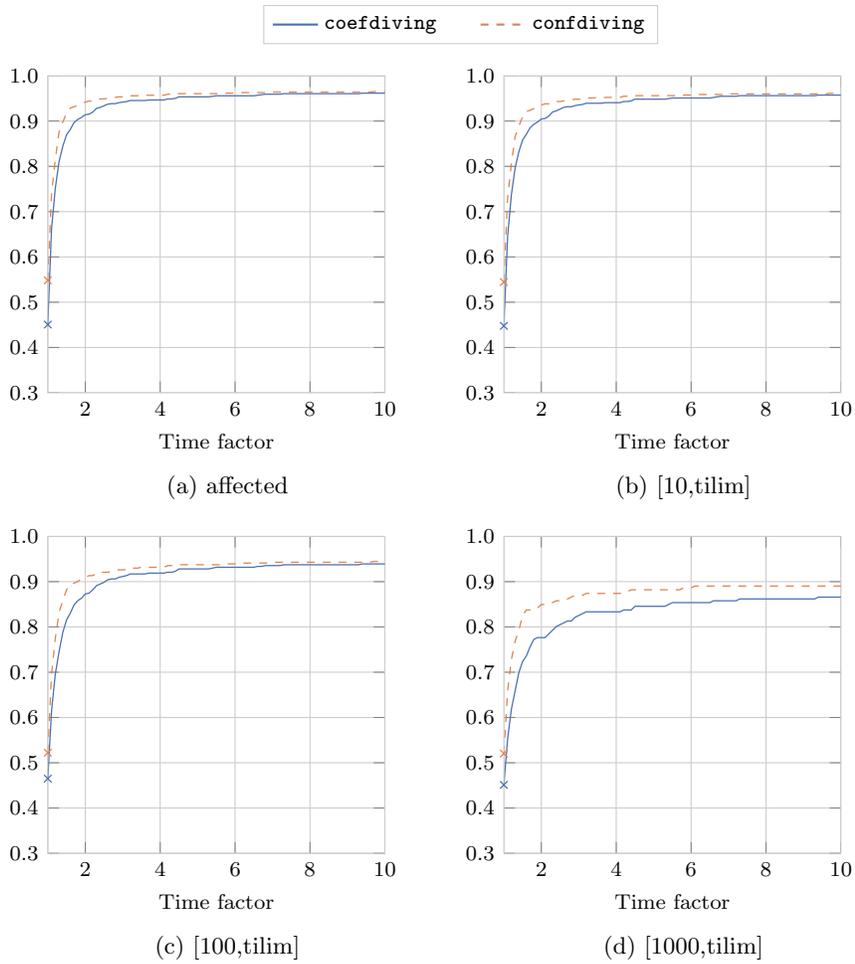

Finally, the performance profiles in Figure~\ref{fig:coef_vs_conf_performance} show that the overall performance improvement was not caused by few instances
but was spread over the complete group of affected instances:
\confdiving dominates \coefdiving on all four groups of increasingly hard instances.
For a time factor of $1.0$ the respective percentage of instances that could be solved fastest with the respective setting is marked with a colored cross.
On all four groups of instances, \confdiving is clearly superior to \coefdiving in the sense that the profiles do not cross.
Between $\percent{52.0}$ and $\percent{54.8}$ of the instances in the respective group were solved fastest by \confdiving,
whereas only between $\percent{44.8}$ and $\percent{46.5}$ of the instances were solved fastest by \default.
 
\paragraph{Length of Diving Paths.}

Both the rounding and scoring functions used by \conflictdiving aim for a fail fast strategy~\citep{haralick1980increasing,Berthold2014}.
In order to analyze whether \conflictdiving indeed achieves this design goal, we additionally measured the average length of diving paths for \coefdiving and \confdiving.
On average, \confdiving exhibits $\percent{30.2}$ shorter diving paths than \coefdiving on affected instances.

The distribution of the length of diving paths over the subsets of affected instances can be seen in Figure~\ref{fig:confdiving_avgdepth}.
\begin{figure}%
  \begin{center}
    \setlength\figureheight{.15\textheight}
    \setlength\figurewidth{\ifreport.75\else.5\fi\textwidth}

  \end{center}
  \caption{Box plots showing the average depth of diving paths generated \coefficientdiving with \coefdiving and \conflictdiving
           with \confdiving on the set of affected instances.}
  \label{fig:confdiving_avgdepth}
\end{figure}
The box plots show for every setting and instance group the 1st and 3rd quartile of all diving path lengths (shaded box)
as well the median (dashed line).
All observations below the 1st or above the 3rd quartile are marked with shaded diamonds.
For all four instance groups, the 1st and 3rd quartile corresponding to \conflictdiving are smaller than for \coefdiving.
The same observation holds for the median in all cases.
For example, on \bracket{100}{tilim} the 3rd quartile ($149$) of \coefficientdiving is $\percent{69.3}$ larger than the 3rd quartile ($87$) of \conflictdiving.
On the same group of instances, the observed median path length of \conflictdiving is $\percent{23.6}$ shorter than the one of \coefficientdiving.

As we already discussed, both \coefficientdiving and \conflictdiving showed only a small number of found solutions per diving path,
which is a good proxy for the number of successful paths, \ie paths without backtracking due to infeasibility (cf.~Line~\ref{line:conflict1} and~\ref{line:conflict2}).
Thus, these statistics confirm that \conflictdiving succeeds better in implementing a fail fast strategy and aborting the many ``unsuccessful'' dives not leading to a feasible solution early.

Moreover, this does not come at the expense of learning less conflict constraints.
As can be seen in column \confs of Table~\ref{tab:conflictdiving_aggregated}, \confdiving created $\percent{30.0}$ more conflict constraints.

\medskip

\paragraph{Impact of Conflict Locks.}

In order to investigate whether \conflictdiving outperforms \coefficientdiving because of the inclusion of conflict locks
or merely because of the difference in the rounding and scoring function, we conducted two further control experiments with
\begin{itemize}
 \item a modified rounding function within \conflictdiving and
 \item different weights $\confweight$.
\end{itemize}
\smallskip

For the first control experiment we modified \conflictdiving to use the same rounding function as \coefdiving and $\confweight = 0.5$ for the scoring function.
The resulting diving heuristic \confdivinglikecoef behaves the same as \coefficientdiving
except for using a convex combination of variable and conflict locks in the variable selection score.

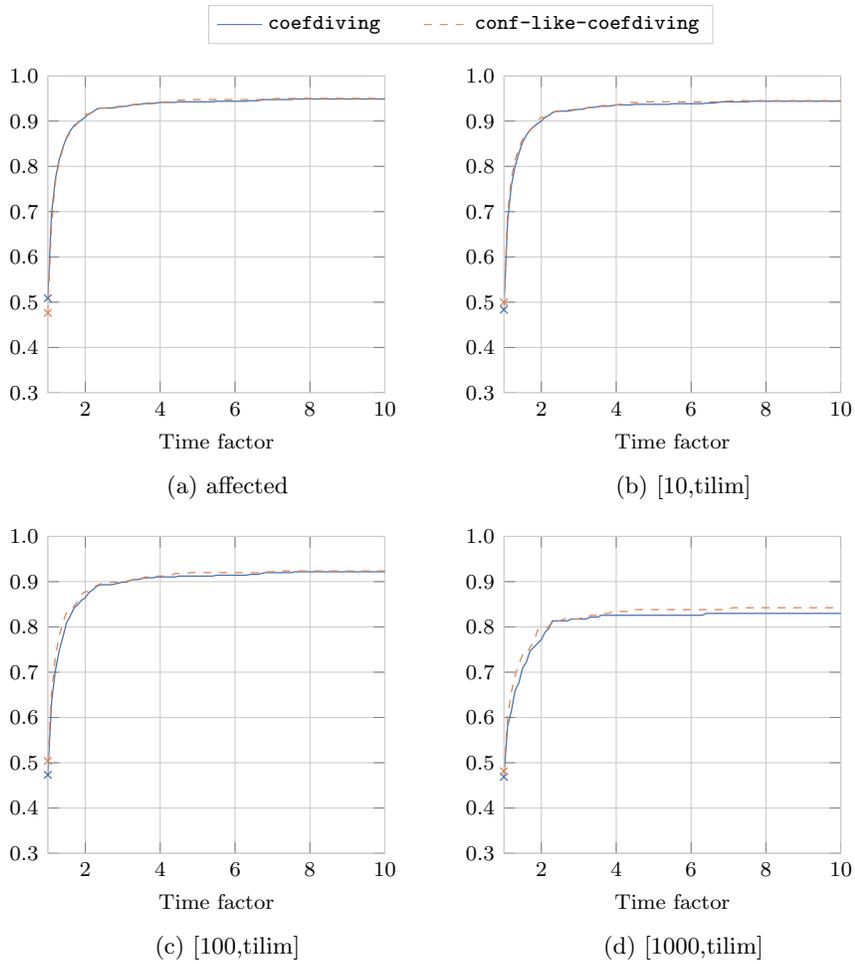
\begin{figure}
  \begin{center}%
    \begin{subfigure}[t][][t]{\ifreport.5\else.25\fi\textwidth}
        {
        \setlength\figureheight{.3\textheight}
        \setlength\figurewidth{\textwidth}
\begin{tikzpicture}

\definecolor{color0}{rgb}{0.298039215686275,0.447058823529412,0.690196078431373}
\definecolor{color1}{rgb}{0.866666666666667,0.517647058823529,0.32156862745098}

\begin{axis}[
axis line style={white!80.0!black},
height=\figureheight,
legend cell align={left},
legend entries={{\coefdiving},{\confdivinglikecoef}},
legend style={at={\ifreport(1.225,1.1)\else(3,1.10)\fi}, anchor=south, draw=white!80.0!black, font=\footnotesize, /tikz/every even column/.append style={column sep=0.5cm}},
legend columns=2,
tick pos=both,
width=\figurewidth,
x grid style={white!80.0!black},
xlabel={Time factor},
xmajorgrids,
xmin=1, xmax=10,
y grid style={white!80.0!black},
ymajorgrids,
ymin=0.3, ymax=1,
ytick={0.2,0.3,0.4,0.5,0.6,0.7,0.8,0.9,1},
yticklabels={,0.3,0.4,0.5,0.6,0.7,0.8,0.9,1.0}
]
\addlegendimage{no markers, color0}
\addlegendimage{no markers, color1, dashed}
\addplot [line width=0.4800000000000001pt, color0, forget plot, mark=x]
table [row sep=\\]{%
1	0.50875 \\
};
\addplot [line width=0.4800000000000001pt, color0, forget plot]
table [row sep=\\]{%
1	0.50875 \\
1.1	0.69625 \\
1.2	0.77125 \\
1.3	0.81375 \\
1.4	0.84125 \\
1.5	0.86375 \\
1.6	0.87875 \\
1.7	0.89 \\
1.8	0.89625 \\
1.9	0.9025 \\
2	0.90875 \\
2.1	0.91625 \\
2.2	0.92 \\
2.3	0.92625 \\
2.4	0.92875 \\
2.5	0.92875 \\
2.6	0.92875 \\
2.7	0.92875 \\
2.8	0.93 \\
2.9	0.93125 \\
3	0.9325 \\
3.1	0.9325 \\
3.2	0.935 \\
3.3	0.93625 \\
3.4	0.93625 \\
3.5	0.9375 \\
3.6	0.93875 \\
3.7	0.93875 \\
3.8	0.93875 \\
3.9	0.94 \\
4	0.94125 \\
4.1	0.94125 \\
4.2	0.94125 \\
4.3	0.94125 \\
4.4	0.94125 \\
4.5	0.9425 \\
4.6	0.9425 \\
4.7	0.9425 \\
4.8	0.9425 \\
4.9	0.9425 \\
5	0.9425 \\
5.1	0.9425 \\
5.2	0.9425 \\
5.3	0.9425 \\
5.4	0.9425 \\
5.5	0.94375 \\
5.6	0.94375 \\
5.7	0.94375 \\
5.8	0.94375 \\
5.9	0.94375 \\
6	0.94375 \\
6.1	0.94375 \\
6.19999999999999	0.94375 \\
6.29999999999999	0.94375 \\
6.39999999999999	0.945 \\
6.49999999999999	0.945 \\
6.59999999999999	0.945 \\
6.69999999999999	0.94625 \\
6.79999999999999	0.9475 \\
6.89999999999999	0.9475 \\
6.99999999999999	0.9475 \\
7.09999999999999	0.9475 \\
7.19999999999999	0.9475 \\
7.29999999999999	0.9475 \\
7.39999999999999	0.9475 \\
7.49999999999999	0.9475 \\
7.59999999999999	0.94875 \\
7.69999999999999	0.94875 \\
7.79999999999999	0.94875 \\
7.89999999999999	0.94875 \\
7.99999999999999	0.94875 \\
8.09999999999999	0.94875 \\
8.19999999999999	0.94875 \\
8.29999999999999	0.94875 \\
8.39999999999999	0.94875 \\
8.49999999999999	0.94875 \\
8.59999999999999	0.94875 \\
8.69999999999999	0.94875 \\
8.79999999999999	0.94875 \\
8.89999999999999	0.94875 \\
8.99999999999999	0.94875 \\
9.09999999999999	0.94875 \\
9.19999999999999	0.94875 \\
9.29999999999998	0.94875 \\
9.39999999999998	0.94875 \\
9.49999999999998	0.94875 \\
9.59999999999998	0.94875 \\
9.69999999999998	0.94875 \\
9.79999999999998	0.94875 \\
9.89999999999998	0.94875 \\
9.99999999999998	0.94875 \\
};
\addplot [line width=0.4800000000000001pt, color1, forget plot, mark=x]
table [row sep=\\]{%
1	0.47625 \\
};
\addplot [line width=0.4800000000000001pt, color1, dashed, forget plot]
table [row sep=\\]{%
1	0.47625 \\
1.1	0.67 \\
1.2	0.76375 \\
1.3	0.80875 \\
1.4	0.8375 \\
1.5	0.86125 \\
1.6	0.875 \\
1.7	0.88625 \\
1.8	0.89625 \\
1.9	0.90875 \\
2	0.915 \\
2.1	0.915 \\
2.2	0.92 \\
2.3	0.925 \\
2.4	0.92875 \\
2.5	0.92875 \\
2.6	0.93 \\
2.7	0.93125 \\
2.8	0.93125 \\
2.9	0.9325 \\
3	0.9325 \\
3.1	0.93375 \\
3.2	0.935 \\
3.3	0.93625 \\
3.4	0.9375 \\
3.5	0.9375 \\
3.6	0.94 \\
3.7	0.94 \\
3.8	0.94 \\
3.9	0.94125 \\
4	0.94125 \\
4.1	0.9425 \\
4.2	0.9425 \\
4.3	0.94375 \\
4.4	0.94625 \\
4.5	0.94625 \\
4.6	0.94625 \\
4.7	0.94625 \\
4.8	0.9475 \\
4.9	0.9475 \\
5	0.9475 \\
5.1	0.9475 \\
5.2	0.9475 \\
5.3	0.9475 \\
5.4	0.9475 \\
5.5	0.9475 \\
5.6	0.9475 \\
5.7	0.9475 \\
5.8	0.9475 \\
5.9	0.9475 \\
6	0.9475 \\
6.1	0.9475 \\
6.19999999999999	0.9475 \\
6.29999999999999	0.9475 \\
6.39999999999999	0.9475 \\
6.49999999999999	0.9475 \\
6.59999999999999	0.9475 \\
6.69999999999999	0.9475 \\
6.79999999999999	0.94875 \\
6.89999999999999	0.94875 \\
6.99999999999999	0.94875 \\
7.09999999999999	0.95 \\
7.19999999999999	0.95 \\
7.29999999999999	0.95 \\
7.39999999999999	0.95 \\
7.49999999999999	0.95 \\
7.59999999999999	0.95 \\
7.69999999999999	0.95 \\
7.79999999999999	0.95 \\
7.89999999999999	0.95 \\
7.99999999999999	0.95 \\
8.09999999999999	0.95 \\
8.19999999999999	0.95 \\
8.29999999999999	0.95 \\
8.39999999999999	0.95 \\
8.49999999999999	0.95 \\
8.59999999999999	0.95 \\
8.69999999999999	0.95 \\
8.79999999999999	0.95 \\
8.89999999999999	0.95 \\
8.99999999999999	0.95 \\
9.09999999999999	0.95 \\
9.19999999999999	0.95 \\
9.29999999999998	0.95 \\
9.39999999999998	0.95 \\
9.49999999999998	0.95 \\
9.59999999999998	0.95 \\
9.69999999999998	0.95 \\
9.79999999999998	0.95 \\
9.89999999999998	0.95 \\
9.99999999999998	0.95 \\
};
\path [draw=white!80.0!black, fill opacity=0] (axis cs:0,0.3)
--(axis cs:0,1);

\path [draw=white!80.0!black, fill opacity=0] (axis cs:1,0.3)
--(axis cs:1,1);

\path [draw=white!80.0!black, fill opacity=0] (axis cs:1,0)
--(axis cs:10,0);

\path [draw=white!80.0!black, fill opacity=0] (axis cs:1,1)
--(axis cs:10,1);

\end{axis}

\end{tikzpicture}
        \caption{\affected}
        }
    \end{subfigure}%
    \hfill%
    \begin{subfigure}[t][][t]{\ifreport.5\else.25\fi\textwidth}
        {
        \setlength\figureheight{.3\textheight}
        \setlength\figurewidth{\textwidth}
\begin{tikzpicture}

\definecolor{color0}{rgb}{0.298039215686275,0.447058823529412,0.690196078431373}
\definecolor{color1}{rgb}{0.866666666666667,0.517647058823529,0.32156862745098}

\begin{axis}[
axis line style={white!80.0!black},
height=\figureheight,
tick pos=both,
width=\figurewidth,
x grid style={white!80.0!black},
xlabel={Time factor},
xmajorgrids,
xmin=1, xmax=10,
y grid style={white!80.0!black},
ymajorgrids,
ymin=0.3, ymax=1,
ytick={0.2,0.3,0.4,0.5,0.6,0.7,0.8,0.9,1},
yticklabels={,0.3,0.4,0.5,0.6,0.7,0.8,0.9,1.0}
]
\addlegendimage{no markers, color0}
\addlegendimage{no markers, color1}
\addplot [line width=0.4800000000000001pt, color0, forget plot, mark=x]
table [row sep=\\]{%
1	0.483562 \\
};
\addplot [line width=0.4800000000000001pt, color0, forget plot]
table [row sep=\\]{%
1	0.483562 \\
1.1	0.678082 \\
1.2	0.757534 \\
1.3	0.8 \\
1.4	0.828767 \\
1.5	0.853425 \\
1.6	0.868493 \\
1.7	0.880822 \\
1.8	0.887671 \\
1.9	0.894521 \\
2	0.9 \\
2.1	0.908219 \\
2.2	0.912329 \\
2.3	0.919178 \\
2.4	0.921918 \\
2.5	0.921918 \\
2.6	0.921918 \\
2.7	0.921918 \\
2.8	0.923288 \\
2.9	0.924658 \\
3	0.926027 \\
3.1	0.926027 \\
3.2	0.928767 \\
3.3	0.930137 \\
3.4	0.930137 \\
3.5	0.931507 \\
3.6	0.932877 \\
3.7	0.932877 \\
3.8	0.932877 \\
3.9	0.934247 \\
4	0.935616 \\
4.1	0.935616 \\
4.2	0.935616 \\
4.3	0.935616 \\
4.4	0.935616 \\
4.5	0.936986 \\
4.6	0.936986 \\
4.7	0.936986 \\
4.8	0.936986 \\
4.9	0.936986 \\
5	0.936986 \\
5.1	0.936986 \\
5.2	0.936986 \\
5.3	0.936986 \\
5.4	0.936986 \\
5.5	0.938356 \\
5.6	0.938356 \\
5.7	0.938356 \\
5.8	0.938356 \\
5.9	0.938356 \\
6	0.938356 \\
6.1	0.938356 \\
6.19999999999999	0.938356 \\
6.29999999999999	0.938356 \\
6.39999999999999	0.939726 \\
6.49999999999999	0.939726 \\
6.59999999999999	0.939726 \\
6.69999999999999	0.941096 \\
6.79999999999999	0.942466 \\
6.89999999999999	0.942466 \\
6.99999999999999	0.942466 \\
7.09999999999999	0.942466 \\
7.19999999999999	0.942466 \\
7.29999999999999	0.942466 \\
7.39999999999999	0.942466 \\
7.49999999999999	0.942466 \\
7.59999999999999	0.943836 \\
7.69999999999999	0.943836 \\
7.79999999999999	0.943836 \\
7.89999999999999	0.943836 \\
7.99999999999999	0.943836 \\
8.09999999999999	0.943836 \\
8.19999999999999	0.943836 \\
8.29999999999999	0.943836 \\
8.39999999999999	0.943836 \\
8.49999999999999	0.943836 \\
8.59999999999999	0.943836 \\
8.69999999999999	0.943836 \\
8.79999999999999	0.943836 \\
8.89999999999999	0.943836 \\
8.99999999999999	0.943836 \\
9.09999999999999	0.943836 \\
9.19999999999999	0.943836 \\
9.29999999999998	0.943836 \\
9.39999999999998	0.943836 \\
9.49999999999998	0.943836 \\
9.59999999999998	0.943836 \\
9.69999999999998	0.943836 \\
9.79999999999998	0.943836 \\
9.89999999999998	0.943836 \\
9.99999999999998	0.943836 \\
};
\addplot [line width=0.4800000000000001pt, color1, forget plot, mark=x]
table [row sep=\\]{%
1	0.5 \\
};
\addplot [line width=0.4800000000000001pt, color1, dashed, forget plot]
table [row sep=\\]{%
1	0.5 \\
1.1	0.69589 \\
1.2	0.773973 \\
1.3	0.815068 \\
1.4	0.839726 \\
1.5	0.860274 \\
1.6	0.869863 \\
1.7	0.879452 \\
1.8	0.889041 \\
1.9	0.90137 \\
2	0.908219 \\
2.1	0.908219 \\
2.2	0.913699 \\
2.3	0.917808 \\
2.4	0.921918 \\
2.5	0.921918 \\
2.6	0.923288 \\
2.7	0.924658 \\
2.8	0.924658 \\
2.9	0.926027 \\
3	0.926027 \\
3.1	0.927397 \\
3.2	0.928767 \\
3.3	0.930137 \\
3.4	0.931507 \\
3.5	0.931507 \\
3.6	0.934247 \\
3.7	0.934247 \\
3.8	0.934247 \\
3.9	0.935616 \\
4	0.935616 \\
4.1	0.936986 \\
4.2	0.936986 \\
4.3	0.938356 \\
4.4	0.941096 \\
4.5	0.941096 \\
4.6	0.941096 \\
4.7	0.941096 \\
4.8	0.942466 \\
4.9	0.942466 \\
5	0.942466 \\
5.1	0.942466 \\
5.2	0.942466 \\
5.3	0.942466 \\
5.4	0.942466 \\
5.5	0.942466 \\
5.6	0.942466 \\
5.7	0.942466 \\
5.8	0.942466 \\
5.9	0.942466 \\
6	0.942466 \\
6.1	0.942466 \\
6.19999999999999	0.942466 \\
6.29999999999999	0.942466 \\
6.39999999999999	0.942466 \\
6.49999999999999	0.942466 \\
6.59999999999999	0.942466 \\
6.69999999999999	0.942466 \\
6.79999999999999	0.943836 \\
6.89999999999999	0.943836 \\
6.99999999999999	0.943836 \\
7.09999999999999	0.945205 \\
7.19999999999999	0.945205 \\
7.29999999999999	0.945205 \\
7.39999999999999	0.945205 \\
7.49999999999999	0.945205 \\
7.59999999999999	0.945205 \\
7.69999999999999	0.945205 \\
7.79999999999999	0.945205 \\
7.89999999999999	0.945205 \\
7.99999999999999	0.945205 \\
8.09999999999999	0.945205 \\
8.19999999999999	0.945205 \\
8.29999999999999	0.945205 \\
8.39999999999999	0.945205 \\
8.49999999999999	0.945205 \\
8.59999999999999	0.945205 \\
8.69999999999999	0.945205 \\
8.79999999999999	0.945205 \\
8.89999999999999	0.945205 \\
8.99999999999999	0.945205 \\
9.09999999999999	0.945205 \\
9.19999999999999	0.945205 \\
9.29999999999998	0.945205 \\
9.39999999999998	0.945205 \\
9.49999999999998	0.945205 \\
9.59999999999998	0.945205 \\
9.69999999999998	0.945205 \\
9.79999999999998	0.945205 \\
9.89999999999998	0.945205 \\
9.99999999999998	0.945205 \\
};
\path [draw=white!80.0!black, fill opacity=0] (axis cs:0,0.3)
--(axis cs:0,1);

\path [draw=white!80.0!black, fill opacity=0] (axis cs:1,0.3)
--(axis cs:1,1);

\path [draw=white!80.0!black, fill opacity=0] (axis cs:1,0)
--(axis cs:10,0);

\path [draw=white!80.0!black, fill opacity=0] (axis cs:1,1)
--(axis cs:10,1);

\end{axis}

\end{tikzpicture}
        \caption{\bracket{10}{tilim}}
        }
    \end{subfigure}%
    \ifreport
        \medskip
        
    \else
        \hfill%
    \fi
    \begin{subfigure}[t][][t]{\ifreport.5\else.25\fi\textwidth}
        {
        \setlength\figureheight{.3\textheight}
        \setlength\figurewidth{\textwidth}
\begin{tikzpicture}

\definecolor{color0}{rgb}{0.298039215686275,0.447058823529412,0.690196078431373}
\definecolor{color1}{rgb}{0.866666666666667,0.517647058823529,0.32156862745098}

\begin{axis}[
axis line style={white!80.0!black},
height=\figureheight,
tick pos=both,
width=\figurewidth,
x grid style={white!80.0!black},
xlabel={Time factor},
xmajorgrids,
xmin=1, xmax=10,
y grid style={white!80.0!black},
ymajorgrids,
ymin=0.3, ymax=1,
ytick={0.2,0.3,0.4,0.5,0.6,0.7,0.8,0.9,1},
yticklabels={,0.3,0.4,0.5,0.6,0.7,0.8,0.9,1.0}
]
\addlegendimage{no markers, color0}
\addlegendimage{no markers, color1}
\addplot [line width=0.4800000000000001pt, color0, forget plot, mark=x]
table [row sep=\\]{%
1	0.473282 \\
};
\addplot [line width=0.4800000000000001pt, color0, forget plot]
table [row sep=\\]{%
1	0.473282 \\
1.1	0.633588 \\
1.2	0.70229 \\
1.3	0.746183 \\
1.4	0.776718 \\
1.5	0.80916 \\
1.6	0.824427 \\
1.7	0.841603 \\
1.8	0.849237 \\
1.9	0.858779 \\
2	0.864504 \\
2.1	0.875954 \\
2.2	0.881679 \\
2.3	0.889313 \\
2.4	0.89313 \\
2.5	0.89313 \\
2.6	0.89313 \\
2.7	0.89313 \\
2.8	0.895038 \\
2.9	0.896947 \\
3	0.898855 \\
3.1	0.898855 \\
3.2	0.902672 \\
3.3	0.90458 \\
3.4	0.90458 \\
3.5	0.906489 \\
3.6	0.908397 \\
3.7	0.908397 \\
3.8	0.908397 \\
3.9	0.910305 \\
4	0.910305 \\
4.1	0.910305 \\
4.2	0.910305 \\
4.3	0.910305 \\
4.4	0.910305 \\
4.5	0.912214 \\
4.6	0.912214 \\
4.7	0.912214 \\
4.8	0.912214 \\
4.9	0.912214 \\
5	0.912214 \\
5.1	0.912214 \\
5.2	0.912214 \\
5.3	0.912214 \\
5.4	0.912214 \\
5.5	0.914122 \\
5.6	0.914122 \\
5.7	0.914122 \\
5.8	0.914122 \\
5.9	0.914122 \\
6	0.914122 \\
6.1	0.914122 \\
6.19999999999999	0.914122 \\
6.29999999999999	0.914122 \\
6.39999999999999	0.916031 \\
6.49999999999999	0.916031 \\
6.59999999999999	0.916031 \\
6.69999999999999	0.917939 \\
6.79999999999999	0.919847 \\
6.89999999999999	0.919847 \\
6.99999999999999	0.919847 \\
7.09999999999999	0.919847 \\
7.19999999999999	0.919847 \\
7.29999999999999	0.919847 \\
7.39999999999999	0.919847 \\
7.49999999999999	0.919847 \\
7.59999999999999	0.921756 \\
7.69999999999999	0.921756 \\
7.79999999999999	0.921756 \\
7.89999999999999	0.921756 \\
7.99999999999999	0.921756 \\
8.09999999999999	0.921756 \\
8.19999999999999	0.921756 \\
8.29999999999999	0.921756 \\
8.39999999999999	0.921756 \\
8.49999999999999	0.921756 \\
8.59999999999999	0.921756 \\
8.69999999999999	0.921756 \\
8.79999999999999	0.921756 \\
8.89999999999999	0.921756 \\
8.99999999999999	0.921756 \\
9.09999999999999	0.921756 \\
9.19999999999999	0.921756 \\
9.29999999999998	0.921756 \\
9.39999999999998	0.921756 \\
9.49999999999998	0.921756 \\
9.59999999999998	0.921756 \\
9.69999999999998	0.921756 \\
9.79999999999998	0.921756 \\
9.89999999999998	0.921756 \\
9.99999999999998	0.921756 \\
};
\addplot [line width=0.4800000000000001pt, color1, forget plot, mark=x]
table [row sep=\\]{%
1	0.503817 \\
};
\addplot [line width=0.4800000000000001pt, color1, dashed, forget plot]
table [row sep=\\]{%
1	0.503817 \\
1.1	0.65458 \\
1.2	0.732824 \\
1.3	0.778626 \\
1.4	0.80916 \\
1.5	0.832061 \\
1.6	0.835878 \\
1.7	0.84542 \\
1.8	0.85687 \\
1.9	0.872137 \\
2	0.877863 \\
2.1	0.877863 \\
2.2	0.885496 \\
2.3	0.891221 \\
2.4	0.895038 \\
2.5	0.895038 \\
2.6	0.896947 \\
2.7	0.898855 \\
2.8	0.898855 \\
2.9	0.898855 \\
3	0.898855 \\
3.1	0.900763 \\
3.2	0.902672 \\
3.3	0.90458 \\
3.4	0.906489 \\
3.5	0.906489 \\
3.6	0.910305 \\
3.7	0.910305 \\
3.8	0.910305 \\
3.9	0.912214 \\
4	0.912214 \\
4.1	0.912214 \\
4.2	0.912214 \\
4.3	0.914122 \\
4.4	0.917939 \\
4.5	0.917939 \\
4.6	0.917939 \\
4.7	0.917939 \\
4.8	0.919847 \\
4.9	0.919847 \\
5	0.919847 \\
5.1	0.919847 \\
5.2	0.919847 \\
5.3	0.919847 \\
5.4	0.919847 \\
5.5	0.919847 \\
5.6	0.919847 \\
5.7	0.919847 \\
5.8	0.919847 \\
5.9	0.919847 \\
6	0.919847 \\
6.1	0.919847 \\
6.19999999999999	0.919847 \\
6.29999999999999	0.919847 \\
6.39999999999999	0.919847 \\
6.49999999999999	0.919847 \\
6.59999999999999	0.919847 \\
6.69999999999999	0.919847 \\
6.79999999999999	0.921756 \\
6.89999999999999	0.921756 \\
6.99999999999999	0.921756 \\
7.09999999999999	0.923664 \\
7.19999999999999	0.923664 \\
7.29999999999999	0.923664 \\
7.39999999999999	0.923664 \\
7.49999999999999	0.923664 \\
7.59999999999999	0.923664 \\
7.69999999999999	0.923664 \\
7.79999999999999	0.923664 \\
7.89999999999999	0.923664 \\
7.99999999999999	0.923664 \\
8.09999999999999	0.923664 \\
8.19999999999999	0.923664 \\
8.29999999999999	0.923664 \\
8.39999999999999	0.923664 \\
8.49999999999999	0.923664 \\
8.59999999999999	0.923664 \\
8.69999999999999	0.923664 \\
8.79999999999999	0.923664 \\
8.89999999999999	0.923664 \\
8.99999999999999	0.923664 \\
9.09999999999999	0.923664 \\
9.19999999999999	0.923664 \\
9.29999999999998	0.923664 \\
9.39999999999998	0.923664 \\
9.49999999999998	0.923664 \\
9.59999999999998	0.923664 \\
9.69999999999998	0.923664 \\
9.79999999999998	0.923664 \\
9.89999999999998	0.923664 \\
9.99999999999998	0.923664 \\
};
\path [draw=white!80.0!black, fill opacity=0] (axis cs:0,0.3)
--(axis cs:0,1);

\path [draw=white!80.0!black, fill opacity=0] (axis cs:1,0.3)
--(axis cs:1,1);

\path [draw=white!80.0!black, fill opacity=0] (axis cs:1,0)
--(axis cs:10,0);

\path [draw=white!80.0!black, fill opacity=0] (axis cs:1,1)
--(axis cs:10,1);

\end{axis}

\end{tikzpicture}
        \caption{\bracket{100}{tilim}}
        }
    \end{subfigure}%
    \hfill%
    \begin{subfigure}[t][][t]{\ifreport.5\else.25\fi\textwidth}
        {
        \setlength\figureheight{.3\textheight}
        \setlength\figurewidth{\textwidth}
\begin{tikzpicture}

\definecolor{color0}{rgb}{0.298039215686275,0.447058823529412,0.690196078431373}
\definecolor{color1}{rgb}{0.866666666666667,0.517647058823529,0.32156862745098}

\begin{axis}[
axis line style={white!80.0!black},
height=\figureheight,
tick pos=both,
width=\figurewidth,
x grid style={white!80.0!black},
xlabel={Time factor},
xmajorgrids,
xmin=1, xmax=10,
y grid style={white!80.0!black},
ymajorgrids,
ymin=0.3, ymax=1,
ytick={0.2,0.3,0.4,0.5,0.6,0.7,0.8,0.9,1},
yticklabels={,0.3,0.4,0.5,0.6,0.7,0.8,0.9,1.0}
]
\addlegendimage{no markers, color0}
\addlegendimage{no markers, color1}
\addplot [line width=0.4800000000000001pt, color0, forget plot, mark=x]
table [row sep=\\]{%
1	0.46888 \\
};
\addplot [line width=0.4800000000000001pt, color0, forget plot]
table [row sep=\\]{%
1	0.46888 \\
1.1	0.580913 \\
1.2	0.614108 \\
1.3	0.659751 \\
1.4	0.676349 \\
1.5	0.709544 \\
1.6	0.721992 \\
1.7	0.746888 \\
1.8	0.755187 \\
1.9	0.763485 \\
2	0.771784 \\
2.1	0.788382 \\
2.2	0.79668 \\
2.3	0.813278 \\
2.4	0.813278 \\
2.5	0.813278 \\
2.6	0.813278 \\
2.7	0.813278 \\
2.8	0.817427 \\
2.9	0.817427 \\
3	0.817427 \\
3.1	0.817427 \\
3.2	0.817427 \\
3.3	0.821577 \\
3.4	0.821577 \\
3.5	0.821577 \\
3.6	0.825726 \\
3.7	0.825726 \\
3.8	0.825726 \\
3.9	0.825726 \\
4	0.825726 \\
4.1	0.825726 \\
4.2	0.825726 \\
4.3	0.825726 \\
4.4	0.825726 \\
4.5	0.825726 \\
4.6	0.825726 \\
4.7	0.825726 \\
4.8	0.825726 \\
4.9	0.825726 \\
5	0.825726 \\
5.1	0.825726 \\
5.2	0.825726 \\
5.3	0.825726 \\
5.4	0.825726 \\
5.5	0.825726 \\
5.6	0.825726 \\
5.7	0.825726 \\
5.8	0.825726 \\
5.9	0.825726 \\
6	0.825726 \\
6.1	0.825726 \\
6.19999999999999	0.825726 \\
6.29999999999999	0.825726 \\
6.39999999999999	0.829876 \\
6.49999999999999	0.829876 \\
6.59999999999999	0.829876 \\
6.69999999999999	0.829876 \\
6.79999999999999	0.829876 \\
6.89999999999999	0.829876 \\
6.99999999999999	0.829876 \\
7.09999999999999	0.829876 \\
7.19999999999999	0.829876 \\
7.29999999999999	0.829876 \\
7.39999999999999	0.829876 \\
7.49999999999999	0.829876 \\
7.59999999999999	0.829876 \\
7.69999999999999	0.829876 \\
7.79999999999999	0.829876 \\
7.89999999999999	0.829876 \\
7.99999999999999	0.829876 \\
8.09999999999999	0.829876 \\
8.19999999999999	0.829876 \\
8.29999999999999	0.829876 \\
8.39999999999999	0.829876 \\
8.49999999999999	0.829876 \\
8.59999999999999	0.829876 \\
8.69999999999999	0.829876 \\
8.79999999999999	0.829876 \\
8.89999999999999	0.829876 \\
8.99999999999999	0.829876 \\
9.09999999999999	0.829876 \\
9.19999999999999	0.829876 \\
9.29999999999998	0.829876 \\
9.39999999999998	0.829876 \\
9.49999999999998	0.829876 \\
9.59999999999998	0.829876 \\
9.69999999999998	0.829876 \\
9.79999999999998	0.829876 \\
9.89999999999998	0.829876 \\
9.99999999999998	0.829876 \\
};
\addplot [line width=0.4800000000000001pt, color1, forget plot, mark=x]
table [row sep=\\]{%
1	0.481328 \\
};
\addplot [line width=0.4800000000000001pt, color1, dashed, forget plot]
table [row sep=\\]{%
1	0.481328 \\
1.1	0.60166 \\
1.2	0.655602 \\
1.3	0.692946 \\
1.4	0.717842 \\
1.5	0.738589 \\
1.6	0.738589 \\
1.7	0.755187 \\
1.8	0.767635 \\
1.9	0.792531 \\
2	0.79668 \\
2.1	0.79668 \\
2.2	0.80083 \\
2.3	0.813278 \\
2.4	0.813278 \\
2.5	0.813278 \\
2.6	0.817427 \\
2.7	0.817427 \\
2.8	0.817427 \\
2.9	0.817427 \\
3	0.817427 \\
3.1	0.821577 \\
3.2	0.821577 \\
3.3	0.825726 \\
3.4	0.825726 \\
3.5	0.825726 \\
3.6	0.829876 \\
3.7	0.829876 \\
3.8	0.829876 \\
3.9	0.834025 \\
4	0.834025 \\
4.1	0.834025 \\
4.2	0.834025 \\
4.3	0.834025 \\
4.4	0.838174 \\
4.5	0.838174 \\
4.6	0.838174 \\
4.7	0.838174 \\
4.8	0.838174 \\
4.9	0.838174 \\
5	0.838174 \\
5.1	0.838174 \\
5.2	0.838174 \\
5.3	0.838174 \\
5.4	0.838174 \\
5.5	0.838174 \\
5.6	0.838174 \\
5.7	0.838174 \\
5.8	0.838174 \\
5.9	0.838174 \\
6	0.838174 \\
6.1	0.838174 \\
6.19999999999999	0.838174 \\
6.29999999999999	0.838174 \\
6.39999999999999	0.838174 \\
6.49999999999999	0.838174 \\
6.59999999999999	0.838174 \\
6.69999999999999	0.838174 \\
6.79999999999999	0.838174 \\
6.89999999999999	0.838174 \\
6.99999999999999	0.838174 \\
7.09999999999999	0.842324 \\
7.19999999999999	0.842324 \\
7.29999999999999	0.842324 \\
7.39999999999999	0.842324 \\
7.49999999999999	0.842324 \\
7.59999999999999	0.842324 \\
7.69999999999999	0.842324 \\
7.79999999999999	0.842324 \\
7.89999999999999	0.842324 \\
7.99999999999999	0.842324 \\
8.09999999999999	0.842324 \\
8.19999999999999	0.842324 \\
8.29999999999999	0.842324 \\
8.39999999999999	0.842324 \\
8.49999999999999	0.842324 \\
8.59999999999999	0.842324 \\
8.69999999999999	0.842324 \\
8.79999999999999	0.842324 \\
8.89999999999999	0.842324 \\
8.99999999999999	0.842324 \\
9.09999999999999	0.842324 \\
9.19999999999999	0.842324 \\
9.29999999999998	0.842324 \\
9.39999999999998	0.842324 \\
9.49999999999998	0.842324 \\
9.59999999999998	0.842324 \\
9.69999999999998	0.842324 \\
9.79999999999998	0.842324 \\
9.89999999999998	0.842324 \\
9.99999999999998	0.842324 \\
};
\path [draw=white!80.0!black, fill opacity=0] (axis cs:0,0.3)
--(axis cs:0,1);

\path [draw=white!80.0!black, fill opacity=0] (axis cs:1,0.3)
--(axis cs:1,1);

\path [draw=white!80.0!black, fill opacity=0] (axis cs:1,0)
--(axis cs:10,0);

\path [draw=white!80.0!black, fill opacity=0] (axis cs:1,1)
--(axis cs:10,1);

\end{axis}

\end{tikzpicture}
        \caption{\bracket{1000}{tilim}}
        }
    \end{subfigure}%
  \end{center}
  \smallskip
  
  \caption{Performance profiles of \coefdiving and the \conflictdiving modification \confdivinglikecoef for four hierarchical groups of increasingly hard affected instances.}
  \label{fig:coef_vs_conf_version_performance}
\end{figure}
 
In our experiment \confdivinglikecoef was superior to \coefdiving and led to a performance improvement of $\percent{2.3}$ on the set of affected instances.
At the same time, the tree size could be reduced by $\percent{3.6}$ and four more instances could be solved.
Similar to the actual version of \conflictdiving evaluated at the beginning of this section, \conflictdiving with \confdivinglikecoef settings generated
more than twice as many conflict constraints and found $\percent{8}$ more improving solutions than \coefficientdiving.
The performance profiles in Figure~\ref{fig:coef_vs_conf_version_performance} show that \conflictdiving with the same rounding function like \coefficientdiving
performs better the harder the instance, \ie whenever many conflict constraints and therefore reliable conflict locks can be expected.
It is not surprising that for a time factor of $1.0$ \coefdiving solves more instances best when easy instances are considered, too,
since the number of conflict constraints can be expected to be small and, thus, the information gained by conflict locks might not be reliable enough.
Even though the profiles are close to each other, they do not cross if both settings need at least $10$, $100$, or $1000$ seconds.
This result indicates that it is the inclusion of conflict locks in general that helps to guide diving better than pure variable locks.

For the second control experiment we varied the weights $\confweight \in \{0,0.25,0.5,0.75,1\}$ in \confdiving in order to quantify the importance of conflict locks.
The results indicate that the reduced length of diving paths
is independent of the weighting of variable and conflict locks.
For every choice of $\confweight$ the average length of diving paths was reduced by $\percent{18}$ to $\percent{34}$.
Hence, the length of the diving paths seems to mainly depend on the rounding function $\roundfunc{C}$,
which always prefers rounding into the ``risky'' direction, \ie the direction with more locks.
We also confirmed that \confdiving is superior to \coefdiving on the affected instances for every evaluated $\confweight$.

Moreover, we observed that on instances with zero objective function a smaller conflict weight~$\confweight$ yields superior performance,
while for instances with nonzero objective function a larger weight~$\confweight$ for conflict information seems to be the better choice.
This observation may be related to how and when the decisions of \conflictdiving align with \scip'{s} default branching rule and when they complement it.
To make this clear, further background information is necessary.
For branching variable selection, \scip combines reliability pseudocost-branching~\citep{achterberg2005branching}, which estimates dual bound improvement, and hybrid branching~\citep{achterberg2009hybrid}, which includes conflict information via VSIDS~\citep{moskewicz2001chaff}.
Both VSIDS and conflict locks approximate the set of variables that frequently appear in conflict constraints.
On problems with nonzero objective function typically the first, objective-based score dominates \scip's branching decisions, while on pure feasibility problems the latter, conflict-based score has more impact.
This explains that on feasibility problems a larger $\confweight$ might align the decisions of \conflictdiving unfavorably with the overall tree search and create redundancy; on the majority of problems with nonzero objective, where conflict information is underweighted in the branching rule, a larger $\confweight$ enables \conflictdiving to contribute more complementary information to the search.
This suggests a further refinement of \conflictdiving by adapting the weight~$\confweight$ dynamically to the problem at hand.

\paragraph{Impact of Primal Solutions and Generated Conflicts.}

Finally, we analyzed the importance of found solutions and generated conflict constraints by comparing \confdiving to two artificially modified versions of \conflictdiving:
one that does not apply conflict analysis (\confdivingnoconfs) and one that discards feasible solutions found by \conflictdiving (\confdivingnosols).
Aggregated results on the set of affected instances are shown in Table~\ref{tab:conflictdiving_no_sols_no_confs}.

\begin{table*}%
\begin{centering}
\footnotesize
\setlength{\tabcolsep}{.75pt}
\caption{\confdiving with and without adding feasible solutions and generated conflict constraints, respectively, on the set of affected instances.}\label{tab:conflictdiving_no_sols_no_confs}
\begin{tabularx}{\textwidth}{L *{6}{R}}
    \toprule
    & \# & \solved & \time & \nodes & \timeQ & \nodesQ \\
    \midrule
    \mbox{\nolockdiving}      & 861 & 822 & 185.62 & 4235 & 1.000 & 1.000  \\
    \mbox{\coefdiving}        & 861 & 831 & 185.25 & 4290 & 0.998 & 1.013  \\
    \mbox{\confdiving}        & 861 & 838 & 176.57 & 3999 & 0.951 & 0.944  \\
    \mbox{\confdivingnoconfs} & 861 & 825 & 181.99 & 4131 & 0.980 & 0.976  \\
    \mbox{\confdivingnosols}  & 861 & 838 & 177.64 & 4006 & 0.957 & 0.946  \\
    \bottomrule\\
\end{tabularx}
\end{centering}
\end{table*}
 
Both variants \confdivingnoconfs and \confdivingnosols are superior to \coefdiving \wrt solving time and number of nodes.
Whereas \confdivingnosols solved the same amount of instances as \confdiving,
disabling conflict analysis led to solving $13$ instances less than \confdiving.
Hence, disabling the addition of feasible solutions found by \conflictdiving only has a marginal impact compared to standard \conflictdiving.
By contrast, disabling conflict analysis during \conflictdiving
leads to a slowdown of $\percent{3}$ compared to \confdiving. %
Thus, the conflict constraints generated during \conflictdiving seem to be the main driver of the heuristic.
This also aligns with our observations regarding the \emph{dual integral}~\citep{Berthold2013}.
The dual integral measures the progress of the dual bound towards the optimal objective value.
This measure increased by $\percent{3.1}$ on average when disabling conflict analysis, which indicates that the increased number of conflict constraints helps to accelerate convergence of the proof of optimality.

\section{Conclusion}
\label{sec:conclusionandoutlook}

In this paper, we presented two new ways how conflict analysis
and primal heuristics can be combined to improve the performance of a MIP solver.
We presented a primal heuristic, called \farkasdiving, that simultaneously aims to construct valid Farkas proofs and feasible solutions.
By design, \farkasdiving is more expensive than previous diving heuristics in \scip.
Therefore, the heuristic is called conservatively, whereby the decision to keep the heuristic enabled for the remainder of the search
is based on its success during the root node.
On the set of instances where \farkasdiving was executed after the root node, it proved to be very successful in generating improving solutions
and conflict constraints.
Regarding both metrics, \farkasdiving outperforms the virtual best of all other diving heuristics on this set of instances.
Moreover, the overall solving time could be improved by $\percent{5.4}$
and the tree size could be reduced by $\percent{14.8}$ on the set of affected instances.

Furthermore, we applied the concept of variable locks to conflict constraints and used this additional information
to guide the search of a second diving heuristic, called \conflictdiving.
In addition, \conflictdiving pursues an aggressive fail fast strategy to prevent fruitless consumption of running time.
This new diving heuristic is an extension and modification of the well-known \coefficientdiving heuristic.
Our computational results indicate that this mix of variable and conflict locks
in combination with an aggressive fail fast strategy outperforms \coefficientdiving on our complete test set.
Conflict diving reduces the overall solving time and tree size on affected instances by $\percent{4.8}$ and $\percent{5.5}$, respectively.

The two ways of combining primal heuristics and conflict analysis presented in this paper have highlighted
that primal and dual solving techniques within a general-purpose MIP solver are not independent and that they can not only interact randomly and create performance variability,
but be combined beneficially in a targeted manner.

\ifreport
\paragraph{Acknowledgments}
\else
\ACKNOWLEDGMENT{%
\fi
Many thanks to Timo Berthold for valuable discussions and Gregor Hendel for his continuous effort in developing IPET.
The work for this article has been conducted within the Research Campus Modal
funded by the German Federal Ministry of Education and Research (grant number 05M14ZAM)
and has received funding from the European Union's Horizon 2020 research and innovation programme under grant agreement No 773897 (plan4res). The content of this paper only reflects the author's views. The European Commission / Innovation and Networks Executive Agency is not responsible for any use that may be made of the information it contains.
\ifreport
\else
}
\fi

\ifreport
    \bibliographystyle{informs2014}
\else
    \bibliographystyle{informs2014}
\fi

\bibliography{Bibliography}

\end{document}